

\input amstex
\documentstyle {amsppt}
\magnification \magstep1
\openup3\jot
\NoBlackBoxes
\pageno=1
\hsize 6 truein
 


 
\predefine\enddimost{\enddemo}
\redefine\enddemo {\penalty 5000\hskip15pt plus1pt minus5pt\penalty1000
		   \qed\enddimost}

\catcode`@=11
\font\s@=cmss10
\predefine\reall{\Re}
\def\Id{{\hbox{\s@ 1\kern-.8mm l}}} 
\redefine\Re{\Bbb R}

\def\Z{\Bbb Z}
\def\To{\Bbb T}


\define\RT {\widetilde{\Cal S}}
\define\RB {\Cal S_\delta}
\define\TB {\widetilde{\Cal S}_\delta}


\define\Q{\Cal Q}
\define\C{\Cal C}
\define\R{\Cal S}
\define\U{\Cal U}  
\define\M{\Cal M}

\define\Lon{\biggr\|_{\lower 2pt\hbox{$\scriptstyle 1$}}}
\define\lon{\bigr\|_{\lower 2pt\hbox{$\scriptstyle 1$}}}
 




\topmatter
 
\title
ERGODICITY IN\\
HAMILTONIAN SYSTEMS.
\endtitle

\author
Carlangelo Liverani,\ Maciej Wojtkowski
\endauthor
\address
Liverani Carlangelo,
Mathematics Department,
University of Rome II, Tor Vergata,
Rome, Italy.
\endaddress
\email
liverani@vaxtvm.infn.it
\endemail
\address
Institute for Mathematical Sciences,
SUNY at Stony Brook,
Stony Brook, NY 11794, USA.
\endaddress
\email
liverani@math.sunysb.edu
\endemail
\address
Maciej Wojtkowki,
Mathematics Department,
University of Arizona,
Tucson, AZ 85721, USA.
\endaddress
\email
maciejw@math.arizona.edu
\endemail
\date
9-20-92
\enddate

\abstract
We discuss the Sinai method of proving ergodicity of a discontinuous
Hamiltonian system with (non-uniform) hyperbolic behavior.
\endabstract

\thanks
{\bf We would like to thank N. Chernov, L. Chierchia, V. Donnay,
A. Katok, N. Sim\'anyi,
D. Sz\'asz and L.-S. Young for helpful and enlightening
discussions.
The first author wishes to thank the the Mathematics Department of the
University of Tucson and the Center for Applied
Mathematics at Cornell University, in particular its director J.
Guckenheimer, where he was visiting during part of this work, he
also acknowledge the partial support received by CNR, grant n. 203.01.52,
and by the GNFM. The second author gratefully acknowledges the hospitality
of Forschungsinstitut f\"ur Mathematik at ETH Z\"urich, where the first draft
of this paper was written. He also acknowledges the partial support from
NSF Grant DMS-9017993.}
\endthanks
 
\endtopmatter
 

\document

\vskip.5cm
\centerline{\bf CONTENT}
\vskip.5cm
\roster
\item"0."  Introduction\dotfill p. \ 3
\item"1."  A Model Problem\dotfill p. \ 4
\item"2."  The Sinai Method\dotfill p. 10
\item"3."  Proof of Sinai Theorem (for piecewise linear maps of
	   the two torus)\dotfill p. 14
\item"4."  Sectors in linear symplectic space\dotfill p. 18
\item"5."  The space of Lagrangian subspaces contained in a sector\dotfill
	   p. 23
\item"6."  Unbounded Sequences of Monotones Maps\dotfill p. 28
\item"7."  Properties of the system and the formulation of the results
	   \dotfill p. 35
\item"8."  Construction of the neighborhood and coordinate system
	   \dotfill p. 45
\item"9." Unstable manifolds in the neighborhood $\U$\dotfill p. 48
\item"10." Local ergodicity in the smooth case\dotfill p. 54
\item"11." Local ergodicity in the discontinuous case\dotfill p. 56
\item"12." Proof of Sinai Theorem\dotfill p. 60
\item"13." `Tail Bound'\dotfill p. 65
\item"14." Applications\dotfill p. 69
\item" "   References\dotfill p. 79
\endroster



\vskip 1cm
\subheading{SYMBOLS USED IN THE PAPER}
\vskip 2cm
$$
\aligned
\qquad\qquad&\\
\alpha\qquad&\text{amount of long leaves in a connecting square}\\
\Cal B(p;\;r)\qquad&\text{Ball of radius $r$ and center $p$}\\
c\qquad& \text{amount of overlap in neighboring squares}\\
\Cal C\qquad&\text{sectors}\\
d\qquad&\text{distance}\\
k(c)\qquad& \text{maximal number of overlapping squares}\\
L\qquad&\text{linear map}\\
\M \qquad&\text{Symplectic manifold}\\
\M^\pm\qquad&\text{Symplectic boxes}\\
\mu\qquad&\text{invariant measure}\\
\omega\qquad&\text{symplectic form}\\
\Cal Q\qquad&\text{quadratic form defining a sector}\\
R\qquad&\text{rectangles}\\
\Cal G\qquad&\text{collection of rectangles}\\
\Cal S^\pm\qquad&\text{singularity sets}\\
T\qquad&\text{map}\\
U\qquad&\text{big neighborhood in the smooth case}\\
\U(x)\qquad&\text{neighborhood of $x$}\\
V\qquad&\text{side of a sector}\\
\Cal W\qquad&\text{linear symplectic space}\\
W\qquad&\text{stable and unstable manifolds}\\
\endaligned
$$
\vskip 1cm
\centerline{In the Figures}
\centerline {\bf the stable direction is vertical}
\centerline{\bf the unstable direction is horizontal}

          \par\newpage   



\vskip.7cm
\subhead \S 0. INTRODUCTION \endsubhead
\vskip.4cm

The notion of ergodicity was introduced by Boltzman as a property
satisfied by a Hamiltonian flow on its energy manifold.  The emergence
of the KAM (Kolmogorov-Arnold-Moser) theory of quasiperiodic
motions made it clear that very few Hamiltonian systems are actually
ergodic. Moreover, those systems which seem to be ergodic do not lend
themselves easily to rigorous methods.

Ergodicity is a rather weak property in the hierarchy of stochastic
behavior of a dynamical system. The study of strong properties
(mixing, K-property and Bernoulliness) in smooth dynamical systems
began from the geodesic flows on surfaces of negative curvature. In
particular, Hopf \cite {H} invented a method of proving ergodicity,
using horocycles, which turned out to be so versatile that it endured a
lot of generalizations. It was developed by Anosov and Sinai \cite
{AS} and applied to Anosov systems with a smooth invariant measure.
With the advances of the theory of Kolmogorov - Sinai entropy the
Hopf method turned out to be also a basis for proving the K-property
of Anosov systems. 

The key role in this approach is played by the 
hyperbolic behavior in a dynamical system. By the hyperbolic 
behavior we mean the
property of exponential divergence of nearby orbits. In the strongest
form it is present in Anosov systems and Smale systems. It leads there to
a rigid topological behavior.  In weaker forms it seems to be a common
phenomenon.

In his pioneering work on billiard systems Sinai \cite{S} showed that
already weak hyperbolic properties are sufficient to establish the 
strong mixing properties. Even the discontinuity of the system can
be accommodated.

The Multiplicative Ergodic Theorem of Oseledets \cite{O} makes
Lyapunov exponents a natural tool to describe the hyperbolic 
behavior of a dynamical system with a smooth invariant measure.

Pesin \cite {P} made the nonvanishing of Lyapunov exponents the
starting point for the study of hyperbolic behavior. He showed that, if
a diffeomorphism preserving a smooth measure has only nonvanishing
Lyapunov exponents, then it has at most countably many ergodic
components and
(roughly speaking) on each component it has the Bernoulli property.

Pesin's work raised the question of sufficient conditions for 
ergodicity or, more modestly, for the openness (modulo sets of measure
zero) of the ergodic components.

In his work, spanning two decades, on the system of colliding balls
(gas of hard balls) Sinai developed a method of proving (local)
ergodicity in discontinuous systems with nonuniform hyperbolic
behavior.  We will refer to it as the Sinai method. It was improved by
Sinai and Chernov
\cite {CS} and by A.Kr\' amli, N.Sim\' anyi and D.Sz\' asz \cite
{KSS}. In both papers the discussion is confined to the realm of 
semidispersing billiards.

The purpose of the present paper is to recover the Sinai method as a
part of the theory of hyperbolic dynamical systems. In the process we
have simplified some of the aspects of the method, and we have revealed
its logical structure and limitations.

We rely on two developments. The first is the work of Katok and Strelcyn
\cite {KS} in which they generalized Pesin Theory to discontinuous 
systems. The other is the development of criteria for nonvanishing of
Lyapunov exponents in Hamiltonian systems in papers \cite {W1},
\cite {W2} and \cite{W3}.
In the language of these criteria 
Burns and Gerber \cite {BG} found 
a sufficient condition for (local) ergodicity in the smooth case of
lowest dimension (3 for flows preserving a smooth measure). 
It was later generalized by Katok \cite {K1} to arbitrary dimension.
As a byproduct of our general approach, which includes discontinuous 
systems, we obtain a similar theorem
(Main Theorem in the smooth case) and a new proof.

Let us give some advice to the reader on how to use our paper.  The
first three Sections demonstrate what the Sinai method is
and how it works. The discussion is conducted in the simplest possible
environment of a linear discontinuous system on the two dimensional
torus. It is reasonable to stop here, especially if the reader is only
interested in two dimensional uniformly hyperbolic systems. But we do
not recommend trying to read the heart of the paper without going
through the first three Sections.

In Sections 4,5 and 6 we develop the linear symplectic language in
which we formulate our results. We suggest that the reader skips
these sections and goes straight to Section 7 where we formulate the
multitude of hypotheses and the two Main Theorems on local ergodicity,
one for smooth systems and the other (much harder) for discontinuous
systems. The reading of Section 7, and the following Sections, will require
numerous trips back to Sections 4-6 for the necessary definitions and
theorems.

If the reader does not care about the discontinuous case, she needs to
read only Sections 8, 9 and 10 with significant leaps (since
everything is simpler in the smooth case). Sections 11 and 12 contain
almost the whole proof of the Main Theorem in the discontinuous case
(it also relies on the results of Sections 8-10). The remaining part of the
proof is contained in Section 13. It stands out by the level of
technical complications.

Section 14 contains some classes of examples where all the hard work
can be put to use, and one class where it cannot. The interest in this
last example comes from the fact that it is multidimensional and all
the Lyapunov exponents are different from zero. Unfortunately, it does
not satisfy an important property (proper alignment of singularity
sets). It points towards the need for a more flexible scheme.

\vskip.7cm
\subhead \S 1. A MODEL PROBLEM \endsubhead
\vskip.4cm
 
We will discuss here a very simple model problem in which the
important features of the Sinai's method are not obscured by technical
details. Our discussion will be very careful so that in the future when
the technical details will cloud the horizon we will be able to refer the
reader to these basic clarifications.
 
We consider a family of linear maps  of the plane
defined by
$$
\aligned
x_1' &= x_1 + a x_2 \\
x_2' &= x_2 ,
\endaligned
$$
where $a$ is a real parameter.  We use these linear maps to define
(discontinuous) maps of the torus by restricting the formulas to the
strip $\{ 0 \leq x_2 \leq 1 \}$ and further taking them modulo 1. In
this way we define a mapping $T_1$ of the torus $\To ^2 =
\Re^2/\Z^2$ which is discontinuous on the circle $\{ x_2 \in \Z
\}$ (except when $a$ is equal to an integer) and preserves the
Lebesgue measure $\mu$.
 
Similarly we define another family of maps depending on the same
parameter $a$ by restricting the formulas
$$
\aligned
x_1' &= x_1 \\
x_2' &= a x_1 + x_2
\endaligned
$$
to the strip $\{ 0 \leq x_1 \leq 1 \}$ and then taking them modulo 1. Thus
for each
$a$ we get a mapping $T_2$ of the torus which is discontinuous on the circle
$\{   x_1 \in \Z \}$ (except when $a$ is equal to an integer) and
preserves  the Lebesgue measure $\mu$.
 
Finally we introduce the composition of these maps $T =
T_2 T_1$ which depends on one real parameter $a$.
An alternative way of describing the map
$T$ is by introducing two fundamental domains for the
torus
$\M^+ = \{ 0 \leq x_1 + a x_2 \leq 1,\, 0 \leq x_2 \leq 1 \}$
and $\M^- = \{ 0 \leq x_1 \leq 1, \,0 \leq - a x_1 + x_2 \leq 1,\, \}$
(see Fig.1).
 
\topinsert
\vskip 3in
\hsize=4.5in
\raggedright
\noindent{\bf Figure 1} The map.
\endinsert
 
The linear map defined by the matrix
$$
\left(\matrix 1 & a \\ a & 1 + a^2\endmatrix\right)=
\left(\matrix 1 & 0 \\ a& 1 \endmatrix\right)
\left(\matrix 1 & a \\ 0 & 1 \endmatrix\right)
$$
takes $\M^+$ onto $\M^-$ thus defining a map of the torus which is
discontinuous at most on the boundary of $\M^+$ and preserves the
Lebesgue measure. This is our map $T$.
 
Let $\Cal S^{\pm}= \partial\M^{\pm}$ be the boundary of $\M^{\pm}$.
Except for integer values of $a$ the mapping $T$ is discontinuous on
$\Cal S^+$ and its inverse $T^{-1}$ is discontinuous on $\Cal S^-$. Let us stress that
the map $T$ is well defined in the closed domain $\M^+$ but two
different points on the boundary $\Cal S^+$ which correspond to the {\it same}
point on the torus will be mapped onto two different points on the
boundary $\Cal S^-$ which correspond to two {\it different} points on the
torus (except for the corner). We adopt the convention that the image
under $T$ of a point from $\Cal S^+$ is the pair of image points in
$\Cal S^-$. With this convention we can apply $T$ or any of its powers to any
subset in the torus.
 
For integer  values of $a \neq 0$ we have a hyperbolic
algebraic automorphism of the torus, a prime example of an Anosov
system. It is thus a Bernoulli system and has a nice
Markov partition \cite{AW}.
We restrict ourselves to the study of ergodicity and we repeat
the proof of ergodicity by the Hopf method, since the Sinai method
is built upon it.
 
Let $f:\To ^2 \to \Re$ be a continuous function. We want to prove that for
 almost
every $x\in \To ^2$ the time averages
$$
\frac{f(x) + f(T x) + \dots + f(T^{n-1}x)}{n}
$$
converge as $n \rightarrow +\infty$ to the average value of $f$, i.e.,
$\int f d\mu$. Once this is established  one can obtain the same
property for all integrable functions by an approximation argument.
From Birkhoff Ergodic Theorem (BET) we know that the time averages
converge almost everywhere to a function
$f^+ \in L^1(\To ^2,\,\mu)$ which is invariant on the orbits of $T$, i.e.,
$f^+\circ T = f^+$, and has the same average value as $f$, i.e.,
$\int f^+ d\mu = \int f d\mu$.
Further applying BET to $f$ and $T^{-1}$ we obtain that the time averages
in the past
$$
\frac{f(x) + f(T^{-1} x) + \dots + f(T^{-n+1}x)}{n}
$$
converge almost everywhere as $n \rightarrow +\infty$ to
$f^- \in L^1(\To ^2,\,\mu)$ for which
$f^-\circ T = f^-$ and $\int f^- d\mu = \int f d\mu$.
 
It is the usual magic of the ergodic theory which forces the functions $f^+$
and
$f^-$ to coincide almost everywhere.
(Let us recall the argument: let
$$
\Cal A_+=\{x\in\To^2\;|\;f_+(x)> f_-(x)\};
$$
by definition
$\Cal A_+$
is an invariant set, hence
$$
\int_{\Cal A_+}\left [f_+(x)-f_-(x)\right ]d\mu(x)=
\int_{\Cal A_+}f(x)d\mu(x)-\int_{\Cal A_+}f(x)d\mu(x)=0
$$
which implies $\mu(\Cal A_+)=0$ and $f_+\leq f_-$ $\mu$-almost
everywhere. The same argument, this time applied to the set $\Cal
A_-=\{x\in\To^2\;|\;f_-(x)> f_+(x)\}$, implies the converse inequality.)
 
For $a \neq  0$  the matrix
$$
\left(\matrix 1 & a \\ a & 1 + a^2\endmatrix\right)
$$
is a hyperbolic matrix with eigenvalues $\lambda = \lambda(a) > 1$ and
$\frac1\lambda < 1$.
For $x\in \To^2$ let us denote by $W^u(x)$ $(W^s(x))$
the line in $\To^2$ passing through $x$ and having the direction of the
unstable
eigenvector (the stable eigenvector), i.e., the eigenvector with eigenvalue
$\lambda$
($\frac1\lambda$). We call $W^u(x)$ $(W^s(x))$ the unstable (stable) leaf of
$x$. The leaves of $x$ have the following property. If $y\in W^{u}(x)$
$(y\in W^s(x))$ then the distance
$$
d(T^ny,\,T^nx) = \lambda^{-|n|} d(y,\,x) \rightarrow 0 \ \
\text{as} \ \ n \rightarrow -\infty (+\infty).
$$
Hence for $y,z \in W^{u(s)}(x)$
$$
|f(T^ny) - f(T^nz) | \rightarrow 0 \ \ \text{as} \ \
n \rightarrow -\infty (+\infty).
$$
It follows that for $y,z \in W^{u(s)}(x)$ either $f^{\pm}(y)$  and
$f^{\pm}(z)$
are both defined and equal or they are both undefined. Lifting the
functions $f^+$
and $f^-$ to $\Re^2$ and using the directions of the eigenvalues as
coordinate
directions we can say that $f^+$ is a function of one coordinate alone
and
$f^-$ is a function of only the other coordinate.
Since the two functions coincide almost everywhere they must be constant.
\par
Let us examine what can be saved of this argument when $a$ is not an integer.
In such a case, we still have the stable and unstable directions but
 a line parallel to, say, the unstable
direction is cut by $\Cal S^-$ into pieces and if $y$ and
$z$ belong to two different pieces
the distance $d(T^ny,\,T^nz)$ does not decrease to zero as
$n \rightarrow -\infty$. Since this last property is of crucial importance in
the Hopf method,
the unstable (and stable) leaves have to be much shorter than before.
 Here is how we construct them. 
For simplicity of notation  we will formulate everything for the
unstable leaves alone.
\par
We proceed inductively. Thus, for $x \in int \M^-$, we define $W_1^u(x)$
as the open segment of the line through $x$ with the direction of the
unstable
eigenvector which contains $x$ and has both endpoints on $\Cal S^-$.
The preimage
$T^{-1}W_1^u(x) $ is by a factor of $\lambda$ shorter than $\Cal
W_1^u(x)$ and, in general, is cut into two or three pieces by $\Cal S^-$.
We pick the piece which contains
$T^{-1}x $ and take its image under $T$; this is our second approximate
unstable leaf  $W_2^u(x)$, i.e.,
$$
W_2^u(x) = T\left(T^{-1}W_1^u(x) \cap
W_1^u(T^{-1}(x))\right).
$$
Unless $T^{-1}x \in \Cal S^-$ the second approximate unstable leaf
$W_2^u(x)$
is again an open segment containing $x$
with endpoints on $\Cal S^- \cup T \Cal S^-$ and naturally $W_2^u(x)
\subset W_1^u(x)$.
Given $W_n^u(x)$, $n=1,2,\dots ,$ we define the $n+1$ approximate unstable
leaf of $x$ $W_{n+1}^u(x)$ by
$$
W_{n+1}^u(x) = T^n\left(T^{-n}W_n^u(x) \cap
W_1^u(T^{-n}(x))\right).
$$
 
If $x \notin \bigcup_{i=0}^{+\infty} T^i\Cal S^-$ then
this inductive procedure will yield a nested sequence of open segments
containing $x$
$$
W_1^u(x) \supset W_2^u(x) \supset \dots
$$
with endpoints on
$$
\bigcup_{i=0}^{+\infty} T^i \Cal S^-\ .
$$
We can also describe this construction in the following way.
First we consider a fairly
long segment $W_1^u(x)$. Then we look at $T \Cal S^-$, if it does not intersect
 $W_1^u(x)$ then we do not change it, if it splits $W_1^u(x)$ into
several
segments, then we keep the segment which contains $x$. We repeat it
with $T^2 \Cal S^-$ and
further images of $\Cal S^-$, so that the segment may be cut shorter infinitely
many times. The property  $x \notin \bigcup_{i=0}^{+\infty} T^i\Cal S^-$
ensures that $x$ stays always strictly inside the segment.
It is quite remarkable that, for almost every $x$, this inductive process
shortens
the segment only finitely many times. More precisely we have
\proclaim{Proposition 1.1}
For almost all $x \in \M^- \setminus \bigcup_{i=0}^{+\infty} T^i\Cal S^-$
the sequence of
approximate unstable leaves of $x$ stabilizes, i.e., there is a
natural $N = N(x)$ such that
$$
\bigcap_{i=1}^{+\infty} W_i^u(x) = \bigcap_{i=1}^N W_i^u(x).
$$
\endproclaim
\demo{Proof}
For $t >0$, let
$$
X_t = \{ x \in \M^- \;|\; d(x,\,\Cal S^-) \leq t\}
$$
where $d(\cdot,\,\cdot)$ is the distance of a point form a set.
Because $\Cal S^-$ is a finite union of segments we have
$$
\mu \left(X_t \right) \leq \text{const} \ t.
$$
Choosing $t_n = \frac 1{n^2}$ we get
$$
\sum_{n=1}^{+\infty} \mu\left( X_{t_n} \right) < +\infty,
$$
hence also
$$
\sum_{n=1}^{+\infty}  \mu\left(T^n X_{t_n}\right)  < +\infty.
$$
It follows by the Borel-Cantelli Lemma that almost every $x$ belongs to
only finitely many of the sets
$$
T X_{t_1},T^2 X_{t_2},\dots  ,
$$
which means that  except for finitely many values of $n$
$$
d(T^{-n}x,\,\Cal S^-) > \frac 1{n^2}.
$$
Choosing  $c(x) > 0$ sufficiently small we can take care of the
finite number of exceptional
values of $n$ so that
$$
d(T^{-n}x,\,\Cal S^-) > \frac {c(x)}{n^2}
$$
for each $n = 1,2,\dots \ .$
Each time $W_{n+1}^u(x)$ is shorter than $W_{n}^u(x)$ we must have
$$
d(T^{-n} x,\,\Cal S^-)  <  \frac{\text{length}\left(\Cal W_{n}^u(x)\right)}
{\lambda^{n}}.
$$
But then 
$$
\frac {c(x)}{n^2} < \frac{\text{length}\left(\Cal
 W_{n}^u(x)\right)}{\lambda^{n}}
\leq \frac{\text{length}\left(W_1^u(x)\right)}{\lambda^{n}},
$$
which can hold for at most finitely many values of $n$.
\enddemo
 
We define the unstable leaf only for points $x$ in the set of full
measure described
in Proposition 1.1,   by taking the intersection
$$
W^u(x) = \bigcap_{i=1}^{+\infty} W_i^u(x).
$$
In view of Proposition 1.1, for each $W^u(x)$, there are
natural numbers $n_l(x)$ and $n_r(x)$ such that $T^{n_l(x)}W^u(x)$
has the left endpoint on $\Cal S^-$ and $T^{n_r(x)}W^u(x)$
has the right endpoint on $\Cal S^-$.
Most importantly we have the exponential contraction of $W^u(x)$, i.e., for
$y \in W^u(x)$ the distance
$$
d(T^{-n}y,\,T^{-n}x) = \frac{d( y,\, x)}{\lambda^n} \rightarrow 0 \ \
\text{as} \ \ n \rightarrow  +\infty .
$$
\par
Everything that we have done to construct the unstable leaves
can be repeated for the stable leaves and they have analogous
properties. Once we have the stable and unstable leaves we 
are ready to do the Hopf argument.

For any continuous function $f:\To ^2 \to \Re$ the forward ergodic average
$f^+$
is constant on the stable leaves and the backward ergodic average $f^-$ is
constant
on the unstable leaves. Let us call a point $x\in \To ^2$ $f$-typical, if
$f^+(x)$, $f^-(x)$, $W^u(x)$ and $W^s(x)$ are well defined
and $f^+(x) = f^-(x)$.
The set of $f$-typical points has full measure,
so a stable (or an unstable) leaf
contains a set of $f$-typical points of full arc-length, except for a family
of leaves of total measure zero. If $W^s(x)$ is not one of those exceptional
leaves, then the set
$$
C_1 = \bigcup_{\Sb y \in W^s(x)\\y\ \text{is}\ f-\text{typical} \endSb}
W^u(y)
$$
has positive measure and $f^- = f^+ = const$ on $C_1$. We can proceed by
adding
all the
stable leaves through $f$-typical points in $C_1$ to obtain $C_2$, etc., but
a priori there is no reason to expect that we will
be able to cover all of the torus in this way.
(Indeed one can imagine that there is a dividing line between two ergodic
components of our system and that all the stable and unstable leaves stop
short of crossing this line.)
That is where the Hopf method breaks down. It can only tell us that the
ergodic
components have positive measure and, therefore, that there are at most
countably many of them.
(To be more precise, we cannot really claim   
that $C_1$ belongs to one ergodic
component. To argue this we have to modify our argument by
taking a sequence of continuous functions dense in $L^1$ and considering
the set of
points which are $f$-typical for all the functions $f$ in the
sequence. This set, as the
intersection
of countably many sets of full measure, has full measure. We can then use it
in the definition of $C_1$ and claim that $f^- = f^+ = const$ on $C_1$ for
all the functions in our dense sequence. This implies that such $C_1$ does
belong to one ergodic component.
It follows easily that every invariant subset of positive measure
contains an ergodic component of positive measure. Hence all ergodic
components have positive measure.)

\vskip.7cm
\subhead \S 2. THE SINAI METHOD \endsubhead
\vskip.4cm
We have seen, in the previous section, that the Hopf method is not
sufficient to prove the ergodicity of a discontinuous map
because the stable and unstable leaves may be short. The
Sinai method amounts to establishing that most of the stable and
unstable leaves are, in a certain sense, sufficiently long. The first
(highly nontrivial) step in this method is to formulate precisely what is
meant by ``sufficiently long".  As before, we do it only for
 the unstable leaves; the changes necessary in the case of 
stable leaves are automatic.
 
Let $\U \subset \To^2$ be a (small) square with the sides parallel to unstable
and stable directions respectively (to make the geometry simpler let us
think that the unstable direction is horizontal and the stable
direction vertical). For any $0 < c <1$ we construct a sequence
$\Cal G_n(c), n=1,2,\dots ,$ of coverings of $\U$ in the following
way. Without loss of generality we can let 
$$\U =\{(u,v)\;|\;-b < u < b,\, -b < v < b\}.
$$
We consider the net
$\Cal N(n,c)$ defined by
$$
\Cal N(n,c)=
\{\frac cn (m,\,k)\in \U\;|\; m,k\in\Z\} .
$$
 Now the covering $\Cal G_n(c)$ is the collection of squares
having centers at points from
$\Cal N(n,c)$ and sides, of length $\frac 1n$, parallel to the sides of $\U$.
If $c< \frac 12$ then $\Cal G_n(c)$ is a
covering of $\U$ (otherwise $\Cal G_n(c)$ may cover only a
smaller square). The parameter $c$
will be chosen later to be very small, so that many squares in
$\Cal G_n(c)$ overlap. However, once $c$ is fixed, a point
in $\U$  may belong, at most, to
a fixed number, independent of $n=1,2,\dots$, of squares in $\Cal G_n(c)$;
we denote this number by
$k(c)$ (one can easily establish that $k(c) \leq (\frac 1{2c} + 1)^2$,
but we will not use any explicit estimate).

\topinsert
\vskip 3in
\hsize=4.5in
\raggedright
\noindent{\bf Figure 2} The covering.
\endinsert

We call two squares, in $\Cal G_n(c)$, immediate neighbors if the
distance between their centers is $\frac cn$.Two immediate neighbors
overlap on $1-c$ part of their areas.

 One can naturally define a column of squares and a row of squares as
special collections of squares in $\Cal G_n(c)$ (see Figure 2). For
example, a sequence $\{R_i\}_{i=1}^l$ of squares from $\Cal G_n(c)$ is
called a column of squares if, for every $i = 1,\dots,l-1$, $R_i$ and
$R_{i+1}$ are immediate neighbors, $R_{i+1}$ is above $R_i$, and there
is no square in $\Cal G_n(c)$ below $R_1$ or above $R_l$.

For each square $R \in \Cal G_n$ we introduce the stable, $\partial _s R$,
and unstable, $\partial _u R$, boundaries of $R$; $\partial _s R$ is
the union of the two boundary segments of $R$ which have
the stable (vertical) direction and $\partial _u R$ is the union of
the two
boundary segments of $R$ which have the unstable
(horizontal) direction.
Given a point $x \in R$, the unstable leaf $W^u(x)$ may intersect
both segments in $\partial _s R$ or it may be too short to reach one
of them (or both).  In the first case we say that $W^u(x)$ is
long in $R$, or that it is connecting in $R$ , in the second that
it is short in $R$ or that it is not connecting in  $R$.
 
\proclaim{Definition 2.1}
Given $\alpha,\, 0 < \alpha < 1,$ we call a square $R \in \Cal G_n(c)$
$\alpha$-connecting if
the measure of the set of points $x \in R$ whose unstable leaf
$W^u(x)$ is long in $R$ is at least $\alpha $ part of
the total area of $R$.
\endproclaim
 
Sinai formulates the property that most of unstable leaves are
sufficiently long in the following way.

\proclaim{Sinai Theorem 2.2}
There is $\alpha_0 <1$ such that for any $\alpha,\, 0 < \alpha \leq
\alpha_0$ and any $c,\, 0<c<1$,
$$
\lim_{n \rightarrow +\infty}
n \ \mu \left(\bigcup\{R \in \Cal G_n(c)\;|\; \,R
\text{ is not } \alpha\text{-connecting }\} \right) = 0.
$$
\endproclaim
In other words, the theorem says that if $\alpha$ is sufficiently small,
then the union of the squares in $\Cal G_n(c)$ which are not
$\alpha$-connecting 
has measure $\text{\it o}(\frac 1n)$.
\par
Before proving the Sinai Theorem let us show how it can
be used to get information about ergodic components.  Notice that
Definition 2.1 and the Sinai Theorem can be repeated for stable leaves.
 
\proclaim{Proposition 2.3}
The square $\U \subset \To^2$ (for which the Sinai Theorem holds
for both unstable leaves and stable leaves) belongs to one ergodic
component of $T$.
\endproclaim
 
In view of the arbitrariness of the square $\U$ to which we can apply
this Theorem we obtain immediately
\proclaim{Corollary 2.4}
The map $T$ is ergodic.
\endproclaim
 
\demo{Proof of Proposition 2.3}
Let us fix $\alpha$ sufficiently small so that the Sinai
Theorem holds for $\alpha$-connecting squares both in the unstable and
stable versions.  Next we fix $c$  smaller than $ \alpha$. As a
consequence two $\alpha$-connecting squares in $\Cal G_n(c)$, which
are immediate neighbors, contain in their intersection a
set of connecting leaves of positive measure.  The reason is that immediate
neighbors intersect over $1-c$ part of their areas and hence the
guaranteed $\alpha$ part of the square covered by connecting leaves
cannot fit into the remaining $c$ part of the square.
In the following we will not change the values of $\alpha$ or $c$
and, for simplicity, we will call an $\alpha$-connecting square simply
a connecting square. Thus a connecting square is $\alpha$-connecting
both with respect to stable and unstable leaves.
\par
Consider any continuous function $f$ on the torus.
We  call a point $y \in \Bbb T ^2$ $f$-typical if
the forward  time average $f^+$ and the backward
time average $f^-$ are well defined at $y$ and $f^+(y) = f^-(y)$.
The set of $f$-typical points has full measure.
We call a stable (unstable) leaf $f$-typical if its points, except for
a subset of zero arc-length, are $f$-typical. The union of leaves
which are not $f$-typical is a set of measure zero.

 For any connecting square $R$ let us define
$$
W^{u(s)}(R)= \{x\in R | W^{u(s)}(x)\  \text{is}\ f\text{-typical and long
in}\ R \}.
$$

Although we cannot apply the Hopf argument to the whole torus we can
use it in a connecting square $R$ to claim that $f^+$ is constant on
all of $W^s(R)$ and $f^-$ is constant on all of $W^u(R)$ with the two
constants coinciding. Note that we say here (and we mean it) ``all of
$W^{s(u)}$'' and not almost all. Indeed, first of all $f^+$ is
constant on each of the stable leaves in $W^s(R)$. Further let us fix
an ustable leaf in $W^u(R)$. The stable leaves from $W^s(R)$ intersect
this unstable leaf in $f$-typical points, except for a set of stable
leaves of total measure zero. Hence excluding these exceptional stable
leaves the value of $f^+$ on the stable leaves has to coincide with the
constant value of $f^-$ on the distinguished unstable leaf. We
conclude that $f^+$ is constant almost everywhere on $W^s(R)$ and the
constant is equal to the constant value of $f^-$ on the unstable leaf.
Since we could have used any other unstable leaf in $W^u(R)$ it
follows that $f^-$ is constant on all of $W^u(R)$. By symmetry $f^+$
is constant on all of $W^s(R)$. (The reader must have noticed the
implicit use of the Fubini Theorem in the arguments above.  It is only
natural since the stable and unstable leaves are parallel segments. In
the nonlinear case one has to use the ``absolute continuity" of the
foliations into stable and unstable manifolds.  This property is all
that we need, to make the present argument work.)

 Further for two connecting squares $R_1$ and $R_2$ which are
immediate neighbors $f^+$ is constant on $W^s(R_1)\cup W^s(R_2)$ and
$f^-$ is constant on $W^u(R_1)\cup W^u(R_2)$ with the two constants
coinciding.  Indeed at least one of the intersections $W^u(R_1)\cap
W^u(R_2)$ (if one square is above the other) or $W^s(R_1)\cap
W^s(R_2)$ (if one square is next to the other) must have positive measure
and hence is nonempty, forcing the constant value of 
$f^+$ or $f^-$ to be the same for both squares.

After this observation we proceed to prove that the time average of
$f$ is almost everywhere constant in $\U$. To that
end let $y,z \in \U$ be two $f$-typical points with $f$-typical leaves,
$W^u(y)$ and $W^s(z)$ respectively.  Our goal is to prove that $f^-(y)
= f^+(z)$.

 We say that $W^u(y)\left(\Cal W^s(z)\right)$ intersects completely a
column (row) of squares in $\Cal G_n(c)$ if it is connecting in one of
the squares of the column (row).  The Sinai Theorem allows us to claim
that, for sufficiently large $n$, $W^u(y)$ intersects completely at
least one column of connecting squares in $\Cal G_n(c)$,
i.e. a column in which all the squares are connecting, and $W^s(z)$
intersects completely at least one row of connecting squares.  Indeed,
suppose to the contrary that every column of squares in $\Cal G_n(c)$
intersected completely by $W^u(y)$ contains at least one non-connecting
square.  Since the number of columns intersected completely by
$W^u(y)$ grows linearly with $n$ and the measure of one square in
$\Cal G_n(c)$ is $\frac 1{n^2}$, we obtain that the measure of the
union of non-connecting squares would be $\text{\it O}(\frac 1n)$ which
contradicts the Sinai Theorem.  (Here we have used the fact
that the squares in $\Cal G_n(c)$ cannot overlap more than $k(c)$ times.)
\par
Let us fix a column and a row of connecting squares which are
intersected completely by $W^u(y)$ and $W^s(z)$ respectively.  Let $R$
be the (unique) square which belongs both to the column and the row.  Let
further $R_1$ denote a square in which $W^u(y)$ is connecting and
$R_2$ denote a square in which $W^s(z)$ is connecting. By the
construction $y \in W^u(R_1)$ and $f^-$ is constant on the, possibly
disjoint, set $W^u(R_1)\cup
W^u(R)$.  Similarly $z \in W^u(R_2)$ and $f^+$ is constant on
$W^s(R_2)\cup W^s(R)$.  It follows that $f^-(y) = f^+(z)$.  In view of
the arbitrariness in the choice of the $f$-typical leaves $W^u(y)$ and
$W^s(z)$ we obtain that the time average of $f$ must be constant in
$\U$.

\par
To finish the proof let us consider a $T$-invariant measurable subset $A$.
Let $g$ be the indicator function of $A$ and
$$
f_n \rightarrow g \ \ \text{in} \ \ L^1(\Bbb T^2,\,\mu)
$$
be a sequence of uniformly bounded continuous approximations
to the indicator function.
We will use  the fact that the time average is continuous with
respect to the $L^1$ norm to establish that the time average of $g$
must be constant on $\U$. Indeed, if we denote by $\|\cdot\|_{\lower
2pt\hbox{$\scriptstyle 1$}}$ the $L^1(\To^2,\,\mu)$ norm, then
$$
\aligned
\|f_n^+-g^+\lon=&
\biggl\|\lim_{N\to\infty}\frac 1N \sum_{i=1}^N\left (
f_n\circ T^i-g\circ T^i\right )\Lon\\
=&\lim_{N\to\infty}\frac 1N
\biggl\|\sum_{i=1}^N\left(f_n\circ T^i-g\circ T^i\right)\Lon
\endaligned
$$
by the Lebesgue Dominated Convergence Theorem.
 
Using the invariance of the measure we get
$$
\|f_n^+-g^+\lon \leq
\lim_{N\to\infty}\frac 1N
\sum_{i=1}^N\biggl\|\left(f_n\circ T^i-g\circ T^i\right)\Lon
= \bigl\|f_n-g\lon
$$
Since the time averages $f_n^+$ of $f_n$ are all constant 
(almost everywhere) on
$\U$ the
above
inequality implies that the  time average $g^+$ is constant
(almost everywhere) on
$\U$. But the invariance of $A$ forces $g^+ = g$  so that either
$\U\setminus A$ or $\U\cap A$ has measure zero. In view of the
arbitrariness of the invariant set $A$ it follows that $\U$ must belong
to one ergodic component.
\enddemo
 
\vskip.7cm
\subhead \S 3. PROOF OF THE SINAI THEOREM \endsubhead
\vskip.4cm
 
The proof of the Sinai Theorem does not require a rigid geometric
structure of the
coverings $\Cal G_n(c)$; it holds for any sequence of coverings by
squares with side $\frac 1n$ as long as there is a uniform bound on
the number of squares covering one point. However, the
lattice structure of the centers of the squares in $\Cal G_n(c)$ allows to
work with columns and rows of squares, as we did in the
above application of the Sinai Theorem.
\par
The first step in the proof is the choice of $\alpha_0$.
To that end we consider the smallest  sector $\Cal C$ in $\Bbb R^2$
symmetric about the horizontal (unstable) line which contains the
lines with the two directions of the sides of $\M^-$, i.e., the
directions of the segments in $\Cal S^-$. Let
$$
\Cal C = \{(\xi,\,\eta)\;|\;\ |\eta| \leq \kappa (a) |\xi|\}.
$$
It can be checked that $\kappa (a) < 1$ for any $a \neq 0$. We put
$\alpha_0 = \frac 12 (1-\kappa (a)).$ The reason for this choice is
that, for any square with vertical and horizontal sides
crossed by a line with the direction contained in $\Cal C$,
the shaded area in Figure 4 does not exceed $1-2\alpha$ part of the area
of the square.

\topinsert
\vskip 3in
\hsize=4.5in
\raggedright
\noindent{\bf Figure 3} Leaves cut by a line with direction contained in the
sector.
\endinsert

Let us observe that all of the segments in
$\bigcup_{i=0}^{+\infty} T^i \Cal S^- $ have directions contained in 
the sector $\Cal C$. Indeed a linear
hyperbolic map pushes lines towards the unstable direction except for
the stable line, which stays put.

It follows from the construction of the unstable leaves (Proposition 1.1)
that an unstable leaf has endpoints on  forward images of $\Cal S^-$ under
$T$. Hence if an unstable leaf is short in a square then the square must
be intersected by
$$
\bigcup_{i=o}^{+\infty} T^i \Cal S^- \ .
$$
Although this does not look like a severe restriction,
since we can expect that
the last set is  dense, it has far reaching consequences. The reason being,
heuristically, that the singularity
lines $T^i \Cal S^- $ become more and more horizontal as $i
\rightarrow +\infty$ and they cannot cut effectively unstable leaves
which are themselves horizontal.
 
We claim that, for any fixed $M \geq 1$, the
singularity lines
$$
\Cal S^-_M = \bigcup_{i=0}^{M} T^i \Cal S^-
$$
by themselves can produce only few squares which are not
$\alpha$-connecting so that their total
measure is $\text{\it O}(\frac 1{n^2})$. To make this precise (and clear)
we introduce an auxiliary notion of an $M$-bad square in a covering
$\Cal G_n(c)$. We say that a square $R \in \Cal G_n(c)$ is $M$-bad if
 the measure of the
set of points $y\in R$ such
that the unstable leaf $W^u(y)$ has an
endpoint in
$R \cap  \Cal S^-_M $
(so that it is short in $R$) is
greater than $1-2\alpha$ part of the measure of the square. (Loosely
speaking a square is $M$-bad if it is not connecting because
of the singularity
lines in $ \Cal S^-_M $.)
 
If a square $R$ intersects only one segment in $\Cal S^-_M $ then the
measure of points in $R$ whose unstable leaves have endpoints on the
intersection of this segment with $R$ does not exceed $1-2\alpha_0 =
\kappa (a)$ part of the measure of the square since the direction of
the segment is in the sector $\Cal C$. Hence an $M$-bad square has to
intersect at least two segments in $ \Cal S^-_M $.  But the singularity set
$ \Cal S^-_M $ is a fixed finite collection of closed
segments with only fixed finite number of intersection points (i.e.,
belonging to several segments).  Away from the intersection points the
segments are fairly wide apart and a small square cannot extend from
one to another, see Figure 5.  Hence, for sufficiently large $n$, an
$M$-bad square in $\Cal G_n(c)$ cannot be farther from one of the
intersection points than $\frac {const}n$. It follows that the total
measure of $M$-bad squares does not exceed $\frac {const}{n^2}$, where
the constant depends only on $a,\,c,\,\alpha$ and $M$.

\topinsert
\vskip 3in
\hsize=4.5in
\raggedright
\noindent{\bf Figure 4} Singularity lines.
\endinsert

In this way we took care (in some sense) of the finite number of
singularity lines in  $\Cal S^-_M$;
we now face the problem of controlling
the effects of the `tail' $\bigcup_{i=M+1}^{+\infty} T^i \Cal S^-$.
 
Let us suppose that a square $R \in \Cal G_n(c)$ is not
$\alpha$-connecting and it is not $M$-bad.
Hence at least $\alpha$ part of its area is covered by short leaves with
endpoints in 
$$
R \cap \bigcup_{i=M+1}^{+\infty} T^i \Cal S^-.
$$
Let $\Cal
W^u(y)$ be such a leaf short in $R$ with an endpoint on $T^i \Cal S^-$.
Then
$$
T^{-i}\left(W^u(y) \cap R \right) \subset X_{t_i}
$$
where $t_i = n^{-1}\lambda^{-i}$ and, as before,
$X_t = \{ x \in \M^- \;|\; d(x,\,\Cal S^-) \leq t\}$. Indeed, under
the action of $T^{-1}$, an
unstable leaf contracts by a factor of $\lambda$ and the length of the
part of $W^u(y)$ in $R$ does not exceed $\frac 1n$.
 
In view of this observation we can claim that each  square which is
not
$\alpha$-connecting and which is not $M$-bad
has at least $\alpha$ part of its area covered by
$$
\bigcup_{i=M+1}^{+\infty} T^i X_{t_i}.
$$
Since each point in $\U$ is covered by, at most,
$k(c)$ squares from $\Cal G_n(c)$, then the measure
of the union of squares in $\Cal G_n(c)$ which are not
$\alpha$-connecting and which are not $M$-bad does not exceed
$$
k(c)\times \frac 1\alpha\sum_{i=M+1}^{+\infty} \frac {const}{n\lambda^i} =
\frac 1n \left(\frac {k(c)}{\alpha}
\sum_{i=M+1}^{+\infty} \frac {const} {\lambda^i}\right),
$$
(here the constant is equal to the total length of $\Cal S^-$).
We have thus estimated the measure of the union of  squares
in $\Cal G_n(c)$, which are not $\alpha$-connecting and which are not
$M$-bad, by the size of an individual square times the $M$-tail of a
fixed convergent
series. Some of the readers may have noticed that this completes the
proof. For clarity, let us do it explicitly.
 
Let us take an arbitrary $\epsilon > 0$.
We choose and fix $M= M(\epsilon)$ so large
that the last series does not exceed $\frac \epsilon {2n}$, i.e.,
$$
\frac {k(c)}{\alpha} \sum_{i=M+1}^{+\infty} \frac {const} {\lambda^i}
< \frac \epsilon 2 .
$$
Given $M$ we can still choose $n_0 = n_0(\epsilon,\,M)$
so large that, for any $n \geq n_0$,
the measure of the union of $M$-bad squares in $\Cal G_n(c)$ is less than
$\frac \epsilon {2n}$. To estimate the measure of the union of squares
in $\Cal G_n(c)$,
for $n \geq n_0$, which are not $\alpha$-connecting we split them into
those which are $M$-bad and those which
are not. For both families of squares the measure of their union is less than
$\frac \epsilon {2n}$. This proves our claim.
\qed
 
\proclaim {Remark 3.6}
\endproclaim
Let us point out that the property that the sector $\Cal C $, defined
by the directions of the segments in $\Cal S^-$, is
sufficiently narrow ($\kappa(a)<1$) can be relaxed.
For a general hyperbolic piecewise linear map
it is sufficient that the
segments in $\Cal S^-$ are not parallel to the stable direction. In such a case
we
can find a natural $N$ such that all the segments in
$\bigcup_{i=N+1}^{+\infty} T^i \Cal S^-$ have directions contained in
a chosen narrow sector
$\Cal C$ ( $N$ is the number of iterates of $T$ which do not
put the singularity lines $\Cal S^-$ into the chosen sector $\Cal C$).
Then the argument above applies to any square neighborhood $\U$
which does not intersect
$$
\Cal S^-_N = \bigcup_{i=0}^{N} T^i \Cal S^-\ .
$$
Similarly in the version of the Sinai Theorem for the stable leaves we
would have arrived at a natural $N'$ such that the claim holds for
any square $\U$ which does not intersect
$$
\Cal S^+_{N'} = \bigcup_{i=0}^{N'} T^{-i} \Cal S^+\ .
$$
 
Hence, it follows from Proposition 2.3 that any open
square, with horizontal and vertical sides, which does not intersect
$ \Cal S^-_N \cup \Cal S^+_{N'} $
belongs to one ergodic component.
This implies that the partition of $\Bbb T^2$ into
ergodic
components is coarser than the partition into (open) connected components of
$$
\Bbb T^2 \setminus \left( \Cal S^-_N \cup \Cal S^+_{N'} \right)\  .
$$
Since $ \Cal S^-_N \cup \Cal S^+_{N'} $ is a finite collection of
segments we
obtain that there are at most finitely many ergodic components. To argue that
there  is only one component let us note that $ \Cal S^-_{N-1} \cup \Cal S^+_{N'} $
and $T^{N} \Cal S^- $ intersect in at most finitely many points which split the
segments in
$T^{N} \Cal S^- $ into finitely many segments  $\{ I_k\}_{k=1}^{K_N}$ so that
the interior of every $I_k$ lies in the boundary of at most two connected
components
of $\Bbb T^2 \setminus \left( \Cal S^-_N \cup \Cal S^+_{N'}  \right)$, i.e., 
it has only one
connected component on each side. Suppose that for such a segment $I_k$ is
in the boundary of two
different ergodic components. Then $T I_k$ is also in the boundary of two
different ergodic components. But $T I_k$ and
$\Cal S^-_N \cup \Cal S^+_{N'}$ have only finitely many points of
intersection, so
that whole open sub-intervals of $T I_k$ must end up inside one connected
component of
$ \Bbb T^2 \setminus  \left( \Cal S^-_N \cup \Cal S^+_{N'}  \right)$ and thus it must
have the same ergodic
component on both sides. This contradiction implies that $I_k$ does
not take part in the splitting of
$\Bbb T^2$ into ergodic components so we can drop it. In this way we can drop
all of
$T^{N} \Cal S^- $ and claim that the partition into ergodic components is
coarser than the partition into connected components of
$$
\Bbb T^2 \setminus \left( \Cal S^-_{N-1} \cup \Cal S^+_{N'} \right) \  .
$$
It is now clear that we can proceed by dropping  $T^{N-1} \Cal S^- $ and
$T^{-N'} \Cal S^+ $ as possible boundaries for the ergodic components and
arriving eventually at $\Cal S^+ \cup \Cal S^-$ as the only possible
boundaries we see that even
these can be dropped. Hence there is only one ergodic component.
 
Let us spell out the property of $T$ which is basic in this argument:
 
{\it Although some points of $\Cal S^-$ return to $\Cal S^-$ under iterates of
$T$, no interval in $\Cal S^-$ can do it.}

\vskip.7cm
\subhead \S 4. SECTORS IN A LINEAR SYMPLECTIC SPACE \endsubhead
\vskip.4cm
For the convenience of the reader we will repeat here some of the
material from \cite {W3} and  \cite {LW}.
 
Let $\Cal W$ be a linear symplectic space of dimension $2d$ with the
symplectic form $\omega$. For instance we call $\Cal W =
\Bbb{R}^{d}\times\Bbb{R}^{d}$ the standard linear symplectic space if
$$ \omega (w_1,w_2) = \langle\xi^1,\eta^2\rangle-\langle\xi^2,\eta^1\rangle,$$
$$\text{where} \ w_i =(\xi^i,\eta^i),\ i =1,2, \ \ \text{and}
\ \ \langle\xi,\eta\rangle = \xi_1\eta_1 + \dots + \xi_d\eta_d .$$
 
The symplectic group $Sp\left(d,\Bbb{R}\right)$ is the group of linear maps
of $\Cal W$ ($2d\times 2d$ matrices if $\Cal W = \Bbb{R}^{d}\times
\Bbb{R}^{d}$) preserving the symplectic form i.e., $L\in
Sp\left(d,\Bbb{R}\right)$ if $$\omega(Lw_1,Lw_2)=\omega(w_1,w_2)$$ for every
$w_1,w_2 \in \Cal W$.
 
By definition a Lagrangian subspace of a linear symplectic space $\Cal W$
is a $d$-dimensional subspace on which the restriction of $\omega$ is
zero (equivalently it is a maximal subspace on which $\omega$
vanishes).
 
\proclaim{Definition 4.1}
Given two transversal Lagrangian subspaces $V_1$ and $V_2$ we define
the sector between $V_1$ and $V_2$ by
$$ \Cal{C} = \Cal{C}\left(V_1,V_2\right) =
\{w\in \Cal W\ |\ \omega(v_1,v_2)\geq 0
\  \text{for}\ w=v_1+v_2, v_i\in V_i, i=1,\, 2\}$$
\endproclaim
 
Equivalently, if we define the quadratic form associated with an
ordered pair of transversal Lagrangian subspaces,
$$\Cal{Q}(w) = \omega (v_1,v_2) $$
where $ w =v_1 + v_2$, is the unique decomposition of $w$ with the
property $v_i\in V_i, i=1,2,$
then we have
$$\Cal{C} =\{w\in \Cal W\ |\ \Cal{Q}(w)\geq 0\}.$$
 
In the case of the standard symplectic space,
$V_1=\Bbb{R}^d\times\{0\}$ and $V_2=\{0\}\times\Bbb{R}^d$ we get
$$\Cal{Q}\left((\xi,\eta)\right)=\langle\xi,\eta\rangle$$ and $$\Cal{C} =\{
 (\xi,\eta)\in
\Bbb {R}^{d}\times \Bbb{R}^{d}\ |\ \langle\xi,\eta\rangle\geq 0\}.$$
We will refer to this $\C$ as the standard sector. Since any two pairs
of transversal Lagrangian subspaces are symplectically equivalent we
may consider only this case without any loss of generality. In the
following we will alternate between the coordinate free geometric
formulations and this special case. On the one hand,
coordinate free formulations are
important because we need to apply these concepts to the case of the
derivative map which in general acts between two different tangent
subspaces, each one with its
preferred sector.  On the other hand, it turns out that many arguments
are greatly simplified by resorting to these special coordinates.
 
It is natural to ask if a sector determines uniquely its sides. It is
not a vacuous question since, for $d>1$, there are
many Lagrangian subspaces in the boundary of a sector. The answer is positive.
 
\proclaim{Proposition 4.2}
For two pairs of transversal Lagrangian subspaces $V_1,V_2$ and
$V_1',V_2'$ if $$ \Cal{C}\left(V_1,V_2\right) = \Cal{C}\left(V_1',V_2'\right) $$
 then $$V_1
= V_1' \ \ \text{and} \ \ V_2 = V_2'.$$ Moreover $V_1$ and $V_2$ are
the only isolated Lagrangian subspaces contained in the boundary of
the sector $ \Cal{C}\left(V_1,V_2\right)$.
\endproclaim
The proof of this Proposition can be found in \cite {W3}.
 
Based on the notion of the sector between two transversal Lagrangian
subspaces (or the quadratic form $\Cal{Q}$) we define two monotonicity
properties of a linear symplectic map. By $int\C$ we denote the interior of the
 sector,
i.e.,
$$
int\C =\{w \in \Cal W | \Q(w) > 0\}.
$$
\proclaim{Definition 4.3}
Given the sector $\Cal C$ between two transversal Lagrangian subspaces
we call a linear symplectic map $L$ monotone if $$ L\Cal{C}\subset\Cal{C} $$
and strictly monotone if $$
L\Cal{C}    \subset int\Cal{C} \cup \{0\}. $$
\endproclaim
 
A very useful characterization of monotonicity is given in the following
 
\proclaim{Theorem 4.4}
$L$ is (strictly) monotone if and only if $\Cal{Q}\left(Lw\right)\geq
\Cal{Q}\left(w\right)$ for
every $w\in \Cal W$ ($\Cal{Q}\left(Lw\right) > \Cal{Q}\left(w\right)$ for
every $w\in \Cal W ,\ w \neq 0$).
\endproclaim
 
The fact that monotonicity implies the increase of the quadratic form
defining the cone is a manifestation of a very special geometric
structure of a sector and does not hold for cones defined by general
quadratic forms. The proof of the theorem relies on the factorization
(4.7), we postpone then the proof until such factorization has
been established.
 
For a pair of transversal Lagrangian subspaces $V_1$ and $V_2$ and a linear map
$L:\Cal W\to \Cal W$ we can define the following `block' operators:
$$\aligned
A:V_1\to V_1 ,\ & B:V_2\to V_1 \\
C:V_1\to V_2 ,\ & D:V_2\to V_2.
\endaligned$$
They are uniquely defined by the requirement that for any $v_1\in V_1, v_2\in
 V_2$
$$ L\left(v_1+v_2\right)= Av_1 + Bv_2 + Cv_1 + Dv_2 .$$
 
We will need the following Lemma.
 
\proclaim{Lemma 4.5}
If $L$ is monotone with respect to the sector defined by $V_1$ and $V_2$ then
$LV_1$ is transversal to $V_2$ and $LV_2$ is transversal to $V_1$.
\endproclaim
\demo{Proof}
Suppose that, to the contrary, there exists $0\neq \bar{v}_{1} \in V_1$ such
that $L\bar{v}_{1} \in V_2 .$ We choose $\bar{v}_{2} \in V_2$ so that
$$\Cal{Q}\left(\bar{v}_{1} + \bar{v}_{2}\right) =
 \omega\left(\bar{v}_{1},\bar{v}_{2}\right) > 0.$$
 
We have also $$
\omega\left(\bar{v}_{1},\bar{v}_{2}\right)=
\omega\left(L\bar{v}_{1},L\bar{v}_{2}\right)=
\omega\left(L\bar{v}_{1},B\bar{v}_{2}+ D\bar{v}_{2}\right)=
\omega\left(L\bar{v}_{1},B\bar{v}_{2}\right).$$
 
Let $v_{\epsilon} = \bar{v}_{1}+\epsilon \bar{v}_{2}.$ We have that for
$\epsilon >0$\ $ v_{\epsilon} $belongs to $ \text{int}\Cal{C}.$ Hence also
$\Cal{Q}\left(Lv_{\epsilon}\right)\geq 0 $ for $\epsilon >0.$ On the other hand
 
$$ \Cal{Q}\left(Lv_{\epsilon}\right) = \epsilon ^2
 \omega\left(B\bar{v}_{2},D\bar{v}_{2}\right) -
\epsilon\omega\left(L\bar{v}_{1},B\bar{v}_{2}\right)$$
which is negative for sufficiently small positive $\epsilon.$
 
This contradiction proves the Lemma.
\enddemo
 
It follows, from Lemma 4.5, that the operators $A:V_1\to V_1$ and
$D:V_2\to V_2$ are invertible.
 
We switch now to coordinate language.
Let $$ L = \left(\matrix A& B\\ C& D\endmatrix\right) $$
be a symplectic map of the standard symplectic space
$\Bbb{R}^{d}\times\Bbb{R}^{d}$ monotone with respect to the standard sector.
$A,B,C,D$ are now just $d\times d$ matrices.
 
Let us describe those symplectic matrices which are
monotone in the weakest sense, namely they preserve the quadratic form
$\Cal{Q}.$ We will call such matrices $\Cal{Q}$-isometries. Obviously
a $\Cal{Q}$-isometry maps the sector onto itself. The converse is also
true.
 
\proclaim {Proposition 4.6}
If $L$ is a linear symplectic map and $$
L\Cal{C} = \Cal{C} $$
then $$ L = \left(\matrix A& 0\\ 0& A^{*-1}\endmatrix\right) . $$
In particular it preserves the quadratic form $\Cal{Q}$
$$\Cal{Q}\circ L = \Cal{Q}  .$$
\endproclaim
\demo{Proof}
If $L\Cal{C} = \Cal{C} $ then $L$ maps  also the boundary of the sector
$\Cal{C}$ onto itself. It follows from Proposition 4.2 that both sides
of the sector stay put under $L$. Hence $B=C=0.$
By symplecticity $D = A^{*-1}.$
\enddemo
 
By Lemma 4.5 given a monotone $L$ we can always factor out the following
$\Cal{Q}$-isometries on the left
$$
L =
\left(\matrix A& B\\ C& D\endmatrix\right) =
\left(\matrix A& 0\\ 0& A^{*-1}\endmatrix\right)\left(\matrix
I& R \\ P &\cdot \endmatrix\right)
$$
 
($P$ and $R$ are uniquely determined).
Symplecticity of $L$ forces $R,\,P$ symmetric and $RP-A^*D=I$,
which allows the further unique factorization
$$
L =
\left(\matrix A& 0\\ 0& A^{*-1}\endmatrix\right)\left(\matrix I& 0 \\ P &I
\endmatrix\right)\left(\matrix I& R \\ 0 &I \endmatrix\right).
\tag {4.7}
$$
Moreover monotonicity forces $P$ and $R$ to be positive semidefinite
($P\ge 0,\,R\ge 0$).
Strict monotonicity means that $P$ and $R$ are positive definite
($P> 0,\, R> 0$).
These claims follow from the following
\demo{Proof of Theorem 4.4}
Using the above factorization we get for $w = (\xi,\,\eta)$
$$
\Cal{Q}(Lw) = \langle\xi,\,\eta\rangle + \langle R\eta,\,\eta\rangle +
\langle P(\xi+R\eta),\,\xi+R\eta\rangle.
$$
Putting $\eta = 0$ we obtain that $P\geq 0$.
To show that also $R \geq 0$ let us consider an eigenvector
$\eta_0$ of $R$ with eigenvalue $\lambda$ and let $\xi = a\eta_0$.  We
get that if $a\geq 0$ then $w = (\xi,\,\eta_0) \in \Cal C$ so that
$\Cal{Q}(Lw)  \geq 0$. It follows that
$$
(a+\lambda)\langle\eta,\,\eta\rangle +
(a+\lambda)^2\langle P\eta,\,\eta\rangle\geq 0 .
$$
This implies immediately that
$\lambda \geq 0$. This proves the monotone version of the Theorem.
The strictly monotone version is obtained in a similar way.
\enddemo
 
As a byproduct of the proof we get the following useful observation
 
\proclaim{Proposition 4.8}
A monotone map $L$ is strictly monotone if and only if
$$LV_i  \subset \text{int }\Cal{C} \cup \{ 0\}  ,\  i=1,2. $$
\endproclaim
\qed
 
The following Proposition simplifies computations with monotone maps.
\proclaim {Proposition 4.9}
If
$$
L = \left(\matrix A& B\\ C& D\endmatrix\right)
$$ is a strictly monotone map then by multiplying it by $\Cal{Q}$-isometries
on the left and on the right we can bring it to the form
$$
\left(\matrix I& I\\ T& I+T\endmatrix\right)
$$
where $T$ is diagonal and has the same eigenvalues as $C^*B$.
\endproclaim
\demo{Proof}
The factorization of the monotone map $L$ yields
$$
\left(\matrix A& 0\\0& A^{*-1}\endmatrix\right)L =
\left(\matrix I& R\\ P& I+PR\endmatrix\right) $$
where $P > 0$, $R > 0$ and $PR = C^*B$.
 
We have further
$$\left(\matrix R^{-\frac{1}{2}}& 0\\ 0&R^{\frac{1}{2}}\endmatrix\right)
\left(\matrix I& R\\ P& I+PR\endmatrix\right)
\left(\matrix R^{\frac{1}{2}}& 0\\ 0& R^{-\frac{1}{2}}\endmatrix\right) =
\left(\matrix I& I\\ K& I+K\endmatrix\right)
$$
where $K= R^{\frac{1}{2}}PR^{\frac{1}{2}}$ has the same eigenvalues as $C^*B=
PR.$
 
Finally if $F$ is the orthogonal matrix which diagonalizes $K$, i.e.,
$F^{-1}KF$ is diagonal, then
$$
\left(\matrix F^{-1}& 0\\ 0& F^{-1}\endmatrix\right)\left(\matrix I& I\\K&
I+K\endmatrix\right)\left(\matrix F& 0\\ 0& F\endmatrix\right) =
\left(\matrix I& I\\ T& I+T\endmatrix\right)
$$
has the desired form with $T = F^{-1}KF$ having the
same eigenvalues as $C^*B$.
\enddemo

Let us note that in the last Proposition we can ask for the diagonal
entries of $T$ to be ordered because any permutation  of the entries
can be accomplished by an appropriate $\Cal Q$-isometry.
 
\vskip.7cm
\subhead \S 5. THE SPACE OF LAGRANGIAN SUBSPACES CONTAINED IN A SECTOR
\endsubhead
\vskip.4cm
 
Let us fix a sector $\C = \C(V_1,\,V_2)$ between two transversal
Lagrangian subspaces $V_1$ and $V_2$. We say that a Lagrangian
subspace $E$ is strictly contained in $\C$ if
$$
E \subset int\,\C \cup \{0\}.
$$
We denote by $Lag(\C)$ the manifold of all such Lagrangian subspaces and by
$\widehat{Lag}(\C)$ its closure in the Lagrangian Grassmanian, i.e.,
$\widehat{Lag}(\C)$ is
the set of all Lagrangian subspaces contained in $\C$.
 
We will introduce a metric and a partial order into $Lag(\C)$.
This will allow us to extend to the multidimensional case ($d>1$) the most
relevant features of the two dimensional case ($d=1$).
Let
$$
\pi_i : \Cal W \to V_i, \ i = 1,2,
$$
be the natural projections, i.e.,
$$
w = \pi_1w +\pi_2w\ \ \text{ for every}\ \ w \in \Cal W.
$$
 If a Lagrangian
subspace $E$ is strictly contained in $\C$ then $\pi_iE = V_i, \ i = 1,2$,
so $\pi_i|_E$ (the restriction of $\pi_i$ to the subspace $E$) is a one
to on map of $E$ onto $V_i$.
 
With every subspace $E \in Lag(\C)$
we can associate a positive definite quadratic form on $V_1$ obtained by the
 formula
$$
\Q \circ \left(\pi_1|_E\right)^{-1}.
$$
It will turn out that this is actually a one-to-one correspondence between
 positive
definite quadratic forms on $V_1$ and Lagrangian subspaces contained
strictly in $\C$.
\proclaim{Definition 5.1}
For two Lagrangian subspaces $E_1,E_2 \in Lag(\C)$ we define the relation $E_1
 \leq E_2$
($E_1 < E_2$) by the inequality of the corresponding quadratic forms
$$
\Q \circ \left(\pi_1|_{E_1}\right)^{-1} \leq (<) \Q \circ
 \left(\pi_1|_{E_2}\right)^{-1}.
$$
We define the distance of two Lagrangian subspaces $E_1,E_2 \in Lag(\C)$ by
$$
d(E_1,\,E_2) =
\frac 12 \sup_{0\neq v \in V_1}
|\ln\Q \circ \left(\pi_1|_{E_1}\right)^{-1}(v) -
\ln\Q \circ \left(\pi_1|_{E_2}\right)^{-1}(v)|.
$$
\endproclaim
It is easy to see that $d(\cdot,\,\cdot)$ is indeed a metric.
 
There are  other ways to introduce the partial order and the metric.
The coordinate free definitions simplify some of the arguments
in the following.
For equivalent definitions of the metric see \cite {LW}, \cite {Ve}.
Theses definitions are justified by the following theorem.
\proclaim{Theorem 5.2}
For two transversal Lagrangian subspaces $E_1,E_2 \in Lag(\C)$
$$
E_1 < E_2 \ \ \text{ if and only if } \ \ \C(E_1,\,E_2) \subset \C(V_1,\,V_2).
$$
Further if $E_1 < E_2$ then for a Lagrangian subspace $E \in Lag(\C)$
$$
E \subset \C(E_1,\,E_2) \ \ \text{ if and only if } \ \ E_1  \leq E \leq  E_2.
$$
\endproclaim
 
\proclaim{Corollary 5.3}
If $E_1,E_2 \in Lag(\C)$ and $E_1 < E_2$ then the diameter of the set
$\widehat{Lag}\left(\C(E_1,\,E_2)\right)$ in $Lag(\C)$
is equal to the distance of $E_1$ and $E_2$.
\endproclaim \qed
 
We will prove Theorem 5.2 at the end of this Section.
 
Let us introduce a convenient parametrization of $Lag(\C)$ by
symmetric positive
definite matrices. We consider the standard sector $\C$ in
$\Bbb{R}^{d}\times \Bbb{R}^{d}$ with
$V_1 = \Bbb R^d  \times \{0\}$ and $V_2 = \{0\} \times \Bbb R^d$.
Let $U : \Bbb{R}^d \to \Bbb{R}^d$ be a linear map and
$$
gU =\{\left(\xi,\eta\right)\in \Bbb{R}^{d}\times\Bbb{R}^{d}\;|\;\eta = U\xi\}
$$
be its graph. The linear subspace $gU$ is a Lagrangian subspace if and only if
$U$ is symmetric and further for a symmetric $U$ its graph
$gU \subset \Cal{C}$
if and only if $U \ge 0.$ Every Lagrangian subspace in $Lag(\C)$
is transversal
to $V_2$ so that it is a graph of a linear map as above.
We will find the following Lemma useful.
\proclaim {Lemma 5.4}
If a Langrangian subspace $E \subset \Cal C(V_1,\,V_2)$
is transversal to both $V_1$ and $V_2$ then it is strictly contained
in the sector.
\endproclaim
\demo{Proof}
We use the coordinate description of the standard sector. Thus the
Lagrangian subspace $E$ being transversal to $V_2$ is the graph of a
symmetric positive semidefinite matrix. Since $E$ is also transversal
to $V_1$ the matrix is nondegenerate and hence positive definite.
It follows immediately that $E$ is strictly contained in the sector.
\enddemo
 
We have obtained a
one-to-one correspondence between Lagrangian subspaces in $Lag(\C)$
and symmetric positive
definite matrices. The quadratic form on $V_1$ introduced in Definition 5.1
becomes the form defined by the positive definite matrix.
The partial order becomes the familiar partial order between symmetric
matrices.

The image of a Lagrangian subspace under a symplectic linear map is again a
Lagrangian subspace. Moreover monotone maps take Lagrangian subspaces
strictly contained in $\C$ into Lagrangian subspaces
strictly contained  in $\C$.
Hence a monotone map $L$ defines a map of $Lag(\C)$ into itself. We
will denote it
again by $L: Lag(\C) \to Lag(\C)$. To simplify notation we will also write $U$
instead of $gU$. We have that
$$ L = \left(\matrix A& B\\ C& D\endmatrix\right) $$
acts on Lagrangian subspaces by the following M\"obius
transformation
$$
LU=\left(C+DU\right)\left(A+BU\right)^{-1}.
$$
 
In particular the action of a $\Q$-isometry
$$ L = \left(\matrix A& 0\\ 0& A^{*-1}\endmatrix\right) $$
is given by
$$
LU = A^{*-1}UA^{-1}.
$$
By putting $A = U^{\frac 12}$ we see that any $U > 0$ can be
mapped onto identity
matrix $I$. Thus $\Q$-isometries act transitively on $Lag(\C)$. Moreover it is
not hard to see that
\proclaim{Proposition 5.5}
The action of a $\Q$-isometry on $Lag(\C)$ preserves the partial order and
the metric.
\endproclaim \qed
 
Let $E_0 = \{(\xi,\,\eta)\ |\  \xi = \eta\}$. By straightforward computations we
find that
$$
\aligned
\C(V_1,\,E_0) = \{(\xi,\,\eta)\ |\  \langle\xi,\,\eta\rangle -
\langle\eta,\,\eta\rangle
\geq 0\},\\
\C(E_0,\,V_2) = \{(\xi,\,\eta)\ |\  \langle\xi,\,\eta\rangle -
\langle\xi,\,\xi\rangle
\geq 0\}.
\endaligned
\tag{5.6}
$$
We get that
$$
\aligned
\C(V_1,\,E_0) \subset \C(V_1,\,V_2), \\
\C(E_0,\,V_2)\subset \C(V_1,\,V_2), \\
\C(V_1,\,E_0) \cap \C(E_0,\,V_2) = E_0.
\endaligned
\tag 5.7
$$
Because the group of $\Q$-isometries acts transitively on $Lag(\C)$
\thetag {5.7} holds not just for the special Lagrangian subspace
$E_0$ from \thetag {5.6}
but for any  Lagrangian subspace from $Lag(\C)$. (It just happens
that the easiest way to establish \thetag{5.7} is to do the
calculation in the standard sector.)
 
\proclaim{Proposition 5.8}
For two Lagrangian subspaces $E_1,E_2 \in Lag(\C)$ the following are equivalent
$$
\aligned
&(1) \ \ E_1 \leq E_2,\\
&(2) \ \ E_2 \subset \C(E_1,\,V_2),\\
&(3) \ \ E_1 \subset \C(V_1,\,E_2).
\endaligned
$$
\endproclaim
\demo{Proof}
We will be using the coordinate description of the standard sector.
Since the group of $\Q$-isometries acts
transitively on $Lag(\C)$ we can assume that $E_1$ is equal to $E_0$ from
\thetag{5.6}. Let $U_2$ be the positive definite matrix defining $E_2$.
We get from \thetag{5.6} that $E_2 \subset \C(E_0,\,V_2)$ if and only if
$U_2 \geq I$. Hence (1) is equivalent to (2). Similarly let $E_2$ be
equal to $E_0$ and $U_1$ be the positive definite matrix defining $E_1$.
Using \thetag{5.6} again we get that $E_1 \subset \C(V_1,\,E_0)$ if and only
if $U_1 -U_1^2 \geq 0$ which is equivalent to $U_1 \leq I$. This proves the
equivalence of (1) and (3).
\enddemo
 
\demo{Proof of Theorem 5.2}
If $E_1 < E_2$ then, by Proposition 5.8 and Lemma 5.4,
$E_2$ is strictly contained
in $\C(E_1,\,V_2)$. Using \thetag{5.7} we get
$$
\C(E_1,\,E_2) \subset \C(E_1,\,V_2) \subset \C(V_1,\,V_2).
$$
 
Suppose now that $\C(E_1,\,E_2) \subset \C(V_1,\,V_2)$. By Proposition 5.8
it suffices to show that $E_2 \subset \C(E_1,\,V_2)$. If it is not so then
there is $e_2 \in E_2$ which does not belong to $\C(E_1,\,V_2)$. Let us
consider  $v_1 = \pi_1e_2$ where $\pi_1 : \Cal W \to V_1$ is the
projection onto $V_1$ in the direction of $V_2$. Let further $e_1$ be
the unique element in $E_1$ such that $\pi_1e_1 = v_1$ (i.e.,
$e_1 = \left(\pi_1|_{E_1}\right)^{-1}v_1$). Clearly the difference
between the two vectors $v_2 = e_2 - e_1$ belongs to $V_2$.
Because $e_2 = e_1 + v_2$ and $e_2 \notin \C(E_1,\,V_2)$ we have
$\omega(e_1,\,v_2) < 0$ so that
$\omega(-e_1,\,e_2) > 0$. It follows that
$v_2 =e_2 - e_1 \in int\,\C(E_1,\,E_2) \subset int\,\C(V_1,\,V_2)$.
We have then reached a contradiction, since $v_2$ cannot belong
simultaneously to $V_2$ and to $ int\,\C(V_1,\,V_2)$.
The above contradiction proves that indeed
$E_2 \subset \C(E_1,\,V_2)$ which by Proposition 5.8 implies that
$E_1 < E_2$ (remember that $E_1$ and $E_2$ are assumed to be
transversal).
The first part of the Theorem is proven.
 
To prove the second part let $E_1 < E_2$ and $E \subset \C(E_1,\,E_2)$.
By Proposition 5.8 we get
$E_2 \subset \C(E_1,\,V_2)$. It follows in view of \thetag{5.7} that
$\C(E_1,\,E_2) \subset \C(E_1,\,V_2)$ and hence $E \subset \C(E_1,\,V_2) $
which is equivalent (again by Proposition 5.8) to $E_1 \leq E$. Similarly
we get $E \leq E_2$.
 
In the opposite direction if $E_1 \leq E < E_2$ then
by Proposition 5.8 $E_1$ and $E$ are strictly contained in $\C(V_1,\,E_2)$
and $E_1 \subset \C(V_1,\,E)$. Applying now the equivalence of (2) and
(3) in Proposition 5.8 to the case of $E_1,E \in Lag(\C(V_1,\,E_2))$
we get immediately $E \subset \C(E_1,\,E_2)$. The case of
$E_1 \leq E \leq E_2$ can be now treated by continuity.
\enddemo
 
Let us consider a special family of Lagrangian subspaces in the
standard sector:
the graphs of multiples of the identity matrix, i.e.,
for a real number $u$ let
$$
Z_u = \{(\xi,\,\eta)\ |\  \eta = e^u \xi\}.
$$
We have that
$$
d(Z_{u_1},\,Z_{u_2}) = \frac 12  | u_1 - u_2|.
$$
In the next Lemma we have chosen two numbers $u_2>u_1$.
\proclaim{Lemma 5.9}
If for a Lagrangian subspace $E \in Lag(\C)$
$$
d(Z_{u_1},\,E) \leq \frac 12 (u_2 - u_1)
$$
then
$$
E \leq Z_{u_2}.
$$
\endproclaim
\demo{Proof}
Let the Lagrangian subspace $E$ be the graph of a positive definite matrix $U$.
For every nonzero $\xi \in \Bbb R^d$, we have
$$
\ln\langle\xi,\,U\xi\rangle -\ln\langle\xi,\,e^{u_1}\xi\rangle \leq u_2 -u_1.
$$
It follows that, for every nonzero $\xi \in \Bbb R^d$,
$$
\ln\frac{\langle\xi,\,U\xi\rangle}{\langle\xi,\,\xi\rangle} \leq u_2 .
$$
We conclude that $U \leq e^{u_2}I$.
\enddemo
 
We will use the following consequence of the last Lemma.
\proclaim{Proposition 5.10}
Let $E_1 < E_2 $ be two Lagrangian subspaces contained strictly in
 $\C(V_1,\,V_2)$.
There is a symplectic map which maps the sector $\C(V_1,\,V_2)$
{\bf onto} the standard
sector $\C$ and the sector $\C(E_1,\,E_2)$ {\bf into} the sector
 $\C(Z_{-u},\,Z_u)$ if and only if
$d(E_1,\,E_2) \leq u$.
\endproclaim
\demo{Proof}
By a symplectic map we can map the subspace
$V_1$ onto $\Bbb R^d \times \{0\}$, the
subspace $V_2$ onto $\{0\} \times \Bbb R^d$ and $E_1$ onto $Z_{-u}$
(because $\Q$-isometries act transitively on $Lag(\C)$).
It follows from Lemma 5.9
that the sector $\C(E_1,\,E_2)$ will be then automatically mapped into
$\C(Z_{-u},\,Z_u)$.
 
The converse follows from the Corollary 5.3.
\enddemo
 
For aesthetical reasons we will be using Proposition 5.10 in a
different coordinate system obtained by the following linear
symplectic coordinate change
$$
\aligned
\xi' &= \frac 1{\sqrt 2}(\xi -\eta),\\
\eta' &= \frac 1{\sqrt 2}(\xi +\eta).
\endaligned
$$
 Let us introduce the family of sectors
$$
\C_{\rho} = \{(\xi,\,\eta)\ |\  \|\eta\| \leq \rho\|\xi\| \}
$$
for any real $\rho > 0$.
\proclaim{Proposition 5.11}
Let $E_1 < E_2 $ be two Lagrangian subspaces contained strictly in
 $\C(V_1,\,V_2)$.
There is a symplectic map which maps the sector $\C(V_1,\,V_2)$
{\bf onto} the sector $\C_{\rho^{-1}}$ and the sector $\C(E_1,\,E_2)$
{\bf into} the sector $\C_\rho$ if and only if
$$
d(E_1,\,E_2) \leq \ln \frac {1+\rho^2}{1-\rho^2},
$$
with $0<\rho<1$.
\endproclaim
\demo{Proof}
It is enough to define the coordinate change $L$, defined by
$$
\aligned
\xi' &= \frac 1{\sqrt 2}(\rho^{-\frac 12}\xi -\rho^{\frac 12}\eta),\\
\eta' &= \frac 1{\sqrt 2}(\rho^{-\frac 12}\xi +\rho^{\frac 12}\eta).
\endaligned
$$
A direct computation shows that, if $\rho<1$, $L\C_{\rho^{-1}}=\C$
and $L\C_{\rho}=\C(Z_{-u},\,Z_{u})$, with $u=\log\frac{1+\rho^2}{1-\rho^2}$.
The result follows then from Property 5.10.
\enddemo
\vskip.7cm
\subhead \S 6. UNBOUNDED SEQUENCES OF LINEAR MONOTONE MAPS
\endsubhead
\vskip.4cm
 
In this section we fix a sector $\C = \C(V_1,\,V_2)$ between two Lagrangian
subspaces. One can think that $\C$ is the standard sector.
We start by computing the coefficient of expansion of $\Q$ under
the action of a monotone symplectic map.
 
For a linear symplectic map $L$ monotone with respect to the sector $\Cal{C}$
we define the coefficient of expansion at $w \in int\Cal{C}$ by
$$\beta\left(w,L\right) =
\sqrt{\frac{\Cal{Q}\left(Lw\right)}{\Cal{Q}\left(w\right)}} .$$
 
We define further the least coefficient of expansion by
 
$$\sigma_{\Cal{C}} \left(L\right)= \inf_{w\in int\Cal{C}}
\beta\left(w,L\right).$$
 
Let us note that, for any two monotone maps $L_1$ and $L_2$,
$$
\sigma_{\Cal{C}} \left(L_2L_1\right) \geq
\sigma_{\Cal{C}} \left(L_2\right)\sigma_{\Cal{C}} \left(L_1\right),
$$
i.e., the coefficient of expansion $\sigma_{\Cal{C}}$ is
supermultiplicative.
 
We will omit the index $\Cal{C}$ in $\sigma_{\C}(L)$ when it is
clear what sector
we have in mind.
 
We want to find the value of the expansion coefficient in coordinates.
We will use the fact that this infimum does not change if $L$ is
multiplied on the left or on the right by $\Cal{Q}$-isometries.
So let
$$
L = \left(\matrix A& B\\ C& D\endmatrix\right)
$$
be a monotone matrix. By the factorization \thetag{4.7} $C^*B = PR$
is equal to the product of two
positive semidefinite matrices and so it has only real non-negative
eigenvalues. Let us denote them by $0\leq t_1 \leq \dots \leq t_d$.
The monotone map $L$ is strictly monotone if and only if $t_1 > 0$.
 
\proclaim {Proposition 6.1}
For a monotone map $L$
$$
\sigma \left(L\right) = \sqrt{1+t_1}+\sqrt{t_1}=
\exp \sinh^{-1}\sqrt{t_1},
$$
moreover, if $L$ is strictly monotone
$$
\sigma \left(L\right) = \beta\left(w,L\right)
$$
for some $w\in \text{int }\Cal{C}$.
\endproclaim
 
\demo{Proof}
Let us put
$$
m \left(L\right) = \sqrt{1+t_1}+\sqrt{t_1} =
\min_{1\le i \le d}\left(\sqrt{1+t_i}+\sqrt{t_i}\right) .
$$
First we prove the inequality $\beta\left(w,L\right) \ge m
\left(L\right)$ for $w\in\text{int}\Cal{C}$.
Since both $\beta\left(w,L\right)$ and $m \left(L\right)$ are
continuous functions of $L$ it is sufficient to prove the inequality for
strictly monotone maps only. In view of Proposition 4.9 we can restrict
ourselves to maps $L$ of the form
$$
L = \left(\matrix I&I\\ T& I+T\endmatrix\right)
$$
with diagonal $T$ and $t_1,\dots ,t_d$ on the diagonal.
We compute $\beta(w,\,L)$ directly, for $w=\left(\xi ,\eta\right)$ such that
$\Cal{Q}\left(w\right)=1$
$$
\aligned
\left(\beta\left(w,L\right) \right)^2 & = \sum_{i=1}^{d}\left(t_i \xi_i^2 +
\left(1+2t_i\right)\xi_i\eta_i + \left(1+t_i\right)\eta_i^2\right)\\
&=\sum_{i:\xi_i\eta_i \ge 0}\left(\left(\sqrt{t_i}\xi_i -
\sqrt{1+t_i}\eta_i\right)^2
+\left(\sqrt{1+t_i} + \sqrt{t_i}\right)^2\xi_i\eta_i \right)\\
&+ \sum_{i:\xi_i\eta_i < 0}\left(\left(\sqrt{t_i}\xi_i +
\sqrt{1+t_i}\eta_i\right)^2
+\left(\sqrt{1+t_i} - \sqrt{t_i}\right)^2\xi_i\eta_i \right) \ge\\
&\ge \sum_{i:\xi_i\eta_i \ge 0}\left(\sqrt{1+t_i} +
\sqrt{t_i}\right)^2\xi_i\eta_i +
\sum_{i:\xi_i\eta_i < 0}\left(\sqrt{1+t_i} +
\sqrt{t_i}\right)^{-2}\xi_i\eta_i\ge\\
&\ge\left(1+\delta\right) m\left(L\right)^2-
\delta m\left(L\right)^{-2} \ge m\left(L\right)^2
\endaligned
$$
where
$$
\delta =\left(\sum_{i:\xi_i\eta_i \ge 0}\xi_i\eta_i\right) -1
=\sum_{i:\xi_i\eta_i<0} \xi_i\eta_i\ge 0
$$
and all the inequalities become equalities for
$$
\xi_1 = \left(\frac{1+t_1}{t_1}\right)^{\frac 14},
\eta_1 = \left(\frac{t_1}{1+t_1}\right)^{\frac 14},
\  \xi_i = 0, \eta_i =0,\
i= 2,\dots,d.
$$
 
Thus the Proposition is proven for strictly monotone matrices
and for all monotone matrices we get the inequality $\sigma(L) \geq
m(L)$.
To get the equality $\sigma(L) = m(L)$
for all monotone matrices we proceed as follows. For any
$\epsilon > 0$ we choose a strictly monotone matrix $L_{\epsilon}$ so close to
 the
identity that $m\left(L_{\epsilon}L\right) < m\left(L\right)
+{\epsilon}.$
Since $L_{\epsilon}L$ is strictly monotone and our Proposition has
been proven for strictly monotone matrices there is
$w_{\epsilon}\in \text{int}\Cal{C}$ such that
$$
\beta\left(w_{\epsilon},L_{\epsilon}L\right) =
m(L_{\epsilon}L)= \sigma\left(L_{\epsilon}L\right).
$$
But $
\beta\left(w,L_{\epsilon}L\right) > \beta\left(w,L\right) $ for any $w \in
\text{int}\Cal{C}.$
Hence
$$
m(L) \leq \sigma\left(L\right) \le \beta\left(w_{\epsilon},L\right) <
 \beta\left(w_{\epsilon},L_{\epsilon}L\right) = m\left(L_{\epsilon}L\right) <
m\left(L\right) +{\epsilon}$$
which ends the proof.
\enddemo

For a given sector $\Cal{C} = \Cal{C}(V_1,V_2)$ let $\Cal{C'} =
\Cal{C}(V_2,V_1)$ be the complementary sector. We have
\proclaim{Proposition 6.2}
If $L$ is (strictly) monotone with respect to $\Cal{C}$ then $L^{-1}$
is (strictly) monotone with respect to $\Cal{C'}$ and
$ \sigma_{\Cal{C}} (L) = \sigma_{\Cal{C'}} (L^{-1})$.
\endproclaim
\demo{Proof}
We have that the union
$$
 \Cal{C}\left(V_1,V_2\right) \cup
\text{int}\Cal{C}\left(V_2,V_1\right)
$$
is equal to the whole linear symplectic space $\Cal W$.
Hence if $$L \Cal{C}\left(V_1,V_2\right) \subset \Cal{C}\left(V_1,V_2\right)$$
 then
$$\Cal{C}\left(V_1,V_2\right) \subset L^{-1} \Cal{C}\left(V_1,V_2\right) $$ and
 finally
$$L^{-1}\text{int}\Cal{C}\left(V_2,V_1\right) \subset
 \text{int}\Cal{C}\left(V_2,V_1\right).
$$
The last property is easily seen to be equivalent to the monotonicity of
 $L^{-1}$.
 
To obtain the equality of the coefficient of least expansion we will
use the standard sector and the block description of $L$. Let (see
\thetag{4.7})
$$
L = \left(\matrix A& 0\\ 0& A^{*-1}\endmatrix\right)\left(\matrix I&
0 \\ P &I \endmatrix\right)\left(\matrix I& R \\ 0 &I
\endmatrix\right).
$$
The linear symplectic map
$
\left(\matrix 0 & I \\ -I & 0 \endmatrix\right)
$
takes the standard sector $\Cal{C}$ onto $\Cal{C'}$ and further
$$
L_1 = \left(\matrix 0 & -I \\  I & 0 \endmatrix\right)L^{-1}
\left(\matrix 0 & I \\ -I & 0 \endmatrix\right)
$$
has the same least coefficient of expansion with respect to $\Cal{C}$
as $L^{-1}$ with respect to $\Cal{C'}$.
Since
$$
L^{-1} = \left(\matrix I& -R \\ 0 &I \endmatrix\right)
\left(\matrix I& 0 \\ -P &I \endmatrix\right)
\left(\matrix A^{-1}& 0\\ 0& A^*\endmatrix\right)
$$
we get
$$
L_1 = \left(\matrix I & P \\ R & I+RP \endmatrix\right)
\left(\matrix A^*& 0\\ 0& A^{-1}\endmatrix\right).
$$
Our claim follows now from the formula in Proposition 6.1 and the fact
that $PR$ has the same eigenvalues as $RP$.
\enddemo
 
The next Proposition is a useful addition to the Corollary 5.3.
\proclaim {Proposition 6.3}
For a strictly monotone map $L$
$$
d(LV_1,\,LV_2) =\ln\frac{\sigma(L)^2 +1}{\sigma(L)^2 -1} .
$$
\endproclaim
\demo{Proof}
Since $\Cal Q-\text{isometries}$ preserve the distance between Lagrangian
 subspaces
it follows from Proposition 4.9 that we can restrict our calculations to
$$
L = \left(\matrix I& I\\ T& I+T\endmatrix\right)
$$
with diagonal $T$. By the Definition 5.1 we have
$$
\aligned
d(LV_1,\,LV_2)&=\frac 12 \sup_{0\neq \xi \in \Bbb R^d} \ |\ \ln\langle
 \xi,\,T\xi \rangle - \ln\langle \xi,\,(T+I)\xi \rangle|\\
&= \frac 12\sup_{0\neq \xi \in \Bbb R^d}\ln\frac{\langle \xi,\,(I+T^{-1})\xi
 \rangle}
{\langle \xi,\,\xi \rangle} = \max_{i}\frac{\ln
\left(1+t_i^{-1}\right)}{2}=
\frac{\ln \left(1+t_1^{-1}\right)}{2}
\endaligned
$$
where $t_1\leq t_2  \leq \dots \leq t_d$ are the eigenvalues of $T$.
The desired formula is now obtained by
a straightforward calculation.
\enddemo
 
We introduce now an important  property of a sequence of monotone
maps. Let us consider a sequence of
linear symplectic monotone maps $\{L_i\}_{i=1}^{+\infty}$.
To simplify notation let us put $ L^n = L_n\dots L_1.$
\proclaim{Definition 6.4}
A sequence $\{L_1,L_2,\dots\}$ of monotone maps is called unbounded
if for all
$w \in \text{int}\Cal{C}$ $$\Cal{Q}(L^nw)
\rightarrow +\infty \ \ \text{as} \ \ n \rightarrow +\infty . $$
It is called strictly unbounded if for all
$w \in \Cal{C}, w \neq 0,$
$$\Cal{Q}(L^nw)
\rightarrow +\infty \ \ \text{as} \ \ n \rightarrow +\infty . $$
\endproclaim
 
\proclaim{Theorem 6.5}
A sequence  $\{L_1,L_2,\dots\}$ of maps monotone with respect
to $\Cal{C}$ is unbounded
if and only if
$$
\bigcap_{n = 1}^{+\infty}L_1^{-1}L_2^{-1}\dots L_n^{-1} \C' =
 \text{one Lagrangian subspace}$$
where $\C'$ is the complementary sector.
\endproclaim
\proclaim{Corollary 6.6}
If a sequence of monotone maps $\{L_1,L_2,\dots\}$ is unbounded then
the sequence $\{L_2,L_3,\dots\}$ is also unbounded.
\endproclaim \qed
 
We were not able to find a proof of Corollary 6.6 independent of
Theorem 6.5.
 
\demo{Proof of Theorem 6.5}
We note that $\{L_1,L_2,\dots\}$ is unbounded
if and only if for any strictly monotone $L$ the sequence
$\{L,L_1,L_2,\dots\}$ is unbounded.
 
The next step is to prove that $\{L_1,L_2,\dots\}$ is unbounded if and
only if for every strictly monotone $L$
$$\sigma_{\Cal C}\left(L^nL\right) \rightarrow +\infty \ \ \text{as} \ \ n
\rightarrow +\infty .\tag{6.7} $$
Indeed the last property implies immediately that $\{L,L_1,L_2,\dots\}$ is
unbounded and so, if it holds for all strictly monotone  $L$, then also
$\{L_1,L_2,\dots\}$ is unbounded. To prove the converse we will need
the following well known fact from point set topology:
\proclaim{Lemma}
Let $f_1 \leq f_2\leq \dots ,$ be a nondecreasing sequence of real-valued
continuous functions defined on a compact Hausdorff space $X$.
If for every
$x\in X$
$$
\lim_{n\to +\infty}f_n(x) = +\infty
$$
then
$$
\lim_{n\to +\infty}\inf_{x\in X}f_n(x) = +\infty .
$$
\endproclaim
 
If $\{L_1,L_2,\dots \}$ is unbounded and $L$ is strictly monotone then
we have
$$
\sigma_{\C}\left(L^nL\right) =
\inf_{w\in int\C} {{\sqrt{\Q(L^nLw)}}
\over {\sqrt{\Q(w)}}} \ge \inf_{0\neq w\in \C}
{{\sqrt{\Q(L^nLw)}} \over {\sqrt{\Q(Lw)}}}\
\sigma_{\C}\left(L\right) . $$
 
Applying the Lemma to
$$
f_n(w)= {{\sqrt{\Q(L^nLw)}} \over {\sqrt{\Q(Lw)}}},\ n=1,2,\dots ,
$$
which can be considered as a sequence of functions on the compact space
of rays in $\Cal C$ we obtain \thetag{6.7}.
 
Now we will be proving that \thetag{6.7} is equivalent to
$$\bigcap_{n=1}^{+\infty}L^{-1}L_1^{-1}L_2^{-1}\dots L_n^{-1}\C' =
\text{one Lagrangian subspace}
$$
where $\C'= \C(V_2,\,V_1)$ is the complementary sector. The sectors
$$
\C'_n = L^{-1}L_1^{-1}L_2^{-1}\dots L_n^{-1}\C' =
L^{-1}\left(L^n\right)^{-1}\C' =
\C(L^{-1}\left(L^n\right)^{-1}V_2,\,L^{-1}\left(L^n\right)^{-1}V_1)
$$
$n=1,2,\dots,$ form a nested sequence. We consider the space
$Lag(\C')$ of all Lagrangian subspaces contained strictly in $\C'$
with the metric defined in Section 5. The sequence of subsets
$\widehat{Lag}(\C'_n) \subset Lag(\C'), n= 1,2,\dots,\dots,$ is a nested
sequence of compact subsets. Hence its intersection contains one point
(= Lagrangian subspace) if and only if their diameters converge to zero.
By Corollary 5.3 the diameter of $\widehat{Lag}(\C'_n)$ is equal to
the distance of the Lagrangian subspaces $L^{-1}\left(L^n\right)^{-1}V_2$
and $L^{-1}\left(L^n\right)^{-1}V_1$.
By Proposition 6.3 this distance is equal to
$$
\ln \frac{s_n^2+1}{s_n^2-1}
$$
where $s_n = \sigma_{\C'}\left(L^{-1}(L^n)^{-1}\right).$
But by Proposition 6.2
$$
\sigma_{\C'}\left(L^{-1}(L^n)^{-1}\right) =\sigma_{\C}(L^nL).
$$
This shows that indeed the set
$$
\bigcap_{n=1}^{+\infty}\widehat{Lag}(C'_n)
$$
contains exactly one point if and only if \thetag{6.7} holds.
\enddemo
 
We will use the following characterization of strict unboundedness.
\proclaim{Theorem 6.8}
Let $\{L_i\}_{i=1}^{+\infty}$ be a sequence of linear
symplectic monotone maps. The
following are equivalent.
$$
\aligned
(1)& \  \text { The sequence } \  \{L_i\}_{i=1}^{+\infty} \  \text
{ is strictly unbounded, } \\
(2)& \ \inf_{0\neq w\in \C} {{\sqrt{\Q(L^nw)}} \over {\|w\|}}
\rightarrow +\infty \ \text{as} \  n \rightarrow +\infty ,\\
(3)& \ \sigma(L^n) \rightarrow +\infty \ \text{as} \  n
\rightarrow +\infty ,\\
(4)& \ \text {the sequence } \ \{L_i\}_{i=1}^{+\infty} \ \text
{ is unbounded and}\  L^{n_0}\
\text{is strictly monotone for some}\\
&\ n_0 \geq 1.
\endaligned
$$
\endproclaim
\demo{Proof}
The Lemma from set topology used in the Theorem 7.5 can also be applied
to the sequence of functions
$$
f_n(w)=  {{\sqrt{\Q(L^nw)}} \over {\|w\|}},\ n=1,2,\dots,
$$
to shows that  (1)  $\Rightarrow$  (2). Further (2) $\Rightarrow$ (3)
because
$$
\sigma(L^n) = \inf_{w\in int\C} {{\sqrt{\Q(L^nw)}} \over
{\sqrt{\Q(w)}}} \geq
\inf_{0\neq w \in \C}{{\sqrt{\Q(L^nw)}} \over {\|w\|}}
\inf_{ w \in int\C}{{\|w\|}\over {\sqrt{\Q(w)}}}.
$$
The implication (3) $\Rightarrow$ (4) is obvious ($\sigma(L^n) > 1$ if
and only if $L^n$ is strictly monotone, cf. Proposition 6.1).
Finally let the sequence $\{L_i\}_{i=1}^{+\infty}$  be unbounded
and  $L^{n_0}$ be strictly monotone. By Corollary 6.6 also the sequence
$\{L_{n_0+1},L_{n_0+2},\dots\}$ is unbounded. It follows that
$\{L_i\}_{i=1}^{+\infty}$ is strictly unbounded.
\enddemo
 
The following example plays a role in the study of special Hamiltonian
systems.
\vskip.5cm
\centerline{\bf  Example.}
\vskip.2cm
Let
$$
L_n = \left(\matrix A_n& 0\\ 0& A_n^{*-1}\endmatrix\right)
\left(\matrix I& 0 \\ P_n &I \endmatrix\right)
\left(\matrix I& R_n \\ 0 &I \endmatrix\right),
$$
$n=1,2,\dots ,$ be a sequence of monotone symplectic matrices with nonexpanding
 $A_n$,
i.e., $\|A_n\xi\| \leq \|\xi\|$ for all $\xi$. We assume further that
the symmetric matrices $R_n$ satisfy
$$
\tau_n' I \geq R_n \geq \tau_n I \ \ \text{and}\ \
\frac{\tau_n'}{\tau_n}
\leq C
$$
for some positive constants $C$ and $\tau_n,\tau_n',\ n= 1,2,\dots .$
We do not make any assumptions about
$P_n$ (beyond $P_n \geq 0$ which is forced by the monotonicity  of
$L_n$).
Note that if a symmetric matrix $R$ satisfies $\tau I \leq R \leq \tau '
I$
then $\tau \|\eta \| \leq \|R \eta \|
\leq \tau '\|\eta \|$. Indeed
$$
\left<R\eta,\,R\eta\right> =
 \frac{\left<RR^{\frac12}\eta,\,R^{\frac12}\eta\right>}
{\left<R^{\frac12}\eta,\,R^{\frac12}\eta\right>}\ \left<R\eta,\,\eta\right>
$$
which yields the estimate.
\proclaim{Proposition 6.9}
If $\sum_{n=1}^{+\infty}\tau_n = +\infty$ then the sequence
$\{L_1,L_2,\dots\}$ is unbounded.
\endproclaim
\demo{Proof}
Let $w_1 = (\xi_1,\,\eta_1) \in int\C$ and $w_{n+1} =
(\xi_{n+1},\,\eta_{n+1}) =
L_nw_n, n=1,2,\dots$.
Our goal is to show that
$$
q_n = \Q (w_n) \rightarrow +\infty\ \ \text{as}\ \ n \rightarrow +\infty.
$$
 
We have $\xi_{n+1} = A_n\left( \xi_n + R_n \eta_n \right) $ so that
$$
\|\xi_{n+1}\| \leq \|\xi_n\| + \|R_n \eta_n \| \leq \|\xi_n\| +
\tau_n'\|\eta_n \| \leq \|\xi_1\| + \sum_{i=1}^{n}\tau_i'\|\eta_i \|.
\tag{6.10}
$$
 
At the same time $ q_n = \left<\xi_n,\,\eta_n\right> \leq \|\xi_n\|\|\eta_n \|$
so that
$$
\|\eta_n \| \geq \frac{q_n}{\|\xi_n\|}
\tag{6.11}
$$
and hence (see also the proof of Theorem 4.4)
$$
q_{n+1} \geq q_n + \left<R_n\eta_n,\,\eta_n\right> \geq
q_n + \tau_n \|\eta_n \|^2 \geq q_n + \tau_n \|\eta_n \|
\frac{q_n}{\|\xi_n\|}.
$$
Using \thetag{6.10} we obtain from the last inequality
$$
\frac{q_{n+1}}{q_n} \geq 1 + \frac{\tau_n \|\eta_n \|}
{\|\xi_1\| + \sum_{i=1}^{n-1}\tau_i'\|\eta_i \|} \geq
1 + \frac1C\frac{\tau_n' \|\eta_n \|}
{\|\xi_1\| + \sum_{i=1}^{n-1}\tau_i'\|\eta_i \|}.
\tag{6.12}
$$
If $\sum_{i=1}^{+\infty}\tau_i'\|\eta_i \| < +\infty $ then by
\thetag{6.10} the sequence $\|\xi_n\|$ is
bounded from above and hence by \thetag{6.11} the sequence $\|\eta_n \|$
is bounded away from zero which is a contradiction
(in view of $\sum_{i=1}^{+\infty}\tau_i' = +\infty $).
 
Hence
$$
\sum_{i=1}^{+\infty}\tau_i'\|\eta_i \| = +\infty.
$$
 Now the claim follows from
\thetag{6.12} and the following
\proclaim{Lemma 6.13}
For a sequence of positive numbers $a_0,a_1,\dots ,$
if
$$
\sum_{n=1}^{+\infty}a_n = +\infty \ \ \text{then } \ \
 \sum_{n=1}^{+\infty}\frac{a_n}
{\sum_{i=0}^{n-1}a_i} = +\infty.
$$
\endproclaim
\demo{Proof of the Lemma}
We have for $1 \leq k \leq l$
$$
\sum_{n=k}^{l}\frac{a_n}{\sum_{i=0}^{n-1}a_i} \geq \frac{\sum_{n=k}^{l}a_n}
{\sum_{n=0}^{l}a_n} \rightarrow 1\ \ \text{as} \ \ l  \rightarrow +\infty .
$$
\enddemo
\enddemo

 
\vskip.7cm
\subhead \S 7. PROPERTIES OF THE SYSTEM AND THE FORMULATION OF THE
RESULTS
\endsubhead
\vskip.4cm
 
In this section we define rigorously the class of systems to
which the present paper applies.
We divide the conditions that the systems must satisfy into several
groups. The multitude of
conditions is justified by the fact that we want to include
discontinuous systems (there is only one way to be continuous but many
ways to be discontinuous !). In the case  of a
symplectomorphism of a compact symplectic manifold most of these
conditions are vacuous. Because of that we will single out this case
and we will refer to
it as the smooth case. The bulk of our effort is devoted to
the discontinuous case.
\medskip
{\bf A. The phase space.}\par
 
In the smooth case the phase space $\Cal M$ is a smooth compact symplectic
manifold.
 
In the discontinuous case it is a disjoint union of nice subsets of
the linear symplectic space.
More precisely, let us consider the standard linear symplectic
space $\Cal W = \Re^d \times \Re^d$ equipped with a Riemannian metric
uniformly equivalent to the standard Euclidean scalar product
and which defines the same volume element (measure) $\mu$. The measure
$\mu$ is also equal to the symplectic volume element.
 
By a submanifold of $\Cal W$ we mean an embedded submanifold
of $\Cal W$. Further we define a piece of a submanifold $\Cal S$ to be
a compact subset of $\Cal S$ which is the closure of its interior (in
the relative topology of the submanifold $\Cal S$).  A piece $X$ of a
submanifold has a well defined boundary which we will denote by
$\partial X$ (it is the set of boundary points with respect to the
relative topology of the submanifold). Notice that at every point of a piece
of a submanifold, including a boundary point, we have a well defined
tangent subspace.
 
 A submanifold carries the measure
defined by the Riemannian volume element, for this measure the  boundary
of a piece of a submanifold is not necessarily of zero measure.
 
The phase space is made up of pieces of $\Cal W$ which have regular
boundaries in the sense of the following definition.
 
\proclaim {Definition 7.1}
A compact subset $X \subset \Cal W$ is called regular if
it is a finite union of pieces $X_i,i=1,\dots,k,$
of $2d-1$-dimensional submanifolds
$$
X = X_1 \cup \dots \cup X_k.
$$
The pieces overlap at most on their boundaries, i.e.,
$$
X_i \cap X_j \subset \partial X_i \cup \partial X_j,\ i,j = 1,\dots k;
$$
and the boundary $\partial X_i$ of each piece
$X_i,\ i=1,\dots k,$ is a finite
union of compact subsets of $2d-2$-dimensional submanifolds.
\endproclaim
 
To picture such sets one can think of the boundary of a $2d$-dimensional
cube. The faces are pieces of $2d-1$-dimensional submanifolds
and they clearly
overlap only at their boundaries. The boundary of each face is a
union of pieces of $2d-2$ dimensional submanifolds (actually it is
a union of $2d-2$ dimensional cubes). Let us stress that in the
definition of a regular set we do
not impose any requirements on the $2d-2$ dimensional subsets in the
boundary.
Due to the generality of the definition one cannot even claim that
the union of two regular sets is regular.
 
As a consequence of Definition 7.1 the natural measures on
the pieces $X_i,i=1,\dots,k,$ of any regular subset $X$ can be
concocted to give a well defined measure $\mu_X$ on $X$ (the $2d-1$
dimensional Riemannian volume). It is so
because the boundaries of the pieces being themselves finite unions of
subsets of
submanifolds of lower dimension have zero measure. Hence if we put
$$
\partial X = \bigcup_{i=1}^k\partial X_i,
$$
then
$$
\mu_X\left(\partial X\right) = 0.
\tag {7.2}
$$
Moreover, by the regularity of the measure $\mu_{X}$,
it follows from \thetag{7.2} that, if we denote by $(\partial X)^\delta$
the $\delta$-neighborhood of $\partial X$ in $X$, then
$$
\lim_{\delta \to 0} \mu_X\left((\partial X)^\delta\right) = 0.
\tag {7.3}
$$
 
Further we have the following Proposition.
 
\proclaim{Proposition 7.4}
For a subset $Y$ of $X\subset \Cal W$  let the $\delta$-neighborhood of
$Y$ in $\Cal W$ be denoted by $Y^\delta$, i.e.,
$$
Y^\delta = \{x \in \Cal W\,|\, \text{d}(x,Y) \leq \delta \}.
$$
If $X$ is a regular ($2d-1$-dimensional) subset of $\Cal W$
and $Y \subset X$ is {\bf closed} then
$$
\lim_{\delta \rightarrow 0} {{\mu(Y^\delta)}\over 2\delta } =
\mu_X(Y).
$$
\endproclaim
Although Proposition 7.4 holds as we formulated it, we will
use only the weaker property
$$
\limsup_{\delta \rightarrow 0} {{\mu(Y^\delta)}\over \delta } \leq
const\mu_X(Y).
\tag 7.5
$$
We leave the proof of the Proposition or of the easier property
\thetag{7.5} to the reader.
 
\proclaim {Definition 7.6}
A compact subset $\Cal M \subset \Cal W$ is called a symplectic box if
the boundary $\partial \Cal M$ of $\Cal M$ is a regular subset of
$\Cal W$ and the interior $int\Cal M$ of $\Cal M$ is connected and
dense in  $\Cal M$.
\endproclaim
 
We can now formulate the requirements on the phase space of
a discontinuous system.\par
{\it The phase space of our system is a finite disjoint union of
symplectic boxes}.
 
To simplify notation we assume that the phase space consists of just one
symplectic box $\Cal M$. It will be quite obvious how to generalize the
subsequent formulations to the case of several symplectic boxes.
\medskip
{\bf B. The map $T$ (the dynamical system).}
 
In the smooth case the map $T$ is a symplectomorphism
$T: \Cal M \to \Cal M$.
 
In the discontinuous case we assume that the symplectic box $\Cal M$ is
partitioned in two ways into unions of equal number of symplectic
boxes
$$
\Cal M = \Cal M_1^+ \cup \dots \cup \Cal M_m^+
= \Cal M_1^- \cup \dots \cup \Cal M_m^-.
$$
Two boxes of one partition can overlap at most on their boundaries, i.e.,
$$
\Cal M_i^{\pm}\cap \Cal M_j^{\pm} \subset \partial \Cal M_i^{\pm}\cap
\partial \Cal M_j^{\pm} ,\ \ i,j = 1,\dots ,m.
$$

The map $T$ is defined separately on each of the symplectic boxes
$\Cal M_i^+, \ i =1,\dots,m$. It is a symplectomorphism of the interior of
each $\Cal M_i^+$ onto the interior $\Cal M_i^-, \ i =1,\dots,m$
and  a homomorphism of
$\Cal M_i^+$ onto  $\Cal M_i^-, \ i =1,\dots,m$. We assume that the
derivative $DT$ is well behaved near the boundaries of
the symplectic boxes. Namely, we assume that it satisfies the
Katok-Strelcyn conditions so that we can apply their results \cite {K-S}
on the existence of the foliation in (un)stable manifolds and
its absolute continuity.
 
We will say that $T$ is a (discontinuous) symplectic map of $\Cal M$.
Formally $T$ is not well defined on the set of points which belong
to the boundaries of several plus-boxes: it has several values.
We adopt the convention that the image of a subset of $\Cal M$ under $T$
contains all such values.
 
Let us introduce the singularity sets $\R^{+}$ and $\R^-$.
$$
\R^{\pm} = \{ p\in \Cal M \;|\; p
\text { belongs to at least two of the boxes }\Cal M_i^{\pm}, i =1,\dots,m\}.
$$
 
The plus-singularity set $\R^+$ is a closed subset and $T$ is
continuous on its complement. Similarly $T^{-1}$ is continuous on the
complement of $\R^-$. Note that most of the points in the
boundary $\partial \Cal M$
of $\Cal M$ do not belong to $\R^-$ or $\R^+$.
 
We have that $\R^+ \cup \partial \Cal M$ is the union of all the boundaries
of the plus-boxes and $\R^- \cup \partial \Cal M$ is the union of all the
boundaries of the minus-boxes, i.e.,
$$
\R^{\pm} \cup \partial \Cal M = \bigcup_{i=1}^m\partial \Cal M_i^{\pm}.
$$

Note that most of the points in the boundary $\partial \Cal M$ of $\Cal M$
do not belong to $\Cal S^-$ or $\Cal S^+$. 
We assume that the singularity sets $\R^{\pm}$
and the union of boundaries $\bigcup_{i=1}^m\partial \Cal M_i^{\pm}$
are regular sets.
 
An important role in our discussion will be played by the singularity
sets of the higher iterates of $T$.
We define for $n\geq 1$
$$
\R^+_n = \R^+ \cup T^{-1}\R^+ \cup \dots
\cup T^{-n+1}\R^+.
$$
and
$$
\R^-_n = \R^- \cup T\R^- \cup \dots
\cup T^{n-1}\R^-.
$$
We have that $T^n$ is continuous on the complement of $\R^+_n$ and
$T^{-n}$ is continuous on the complement of $\R^-_n$.
 
\proclaim{Regularity of singularity sets}
We assume that for every $n \geq 1$ both $\R^+_n$ and $\R^-_n$ are
regular.
\endproclaim
 
We will formulate, in Lemma 7.7, an abstract condition on the first
power of $T$ alone that guarantees the regularity of the singularity
sets but it requires that the map is a diffeomorphism on every
symplectic box up to and including its boundary i.e., it can be extended
to a diffeomorphism of an open neighborhood of $\Cal M_i^+$ onto an
open neighborhood of $\Cal M_i^-,\ i=1,\dots,m$.
 
Hence it is very appealing to restrict the discussion to such maps.
Unfortunately, such a restriction
would leave out important examples: billiard systems where the
derivative may blow up at the boundary.
The  conditions in the work of Katok and Strelcyn \cite{K-S} were
tailored for such systems.
 
Nevertheless the reader is invited  to be generous with the
restrictions on the regularity of $T$, this will make it easier to
follow the main line of the argument.
\medskip
{\bf C. Monotonicity of $T$.}
\par
In the smooth case we assume that two continuous bundles of
transversal Lagrangian subspaces are chosen in an open
subset $U \subset \Cal M$ ($U$ is not necessarily dense). We denote
them by $\{V_1(p)\}_{p\in \U}$ and $\{V_2(p)\}_{p\in \U}$
respectively.
 
In the discontinuous case we assume that  two continuous bundles of
transversal Lagrangian subspaces are chosen in the interior of the
symplectic box $\Cal M$.
Their limits (if they exist at all)  at the boundary  $\partial\M$
are allowed to have nonzero intersection.
 
We consider the bundle of sectors
(see Definition 4.1) defined by these Lagrangian subspaces
$$
\Cal C(p)=\Cal C(V_1(p),\,V_2(p)).
$$
Let
$$
\Cal C'(p)=\Cal C(V_2(p),\,V_1(p))
$$
be the complementary sector.

We require that the derivative of the map and its iterates, where defined, is
monotone, if only monotonicity is well defined (cf. Definition 4.3).
 
More precisely,
in the smooth case we require that, if  $p\in\U$ and  $T^kp\in\U$ for
$k\geq 1$, then
$$
D_pT^k\C(p)\subset\C(T^kp).
$$
 
In the discontinuous case we assume that
$$
D_pT \Cal C(p) \subset \Cal C(Tp)
$$
for points $p$ in the interior of every symplectic boxes
$\Cal M_i^+, i = 1,\dots,m$ .
 
We call a point $p \in int\Cal M$ ($p\in U$ in the smooth case) {\it
strictly monotone in the future}  if there is $n \geq 1$ such that $D_pT^{n}$
is defined and it is strictly monotone
( in the smooth case we require naturally
that $T^np\in U$), i.e.,
$$
D_pT^{n} \Cal C(p) \subset int\Cal C(T^{n}p) \cup \{0\}.
$$
 
Similarly a point $p$ is called {\it strictly monotone in the past}
if there is $n \geq 1$ such that $D_pT^{-n}$
is strictly monotone with respect to the complementary sectors, i.e.,
$$
D_pT^{-n} \C'(p) \subset int\C'(T^{-n}p) \cup \{0\}.
$$
 
It is clear that if $p$ is strictly monotone  in the future then its
preimages are also strictly monotone in the future.
By Proposition 6.2 we also have that if $p$ is strictly monotone
in the future then
there is $n \geq 1$ such that $T^np$ is strictly monotone in the past.
\proclaim{Strict monotonicity almost everywhere}
We assume that almost all points in $\Cal M$ ($U$ in the smooth case)
are strictly monotone.
\endproclaim
This property implies that all Lyapunov exponents are non-zero almost
everywhere
in $\Cal M$ (in $U$ in the smooth case). The proof of this fact is
quite simple and can be found in \cite {W1}. It will also follow
easily from our Proposition 8.4. Thus by the work of Pesin \cite {P} in
the smooth case and of Katok and Strelcyn \cite {K-S} in the discontinuous
case through almost every point there are local stable and unstable
manifolds of dimension $d$ and the foliations into these manifolds are
absolutely continuous.
 
The sectors $\C(p)$ contain the unstable Lagrangian subspaces (tangent
to the unstable manifolds) and the complementary sectors $\C'(p)$ contain the
stable Lagrangian subspaces (tangent to the stable manifolds).
The sectors can be viewed as a priori approximations to the unstable
and stable subspaces.
We will refer to the sectors as unstable sector  and stable
sector respectively.

This ends the list of required properties for the smooth case.
The last three properties of our system are introduced only for
the discontinuous case.
\medskip
{\bf D. Alignment of Singularity sets}
\par
For a codimension one subspace in a linear symplectic space its
characteristic line is, by definition, the skeworthogonal complement
(which is a one dimensional subspace).
 
\proclaim {Proper alignment of $\R^-$ and $\R^+$}
We assume that  the tangent subspace of $\R^-$ at any $p  \in \R^-$
has the characteristic
line contained strictly in the sector $\C(p)$
and  that  the tangent subspace of $\R^+$ at any $p  \in \R^+$
has the characteristic line contained strictly in the complementary
sector $\C'(p)$. We say that the
singularity sets $\R^-$ and $\R^+$ are properly aligned.
\endproclaim
Let us note that if a point in $\R^\pm$ belongs to several pieces of
submanifolds then we require that the tangent subspaces to {\it all}
of these pieces have characteristic lines in the interior of the sector.
\par
It will be clear from the way in which the proper alignment of
singularity sets is used in Section 12 that it is sufficient to assume
that there is $N$ such that $T^N\R^-$ and $T^{-N}\R^+$ are properly
aligned. We will show, in section 14, that for the system of
falling balls even this
weaker property fails. Hence the study of ergodicity of this system
would require some further relaxation of this property.
\par
Let us note that it is helpful in establishing the regularity of
singularity sets $\R^{\pm}_n$ if the boundaries of $\Cal M$
have tangent subspaces characteristic lines
contained in the boundary of the sectors $\C(p)$. It is so in some
examples. More precisely we have the following lemma.
\proclaim {Lemma 7.7}
If the map $T$ is a  diffeomorphism
up to and including the boundaries of the symplectic boxes
$\Cal M^+_1,\dots,\Cal M^+_m$, satisfies properties C, D and the boundary
$\partial\Cal M$ of $\Cal M$  has all the tangent subspaces
with characteristic lines contained
in the boundary of the sectors then the sets $\R^{\pm}_n, n\geq 1,$
are regular (i.e. the property {\bf B} is automatically verified).
\endproclaim
\demo{Proof}
Let us recall that, by assumption, $\R^-$ and
$\bigcup_{i=1}^m\partial\Cal M^+_i$ are properly aligned
regular subsets.
Further the intersection of any properly aligned
regular subset $X$ (the characteristic lines
of its tangent subspaces are contained strictly in the unstable sector
$\C$) with any of the symplectic boxes
$\Cal M^+_1,\dots,\Cal M^+_m$ is a regular subset. Indeed
let $X_1,\dots,X_p$ be the
pieces of  $2d-1$ dimensional manifolds which make up $X$
($X=\bigcup_{i=1}^p X_i$) and $Y_1,\dots,Y_q$ be the
pieces of  $2d-1$ dimensional manifolds which make up the boundary of
say $\Cal M^+_1$ ($\partial\Cal M^+_1=\bigcup_{j=1}^q Y_j$).
By the proper alignment
of the pieces we can assume that any $X_i$ and any $Y_j$ are pieces of
transversal submanifold. Hence the intersection of the submanifolds
is a submanifold of dimension $2d-2$, and therefore $X_i\cap Y_j$ are
disjoint pieces of $2d-2$-dimensional manifolds (allowed to intersect
only at the boundary).
It follows that the intersection of $X_i$ with $\Cal M^+_1$ is
a  piece of the $2d-1$ dimensional manifold
and also a regular subset. The same can be repeated for the other
symplectic boxes $\Cal M^+_2,\dots,\Cal M^+_m$.
 
Moreover we have that any
$(X_i\cap \Cal M^+_1)\cup\partial \Cal M^+_1,\ i = 1,\dots,p,$
 is a regular subset and further
$(X\cap \Cal M^+_1)\cup\partial \Cal M^+_1$
is a regular subset. It follows that
$T\left((X\cap\M^+_1)\cup\partial\M^+_1\right)$
$=(TX\cap \Cal M^-_1)\cup\partial \Cal M^-_1$
is a regular subset and after repeating the argument for the other
symplectic boxes we get that for any regular and properly aligned
subset $X$
$
TX\cup\bigcup_{i=1}^m\partial \Cal M^-_i
$
and therefore
$TX\cup \R^-$
are regular properly aligned subsets.
 
Now the proof can clearly be completed by induction since
$$
\R^-_{n+1} = T\R^-_n\cup \R^-.
$$
The argument for $\R^+$ is completely analogous.
\enddemo
\par
The last two properties are rather technical. They are used
only in Section 14 in the proof of the `tail bound'.
It remains an open question if one can do without them.
\medskip
{\bf E. Noncontraction property.}
\par
There is a constant $a, \ 0 < a \leq 1,$ such that for every $n \geq 1$
and for every  $p \in \Cal M \setminus \R^+_n$
$$
\|D_pT^nv\| \geq a \|v\|
$$
for every vector $v$ in the sector $\C(p)$.
\par
Notably the above condition holds in all the examples to which the other
conditions apply (see \S 14), apart from the case of semi-dispersing
billiards in more then two dimensions
(the case from which this type of strategy originated).
In fact, through a tangent collision a
vector in the unstable direction can shrink by an
arbitrary amount. Instead of the present condition the original article of
Chernov-Sinai \cite{CS} was taking advantage of a special property of
semi-dispersing billiard. Namely the existence of a semi-norm
(the configuration norm) that is increased by the dynamics for vectors
in the unstable direction. Moreover, such norm is well aligned with respect
to the singularity manifolds and with respect to the cone bundle: on the
one hand a $\delta$ neighborhood of the singularity in this semi-norm is
of measure $\Cal O(\delta)$, on the other hand the hyperplane of vectors
on which the seminorm has value zero is not contained in the interior of
the cone (note that this two requirement, together with the requirement of
the proper alignment of the singularities, imply that the singularity
manifold is aligned with the boundary of the cone). It would be possible
to generalize such setting and use the generalization of these properties
instead of the non-contraction property. The bold reader can see
how it would be possible to adapt \S 13 to this setting. We choose not
to do this explicitly for reasons of clarity and also because we do not
know of any example (apart from semi-dispersing billiards) to which such
alternative condition could apply.
 
\medskip
{\bf F. Sinai - Chernov Ansatz.}
\par
This is a property pertaining the derivatives of the iterates of $T$ on the
singularity set itself, of $T^{-1}$ on $\R^+$ and of $T$ on $\Cal
R^{-}$. Namely, we require that, for almost every point in $\Cal
R^-$ with respect to the measure $\mu_{\R}$ ($\mu_{\R}$ is the
$2d-1$ dimensional Riemannian volume on $R^- \cup R^+$),
all iterates of $T$ are
differentiable and for almost every point in $\R^+$ all iterates
of $T^{-1}$ are differentiable. Note that the last requirement holds
automatically under the assumptions of Lemma 7.7. Moreover,
 
\medbreak{\sl
we assume that for almost every point $p \in \R^{-}$ with respect to
the measure $\mu_{\R}$, the sequence of derivatives
$\{D_{T^np}T\}_{n\geq 0}$
is strictly unbounded (cf. Definition 6.4).
Analogous property must hold for $\R^+$ and $T^{-1}$.}
\medbreak
By Theorem 6.8 the forward part of Sinai - Chernov Ansatz is equivalent to
the following property. For almost every point
$p \in \R^{-}$ with respect to the measure $\mu_{\R}$
$$
\lim_{n\to +\infty} \sigma(D_pT^n)  = +\infty,
$$
where  the coefficient $\sigma$ is defined at the beginning of
Section 6.
\par
In several examples unboundedness holds for all orbits
by virtue of Proposition 6.9 but strict monotonicity is hard to establish.

We have completed the formulation of the conditions. Under these
conditions we will prove the following two theorems.
 
\proclaim{Main Theorem (Smooth case)}
For any $n\geq 1$ and any $p \in U$ such that $T^np\in U$ and
$\sigma(D_pT^n) > 1$
(i.e., $p$ is strictly
 monotone) there is a neighborhood of $p$
which is contained in one ergodic component of $T$.
\endproclaim
 
It follows from this theorem that if $U$ is connected and  every
point in it is strictly monotone then
$\bigcup_{i=-\infty}^{+\infty}T^iU$ belongs to one
ergodic component. Such a theorem was first proven by
Burns and Gerber \cite {BG}
for flows in dimension 3. It was later generalized by Katok \cite {K}
to arbitrary dimension and recently also
to a non-symplectic framework \cite{K1}.
Our proof is a byproduct of the preparatory steps
in the proof of the following
 
\proclaim{Main Theorem (Discontinuous case)}
For any $n\geq 1$ and for any $p \in \Cal M \setminus \R^+_n$ such that
$\sigma(D_pT^n) > 3$ there is a neighborhood of $p$ which is contained
in one ergodic component of $T$.
\endproclaim
 
Let us note that the conditions of the last theorem are satisfied
for almost all points  $p\in \Cal M$. Indeed
let
$$
\Cal M_{n,\epsilon} = \{p\in\M\;|\;\sigma(D_pT^n) > \epsilon \}.
$$
Since almost all points are strictly monotone, then
$$
\bigcup_{n=1}^{+\infty}\bigcup_{\epsilon>0}
\Cal M_{n,\epsilon}
$$
has full measure. By the Poincare Recurrence Theorem and the
supermultiplicativity of the coefficient $\sigma$ we conclude that
$$
\bigcup_{n=1}^{+\infty}
\Cal M_{n,3}
$$
has also full measure.
 
Hence the theorem implies in particular that all ergodic components
are essentially open. The theorem allows also to go further since
 we assume that only finitely many iterates of $T$ are differentiable
at $p$ so
that we can apply it to orbits that end up on the singularity sets
both in the future and in the past (e.g. $p\in \R^-$ and $T^np \in
\R^+$).  We need though a specific amount of hyperbolicity on this
finite orbit ($\sigma(D_pT^n) > 3$); note that in the smooth case any
amount of hyperbolicity ($\sigma(D_pT^n) > 1$) is sufficient.
 
This theorem gives a fairly explicit description of points which can
lie in the boundary of an ergodic component. By checking that there
are only few such points (e.g. that they form a set of codimension 2)
one may be able to conclude that a given system is ergodic.
 
Although the techniques used in the proof make it unavoidable to
require more hyperbolicity in the non-smooth case,
we do not know of any examples of non-ergodic systems satisfying all
the conditions above where some points on the boundaries of two
ergodic components are strictly monotone, i.e.,
$\sigma(D_pT^n) > 1$ for some $n\geq 1$.
\par
In all the examples that we know, any point with an infinite orbit
(in the future or in the past) has the unbounded sequence of derivatives
(in the sense of Definition 6.4). In such case, it follows from Theorem
6.8 that for any strictly monotone point with the infinite orbit in the future
the condition $\sigma(D_pT^n)>3$ is satisfied automatically, if
only $n$ is sufficiently large.
\par
There is no need to formulate the Main Theorem separately for a
point $p$ which has only the backward orbit $(p\in\Cal S^+)$. 
We can simply apply the theorem to $T^{-n}p$ (one can appreciate
now the convenience of Proposition 6.2).

\topinsert
\vskip 3in
\hsize=4.5in
\noindent{\bf Figure 5} The Baker Map and the Modified Baker Map.
\endinsert

\topinsert
\vskip 3in
\hsize=4.5in
\noindent{\bf Figure 6} The discontinuity lines of the Modified Baker Map.
\endinsert

Let us finish this Section with an example where the role of
the proper alignment of singularities is exposed.
The well known Baker's Transformation maps the unit square
as shown in Fig.5a and it is ergodic.
Let us consider a variation of this construction where the square
is stretched and squeezed as before but now the middle one half is
left at the bottom and the quarters on the left and right are
translated to the top as shown in Fig.5b.
This time the map $T$ is not ergodic. The ergodic components are
separated by the dotted line although for any point $p$ on the dotted
line we have that
$$
\sigma(D_pT^2) = 4.
$$
Of all the conditions formulated in this Section only the proper
alignment of singularity sets is violated; namely part of $\R^-$
has stable (vertical) direction (all of $\R^+$ has stable
direction which is fine), see Fig.6 where $\R^\pm$ are
indicated by bold lines. For the standard Baker's transformation
the condition of the proper alignment is clearly satisfied.

 
\vskip.7cm
\subhead \S 8. CONSTRUCTION OF THE NEIGHBORHOOD AND THE COORDINATE SYSTEM
\endsubhead
\vskip.4cm
 
We will construct a convenient coordinate system in a neighborhood of
a strictly monotone point $p\in\Cal M $. There are two cases:
strict monotonicity in the past and strict monotonicity in the future
but they are completely symmetric. Therefore, we will discuss
only one of them. Namely we assume that there is $N \ge 1$ such that
$$
\aligned
i)&\qquad
T^{-N} \text{ \ is differentiable at \ }  p\;:\;
p \not\in \R^-_N\cup\partial \M,
\text{\it\ (discontinuous case)}\\
&\qquad T^{-N}p \in U, \text{\it\ (smooth case)}\\
ii)&\qquad D_pT^{-N} \text{ \ is strictly monotone.}
\endaligned
\tag 8.1
$$
We will find a neighborhood $\Cal U(p)$ in which there is an
abundance of ``long'' stable and unstable manifolds. Let us emphasize
that we have assumed only that $p$ (and its $N$ preimages)
does not belong to $\R^-$ but it may very well belong to $\R^+$.
Such a level of generality is crucial in obtaining local ergodicity also for
points in the singularity sets $\R^\pm$.
\par
Our first requirement on the neighborhood is that $T^{-N}$ is a 
diffeomorphism of $\Cal U(p)$ onto a neighborhood of $\bar p =
T^{-N}p$ (and in the smooth case both neighborhoods are contained in $U$).
 
By the Darboux theorem a symplectic manifold looks locally like a piece
of the standard linear symplectic space. Hence reducing
$\Cal U(p)$ further, if necessary, we can identify it with a
neighborhood $\Cal U$ of the standard linear symplectic space
$\Re^{d}\times\Re^{d}$
$$
\Cal U = \Cal U_a =  \Cal V_a \times \Cal V_a,
$$
where
$$
\Cal V_a = \{x = (x^1,\dots,x^d)\in\Re^{d}\;|\;
|x^i| <  a, i = 1,\dots,d  \}.
$$
(In the discontinuous case we have assumed from the very beginning that a
symplectic box is a subset in $\Re^d \times \Re^d$).
We assume that the point $p$ becomes the zero point and
the symplectic structure is the standard one. In particular
all the tangent spaces in $\Cal U(p)$ can be identified with
$\Re^{d}\times\Re^{d}$.
The choice of a cube for the shape of the neighborhood is important only
for some of the arguments in Section 11 otherwise we want to stress that
our neighborhood $\Cal U$ is the cartesian product of neighborhoods
$\Cal V_a$ in the $d$-dimensional linear space and we will not use any special
directions there.
 
Let us further introduce for any positive $\rho $  the following sectors
in the tangent space of $\Cal U$.
$$
\C_\rho = \{(\xi,\,\eta)\in\Re^{d}\times\Re^{d}\;|\;
\|\eta\| \leq \rho\|\xi\| \}
$$
and the complementary sector
$$
\C'_\rho = \{(\xi,\,\eta)\in\Re^{d}\times\Re^{d}\;|\;
\|\xi\| \leq \rho^{-1}\|\eta\| \}.
$$

By the assumption  \thetag {8.1} the sector
$D_{\bar p}T^{N}\C(\bar p)$ is strictly
inside the sector $\C(p)$. We change coordinates in $\Cal U$ in such a way
that for some $\tilde{\rho} < 1$
$$
\C'(p) = \C'_{{\tilde{\rho}}^{-1}}
$$
and
$$
D_{\bar p}T^{N}\C(\bar p) \subset \C_{\tilde{\rho}}.
$$
By Propositions 5.11, 6.2 and 6.3 this can be done with $\tilde{\rho} =
(\sigma(D_{\bar p}T^{N}))^{-1}$ .
  
We pick $\rho,\ \tilde{\rho} < \rho <1$. By the continuity 
of the sector bundle $\C(z),
z\in \Cal U$, and of the derivative $ D_{y}T^{N}, y \in T^{-N}\Cal U $,
if we  reduce the size of $\Cal U$ appropriately,
we can achieve that for any $z \in \Cal U$ (see Figure 7)
$$
\C'(z) \subset \C'_{\rho^{-1}}
\tag {8.2}
$$
and for any $y \in T^{-N}\Cal U$
$$
D_{y}T^{N}\C(y) \subset \C_\rho .
\tag {8.3}
$$
 
The properties \thetag {8.2} and \thetag {8.3} seem to be asymmetric
in time, i.e., $T$ plays in them a different role than $T^{-1}$.
Nevertheless we can obtain from them the following
fundamental Proposition which is perfectly symmetric in time.
 
\topinsert
\vskip 3in
\hsize=4.5in
\noindent{\bf Figure 7} The cones at $\Cal T_z\M$.
\endinsert

We will say that a point $z \in \Cal U$ has $k$ spaced returns in a given
time interval if there are $k$ moments of time in this interval
$$
i_1 < i_2 < \dots < i_k
$$
at which $z$ visits $\Cal U$, i.e.,
$$
T^{i_j}z \in \Cal U \text {\ for \ }j=1,\dots,k,
$$
and the visits are spaced by at least time $N$, i.e.,
$$
i_{j+1} - i_j \geq N \text {\ for \ }j=1,\dots,k-1.
$$
 
\proclaim{Proposition 8.4}
If $T^n$ is differentiable at $z\in \Cal U$ for $n \geq N$ and
$z' = T^nz\in \Cal U$  then
$$
D_zT^n \C_{\rho^{-1}} \subset \C_\rho
\tag 8.5u
$$
and
$$
D_{z'}T^{-n} \C'_\rho \subset \C'_{\rho^{-1}}.
\tag 8.5s
$$
 
Moreover for $ (\xi',\,\eta') = D_zT^n(\xi,\,\eta) $
if $(\xi,\,\eta) \in \C_\rho$ then
$$
\|\xi'\| \geq b\rho^{-k}\|\xi\|
\tag 8.6u
$$
and if $(\xi',\,\eta') \in \C'_{\rho^{-1}}$ then
$$
\|\eta\| \geq b\rho^{-k}\|\eta'\|.
\tag 8.6s
$$
where $k$ is the maximal number of spaced returns of $z$ in the
time interval from $N$ to $n$ and 
$$
b = \sqrt{1-\rho^4}.
$$
\endproclaim
 
\demo{Proof}
It follows from \thetag {8.2} that for any $x \in \Cal U$
$$
\C_\rho \subset \C_{\rho^{-1}} \subset \C(x).
$$
Hence
$$
D_zT^{n-N}\C_{\rho^{-1}} \subset \C(T^{n-N}z).
$$
Now \thetag {8.5u} follows from \thetag {8.3}.
 
Let us further note that 
 \thetag {8.3} implies that for any $x \in \Cal U$
$$
 D_xT^{-N}\C'_{\rho} \subset \C'(T^{-N}x).
$$
We obtain \thetag {8.5s} by applying first $D_{z'}T^{-N}$, then
$D_{T^{-N}z'}T^{-n+N}$ and using \thetag {8.2} again.
 
The properties \thetag {8.6u} and \thetag {8.6s} follow from
\thetag {8.5u} and \thetag {8.5s} respectively in exactly the same way.
We will prove only the unstable version. To measure vectors in $\C_\rho$
we use the form  $\Q$ associated with the sector $\C_{\rho^{-1}}$. It is
equal to
$$
\rho^{-1}\|\xi\|^2 - \rho\|\eta\|^2
$$
and on every spaced return to $\Cal U$ the value of this form on
vectors from $\C_{\rho^{-1}}$ gets increased by at least the factor
$\rho^{-2}$, cf. Propositions 5.11 and 6.3 . 
It remains to compare the value of
this form at $(\xi,\,\eta) \in \C_\rho$ with $\|\xi\|^2$.
We have
$$
\rho^{-1}\|\xi\|^2 \geq \rho^{-1}\|\xi\|^2 - \rho\|\eta\|^2
 \geq (\rho^{-1}-\rho^3)\|\xi\|^2
$$
which immediately yields \thetag {8.6u}.
 
\enddemo
 
Having achieved the symmetry with respect to the direction of time
 we will restrict the discussion in
the next section to the case of unstable manifolds using the unstable
version of Proposition 8.4. It can be then repeated for the stable
manifolds with the use of the stable version.
 
\proclaim{Remark 8.7}
If  $p$ is not a periodic point then by reducing the neighborhood
$\Cal U$ we can guarantee that any successive visits to $\Cal U$ are
spaced by, at least, a time $N$. In such a case the number of spaced
returns becomes simply the number of returns to $\Cal U$.
It is so also if $N=1$.
\endproclaim

 
\vskip.7cm
\subhead \S 9. UNSTABLE MANIFOLDS IN THE NEIGHBORHOOD $\Cal U$
\endsubhead
\vskip.4cm
 
Let us repeat the properties of $T$ and $\Cal U$
established in the previous section which we will rely upon.
Note that the original point $p$ does not appear explicitly.
 
There is a positive number $\rho < 1$ such that
for any $z \in \Cal U$
$$
\C_\rho \subset \C_{\rho^{-1}} \subset \C(z)
\tag 9.1
$$
and for any
$y \in T^{-N}\Cal U$
$$
D_yT^N \C(y) \subset \C_\rho .
\tag 9.2
$$
 
It follows that
if $z\in \Cal U$ and $T^nz\in \Cal U$ for  $n \geq N$ then
$$
D_zT^n \C_{\rho^{-1}} \subset \C_\rho.
\tag 9.3
$$
Moreover  if
$$
(\xi,\,\eta) \in \C_\rho \text{ and } (\xi',\,\eta') = D_zT^n
(\xi,\,\eta)
$$
then
$$
\|\xi'\| \geq b\rho^{-k}\|\xi\|
\tag 9.4
$$
where $k$ is the maximal number of spaced returns to $\Cal U$ between
the times $N$ and $n$ and $b= \sqrt{1-\rho^{4}}$.
 
By the Pesin theory \cite {P} in the smooth case and the
Katok-Strelcyn theory \cite {K-S} in the general case for almost all
$z \in \Cal U$ we have a local
unstable manifold $W_{loc}^u(z)$ through $z$.
Further the tangent spaces of $W_{loc}^u(z) \cap \Cal U$
are Lagrangian subspaces contained in $\C_\rho$.
Unfortunately the general theory does not give us a good hold on their size.
 
Let $\pi_i : \Cal V \times \Cal V \to \Cal V,\ i=1,2,$
be the projection on the first and second component respectively.
We denote by
$\Cal B(c;r)$ the open ball with the center at $c$ and the radius $r$.
 
\proclaim{Definition 9.5 }
We say that an unstable manifold in $\Cal U$ of a point
$z=(z_1,\,z_2) \in \Cal U$ has
size $\varepsilon$ if it contains the graph of a
smooth mapping from
$\Cal B(z_1;\varepsilon)$ to $\Cal V$. We denote such a graph by
$W_{\varepsilon}^u(z)$ and we will call it the unstable manifold of size
$\varepsilon$.
\endproclaim
 
By the definition of an unstable manifold $W_{\varepsilon}^u(z)$ of
size $\varepsilon$ its projection onto the first component is the
open ball with the center at $\pi_1z$ and radius $\varepsilon$.
 
\proclaim{Lemma 9.6}
The projection onto the second component of an unstable manifold through
$z =(z_1,\,z_2)\in \Cal U$ of size $\varepsilon$ lies in the
open ball with the center at  $z_2$ and the radius $\rho\varepsilon$, i.e.,
$$
\pi_2 \left(W_{\varepsilon}^u(z)\right) \subset \Cal B(z_2;
\rho\varepsilon).
$$
\endproclaim
 
\demo{Proof}
Let $W_{\varepsilon}^u(z)$ be the graph of
$$
\psi :  \Cal B(z_1;\varepsilon) \to \Cal V.
$$
The subspace $\{(\xi,\,D\psi\xi)| \xi\in \Re^{d}\}$ is tangent
to $W_{\varepsilon}^u(z)$ and hence is contained in $\C_\rho$. It follows
that
$$
\|D\psi\| \leq \rho.
$$
By the mean value theorem if $z'=(z_1',\,z_2') \in W_{\varepsilon}^u(z)$
then
$$
\|z_2' - z_2\| = \|\psi(z_1') - \psi(z_1)\| \leq
\sup \|D\psi\| \|z_1' - z_1\| < \rho \varepsilon .
$$
\enddemo
 
In contrast to the model problem at the beginning where we had fairly
long initial unstable leaves and then we cut them because of the
discontinuity of our system we start here with small unstable manifolds
and ``grow'' them until they are large or until they hit the singularity
whichever comes first. This is done in the proof of the following
Theorem.
 
\proclaim{Theorem 9.7}
For any $\delta >0 $ almost every point
$z$ in $\Cal U^1_\delta $,
$$
\Cal U^1_\delta = \Cal U_{a_1(\delta)}
$$
where $a_1(\delta) =a-b^{-1}\delta $
($\U_a$ is defined in \S 8 and $b= \sqrt{1-\rho^4}$),
either has an unstable manifold of size $\delta$ or it has an
unstable manifold of size $\delta' < \delta$ such that the
closure of $W^u_{\delta'}(z)$ intersects $\bigcup_{j> N}T^j\R^-$.
\endproclaim
 
\demo{Proof}
Let $\Cal A(\varepsilon) \subset \Cal U^1_\delta $ be the set of points
which have unstable manifolds of size $\varepsilon$. By the Katok-Strelcyn
theory almost all points in $\Cal U^1 _\delta $ belong to
$\bigcup_{\varepsilon > 0}\Cal A(\varepsilon)$. Let us fix
$\Cal A(\varepsilon)$ of positive measure and let $k$ be the smallest
natural number such that
$$
b\rho^{-k}\varepsilon \geq \delta.
$$
Almost all points in $\Cal A(\varepsilon)$ have $k$ spaced returns
to $\Cal A(\varepsilon)$ in the past. Let $z$ be such a point and let
$$
 -N \geq -i_1 > \dots > -i_k = -n
$$ 
be the $k$ times of spaced returns of this point, i.e.,
$$
T^{-i_j}z \in \Cal A(\varepsilon),\ j =1,\dots,k.
$$
 
The geometric idea for growing unstable manifolds is to take
the unstable manifold of size $\varepsilon$ through the point
$T^{-n}z$ and map it forward under $T^n$. The expansion property
\thetag {9.4} guarantees then that the image contains the unstable
manifold of size $\delta$. There are two complications in this argument.
First it may happen that $T^{n}$ is not continuous on
the unstable manifold
$W^u_\varepsilon(T^{-n}z)$, that is
$$
 W^u_{\varepsilon}(T^{-n}z)\cap \R^+_n \neq\emptyset .
$$
The other problem occurs when parts of the images of the unstable manifold
are outside of $\Cal U$ where the expansion property \thetag {9.4}
may fail.
 
To present clearly the core of the argument we ignore for the time
being these two difficulties
and assume that $T^n$ is differentiable on $W^u_\varepsilon(T^{-n}z)$
and that
$$
T^{n-i_j}W^u_\varepsilon(T^{-n}z) \subset \Cal U,\ j = 0,\dots,k,
$$
here we set $i_0=0$.
We can prove then that $z$ has an unstable manifold of size $\delta$.
Indeed let $W^u_\varepsilon(T^{-n}z)$ be the graph of
$$
\psi : \Cal B(\pi_1(T^{-n}z);\varepsilon) \to \Cal V
$$
and let us consider the map
$$
\varphi : \Cal B(\pi_1(T^{-n}z);\varepsilon) \to \Cal V
$$
defined by
$\varphi (x) = \pi_1\left( T^n(x,\psi x)\right)$.
By \thetag {9.4} this map is
an expanding map with the coefficient of expansion not less than
$b\rho^{-k}$, i.e.,
$$
\|D\varphi\xi\| \geq b\rho^{-k}\|\xi\|.
$$
Hence the image of $\Cal B(\pi_1(T^{-n}z);\varepsilon)$ by $\varphi$
contains the ball $\Cal B(\pi_1z;\delta)$. Additional complication
is caused by the fact that $\varphi$ is not necessarily one-to-one.
But since $\varphi$ is a local diffeomorphism we can define $\varphi^{-1}$
on $\Cal B(\pi_1z;\delta)$ as the branch of the inverse for which
$\varphi^{-1}\pi_1z = \pi_1(T^{-n}z)$. Therefore,
$T^nW^u_{\varepsilon}(T^{-n}z)$ contains the graph
of the map
$$\pi_2\circ T^n\circ(id\times\psi)\circ\varphi^{-1}
$$
which defines $W^u_{\delta}(z)$.
 
Let us now address the general case. We will construct the maximal 
subset of $W^u_\varepsilon(T^{-n}z)$ on which $T^n$ is
differentiable and its images at the return times to $\Cal U$ are
contained in $\Cal U$. Our first step is to consider
the connected component of
$$
W^u_\varepsilon(T^{-n}z) \setminus \R^+_n
$$
which contains $T^{-n}z$ and denote it by
$\widetilde{\widetilde {W^u_\varepsilon}}(T^{-n}z)$. Further 
 the connected component of
$$
\bigcap_{j=0}^k T^{i_j-n}
\left(T^{n-i_j}\widetilde{\widetilde {W^u_\varepsilon}}(T^{-n}z)
\cap\Cal U\right)
$$
which contains $T^{-n}z$ will be denoted it by
$\widetilde{W^u_\varepsilon}(T^{-n}z)$.
It is the part of the unstable manifold
which has the desired properties. 
 	
Now we consider the image
$$
T^{n}\widetilde {W^u_\varepsilon}(T^{-n}z)
$$
and we let  $\delta'$ be 
the largest positive number such that $W^u_{\delta'}(z)$ 
is well defined and contained
in  $T^{n}\widetilde {W^u_\varepsilon}(T^{-n}z)$. 

If $\delta' \geq \delta$ then we are done. Let us hence assume that
$\delta' < \delta $.
 
It follows from the maximality of $\delta'$ that the boundary of
$W^u_{\delta'}(z)$ contains, at least, a point from the boundary of
$T^{n}\widetilde{ W^u_\varepsilon}(T^{-n}z)$. Let $z'$ be such a point.
If $z'$  belongs to
$ \bigcup_{i \geq N}^{n-1}T^{i}\R^- $
then we are again done. If not then $T^{-n}$ is differentiable at $z'$
and hence $T^{-n}z'$ belongs to the boundary of
$\widetilde {W^u_\varepsilon}(T^{-n}z)$ and it 
does not belong to $\R^+_n $.
It follows now from the construction of 
$\widetilde {W^u_\varepsilon}(T^{-n}z)$ that 
$T^{-n}z'$ must belong to the boundary of 
$W^u_\varepsilon(T^{-n}z)$ or for some $j, 0\leq j \leq k,$
$T^{-i_j}z'$ belongs to the boundary of $\Cal U$.

We will obtain now a contradiction by using the expansion property
\thetag {9.4} .
Let $W_{\delta'}^u(z)$ be the graph of
$$
\chi :  \Cal B(\pi_1z;\delta') \to \Cal V
$$
and let
$$
\gamma_0 : [0,1) \to \Cal B(\pi_1z;\delta')
$$
be the segment connecting $\pi_1z$ and $\pi_1(z')$. We consider the
preimages of the curve $\{(\gamma_0(t),\,\chi\gamma_0(t))\;|\;0\leq t< 1\}$
and obtain $\gamma_j : [0,1) \to \Cal V,\ j=0,\dots,k$ by the formula
$$
\gamma_j(t) = \pi_1\left( T^{-i_j}(\gamma_0(t),\,\chi\gamma_0(t))\right).
$$
It follows from \thetag {9.4} that the length of $\gamma_0$ is not smaller
than the length of $\gamma_j$ times $b\rho^{-j}$.
If $T^{-n}z'$ belongs to the boundary of $W^u_\varepsilon(T^{-n}z)$ then
the length of $\gamma_k$ is at least $\varepsilon$ and we get the
contradiction
$$
\delta' \geq b\rho^{-k}\varepsilon \geq \delta .
$$
Finally if $T^{-i_j}z'$ belongs to the boundary of
$\Cal U$ for some $j, 0\leq j \leq k,$ then
 $\gamma_j$ which connects $\pi_1(T^{-i_j}z) \in \Cal U^1_\delta$ and
$\pi_1(T^{-i_j}z')$ must have the length at least
$b^{-1}\delta$. We get again the contradiction
$$
\delta' \geq b\rho^{-j} b^{-1} \delta \geq \delta.
$$
\enddemo
 
\proclaim{Definition 9.8}
We say that the unstable manifold of size $\delta$
$\Cal W^u_{\delta}(z)$ is cut by $T^i\R^-, i\geq 0,$ if its boundary 
contains a point from $T^i\R^-$.
\endproclaim
 
By Theorem 9.7 to guarantee that at least some points (and in
the case of a smooth map almost all points) have unstable manifolds of
size $\delta$ we need to step away from the boundary of $\Cal U$ by at
least $b^{-1}\delta$. In the following we fix a sufficiently small
$\delta_0$ and restrict our discussions to $\Cal U^1 = \Cal
U^1_{\delta_0}$.  We can then claim that in $\Cal U^1$ almost
every point has a uniformly large unstable manifold (of size
$\delta_0$) or a smaller unstable manifold cut by some image of the
singularity set $\R^-$.

By $\bar{\Cal B}(c;r)$ we denote the closed ball with the center at
$c$ and the radius $r$.  We define a rectangle $R(z;\delta)$ with
the center at $z= (z_1,\,z_2)$ and the size $\delta$ as the Cartesian
product of closed balls
$$
R(z;\delta) = \bar{\Cal B}(z_1;\frac \delta 2)
\times\bar{\Cal B}(z_2;\frac \delta 2).
$$
 
\proclaim{Definition 9.10}
We say that the unstable manifold $W^u_{\delta'}(z')$ of $z'=(z_1',\,z_2')$
of size $\delta'$ is connecting in
the rectangle $R(z;\delta)$ with the center at $z= (z_1,\,z_2)$
and size $\delta$ if
$$
\bar{\Cal B}(z_1;\frac \delta 2) \subset
\Cal B(z_1';\delta')
$$
and
$$
\pi_2\left(W^u_{\delta'}(z')\cap R(z;\delta)\right)
\subset \Cal B(z_2;\frac \delta 2).
$$
\endproclaim 
\par
We can say equivalently that an unstable manifold $W^u_{\delta'}(z')$
is connecting in the rectangle $R(z;\delta)$ if the intersection
of $W^u_{\delta'}(z')$ with the rectangle is the graph of a smooth
mapping from the closed ball $\bar{\Cal B}(\pi_1z;\frac \delta 2)$
to the open ball $\Cal B(\pi_2z;\frac \delta 2)$.  Clearly it is
necessary that $\delta' > \delta /2$.
 
\proclaim{Definition 9.11}
For a given rectangle $R(z;\delta)$ with the center at $z=
(z_1,\,z_2)$ and size $\delta$ we define its unstable core as the
subset of those points $z'=(z_1',\,z_2') \in R(z;\delta)$ for which
$$
\rho\|z_1' - z_1\| +\|z_2' - z_2\| <
(1-\rho)\frac \delta 2.
$$
 
\endproclaim

The role of an unstable core is revealed in the following Lemma.

\topinsert
\vskip 3in
\hsize=4.5in
\raggedright
\noindent{\bf Figure 8.} The core of a rectangle.
\endinsert

\proclaim{Lemma 9.12}
If an unstable manifold $W^u_{\delta'}(z')$ of size $\delta' >
\|\pi_1z' -\pi_2z\| +\frac \delta 2$
intersects the unstable core of a rectangle $R(z;\delta)$ then it is
connecting in the rectangle.
\endproclaim
\demo{Proof}
 
Let $z=(z_1,\,z_2)$ and $z'=(z_1',\,z_2')$, let $W^u_{\delta'}(z')$
be the graph of $\psi : \Cal B(z_1';\delta')
\to \Cal V$ and let $(x_1,\,\psi x_1)$ be a point in the unstable core of
the rectangle. By the condition on $\delta'$
$$
\bar {\Cal B}(z_1;\frac \delta 2) \subset \Cal B(z_1;\delta').
$$

We have to check only that if $x \in \bar{\Cal
B}(z_1;\frac \delta 2)$ then
 
$$
\|\psi x -z_2\| < \frac \delta 2 .
$$
We have
$$
\aligned
\|\psi x -z_2\| &\leq
\|\psi x -\psi x_1\| + \|\psi x_1 -z_2\|\\
&\leq \sup\|D\psi\|\|x- x_1\| + \|\psi x_1 -z_2\|\\
&\leq \rho\|x-z_1\| + \rho\|x_1-z_1\|
+\|\psi x_1 -z_2\|\\
& < \rho \frac \delta 2 + (1-\rho)\frac \delta 2 = \frac \delta 2.
\endaligned
$$
\enddemo

The point of the above lemma is that a large unstable manifold
may fail to be connecting in a rectangle if it intersects the rectangle
too close to the boundary.

\vskip.7cm
\subhead \S 10. LOCAL ERGODICITY IN THE SMOOTH CASE
\endsubhead
\vskip.4cm
 
Contrary to the title of this section we will consider here several
propositions valid in the general case. Incidentally they will suffice
to obtain local ergodicity in the smooth case.
 
It is important to remember that all of Section 9 can be repeated for
stable manifolds. In this section we will be using both stable and unstable
manifolds.
 
\proclaim{Lemma 10.1}
If an unstable manifold and a stable manifold are connecting in a
rectangle then there is a unique point of intersection of these manifolds
in the rectangle and it belongs to the interior of the rectangle.
\endproclaim
\demo{Proof}
Let the rectangle have the center at $z=(z_1,\,z_2)$ and size $\delta$.
The intersections of the unstable
and stable manifolds with the rectangle $R(z;\delta)$ are the graphs
of the smooth mappings
$$
\psi^u : \bar{\Cal B}(z_1;\frac \delta 2)  \to
\Cal B(z_2;\frac \delta 2)
$$
and
$$
\psi^s : \bar{\Cal B}(z_2;\frac \delta 2)  \to
\Cal B(z_1;\frac \delta 2)
$$
respectively.
 
Since both $\psi^u$ and $\psi^s$ are contractions so is their
composition
$$
\psi^s\psi^u : \bar{\Cal B}(z_1;\frac \delta 2) \to
\Cal B(z_1;\frac \delta 2).
$$
Hence it has a unique fixed point $x \in \Cal B(z_1;\frac \delta 2)$.
The point
$$
\left(x,\,\psi^u x\right) =
\left(\psi^s\psi^u x,\,\psi^u x \right)
$$
is the desired intersection point.
\enddemo
 
For a rectangle $R$ we denote by $W^{(u)s}(R)$ the union of the
intersections with $R$ of all (un)stable manifolds connecting
in $R$, i.e.,
$$
W^{(u)s}(R) = \bigcup\{R\cap W^{(u)s}_{\delta'}(z')\;|\;
\ W^{(u)s}_{\delta'}(z') \ \text { is connecting in } \ R \}.
$$
The union of the unstable core and the stable core of a rectangle will be
in the following called simply the core of the rectangle.
 
\proclaim {Proposition 10.2}
For any rectangle $R \subset \Cal U^1$ if the sets $W^s(R)$
and $W^u(R)$ have positive measure then $W^s(R) \cup
W^u(R)$ belongs to one ergodic component of $T$.
 
\endproclaim
\demo{Proof}
The proof is done by the Hopf method as described in Sections 1 and 2.

Let us fix a continuous function defined on our phase space. For all
points in one (un)stable manifold the (backward) forward time
averages are the same. As shown in Section 1 the forward and backward
time averages have to coincide almost everywhere. Our goal is to
show that they are constant almost everywhere in $W^s(R) \cup W^u(R)$. 

There is a technical difficulty stemming from the fact that the
foliations into stable and unstable manifolds are not smooth in
general.  One has to use the absolute continuity of the foliations
which was proven in \cite {KS} under the conditions which fit our
scheme. (It is by far the hardest fact to prove in their theory.)

It follows from absolute continuity of the foliation into unstable
manifolds that except for the union of unstable manifolds from
$W^u(R)$ of total measure zero almost every point (with respect to the
Remannian volume in the manifold) in an unstable manifold from
$W^u(R)$ has equal forward and backward time averages.  Let us take
such a typical unstable manifold. Again by the property of absolute
continuity the union of stable manifolds in $W^s(R)$ which intersect
the distinguished unstable manifold at points where the forward and
backward time averages exist and are equal differs from $W^s(R)$ by a
set of zero measure. Hence the time average of our function is
constant almost everwhere in $W^s(R)$.  Similarly the time average of
our function is constant almost everywhere in $W^u(R)$.

Finally using the property of absolute continuity for the third time
we can claim that $W^u(R)$ and $W^s(R)$ intersect on a subset of
positive measure. Hence the time average of our function is constant
almost everywhere in $W^s(R) \cup W^u(R)$. 

To prove that $W^s(R) \cup W^u(R)$ belongs to one ergodic component 
we proceed in the same way as at the end of Section 2.
\enddemo 

We are ready to prove the local ergodicity in the smooth case
 
\demo{Proof of Main Theorem (smooth case)}
 
All the constructions started in Section 9 apply to our point $p$.
We will prove that a neighborhood $\Cal U^2$ only slightly smaller
than $\Cal U^1$ belongs to one ergodic component. Indeed according to
Lemma 9.12 all the points in the (un)stable core of a rectangle 
$R \subset \Cal U^1$ which have an (un)stable manifold of sufficiently
large size belong to $W^{(u)s}(R)$. By Theorem 9.7 in the smooth case
almost every point in $\Cal U^1$ has both the unstable manifold and the
stable manifold of size $\delta_0$. Hence by Lemma 9.12 for any rectangle
$R \subset \Cal U^1$ of size $\delta < \delta_0$ the set
$W^s(R)$ contains at least the stable core of $R$ and
$W^u(R)$ contains at least the unstable core of $R$.
Clearly then the sets $W^s(R)$ and $W^u(R)$ have positive measure
and  we can apply Proposition 10.2.
 
To end the proof we consider a family of rectangles of size $\delta
\leq \delta_0$ contained in $\Cal U^1$ whose cores cover a slightly
shrunk neighborhood $\Cal U^2 \subset \Cal U^1$. By Proposition 10.2
we can claim that each core belongs to one ergodic component. Since
the cores form an open cover of the connected set $\Cal U^2$ we can
conclude that $\Cal U^2$ belongs to one ergodic component.
\enddemo

Actually we can claim that under the assumptions of the Main Theorem
the whole neighborhood $\Cal U$ constructed in Section 8 belongs to
one ergodic component. Indeed by taking $\delta \to 0$ the above
argument applies to $\Cal U^2 \to \Cal U^1$ so that actually $\Cal
U^1$ belongs to one ergodic component. Again the $\delta_0$ in the
definition of $\U^1$ can be chosen arbitrarily small so that
also the whole neighborhood
$\Cal U$ belongs to one ergodic component. This does not strengthen
the theorem but it demonstrates the usefulness of coverings with
rectangles of size $\delta \to 0$. It will be crucial in the treatment
of the discontinuous case.
 
Let us outline the plan for proving local ergodicity in the general
case.  We cover the neighborhood $\Cal U^2$ with rectangles of size
$\delta$. At least for some rectangles $R$ the sets $W^s(R)$ and
$W^u(R)$ will have positive measure.  We will be actually interested in the
property that these sets cover certain fixed (but otherwise
arbitrarily small) percentage of the core of the rectangle and we will
call such rectangles connecting. One may then expect to have more
connecting rectangles as $\delta \to 0$. The precise formulation of
such a property is the subject of Sinai Theorem. The method of the
proof requires that the size of the sector satisfies $\rho<\frac 1 3$. In
applying Sinai Theorem it is convenient to work with more structured
coverings, namely the centers of the rectangles will belong to a lattice
with vertices so close that the cores of nearest neighbors rectangles will
overlap almost
completely.  Consequently, if both nearest neighbors $R_1$ and $R_2$ are
connecting then the union of $W^s(R_1) \cup W^u(R_1)$ and $W^s(R_2)
\cup W^u(R_2)$ belongs to one ergodic component (see Preposition 2.3).
It will follows from Sinai
Theorem that the network of connecting rectangles becomes more and
more dense as $\delta \to 0$ so that we will be able to
claim that one ergodic
component reaches from any place in the neighborhood $\Cal U^1$ to any
other place.  We will conclude by using the Lebesgue
Density Theorem to show that $\Cal U^2$ belongs to one ergodic component.

\vskip.7cm
\subhead \S 11. LOCAL ERGODICITY IN THE DISCONTINUOUS CASE
\endsubhead
\vskip.4cm
 
Given $\delta > 0$ we consider a shrunk neighborhood $\Cal
U^2_\delta$ defined by the requirement that a rectangle with the
center in $\Cal U^2_\delta$ and size $\delta$ lies completely in $\Cal
U^1$. (One can easily see that $\Cal U^2_\delta = \Cal U_{a_2(\delta)}$
where $a_2(\delta) = a_1(\delta_0) - \frac \delta 2$). Let us note that
$\Cal U^2_\delta \to \Cal U^1$ as $\delta \to 0$.
 
Let $\Cal N(\delta,c)$ be the net defined by
$$
\Cal N(\delta,c)=
\{c\delta (m,\,k)\in \Cal U^2_\delta\;|\; m,\,k\in\Z^d\}.
$$
We consider the family $\Cal G_\delta$ of all rectangles with the
centers in $\Cal N(\delta,c)$ and size $\delta$
$$
\Cal G_\delta = \{R(z;\delta)\;|\; z \in \Cal N(\delta,c)\}.
$$
If $c$ is sufficiently small the family $\Cal G_\delta$ is a covering
of $\Cal U^2_\delta$.  The parameter $c$
will be chosen later to be very small so that many rectangles in 
$\Cal G_\delta$ overlap. But once $c$ is fixed a point may belong to at most
a fixed number of rectangles, which we denote by $k(c)$ (it does not
depend on $\delta$).
 
\proclaim{Definition 11.1}
Given $\alpha,\, 0 < \alpha < 1,$ we call a rectangle $R \in \Cal G_\delta$ 
$\alpha$-connecting in the (un)stable direction (or simply
connecting) if at least the $\alpha$ part of the measure of the
(un)stable core of $R$ is covered by $W^{(u)s}(R)$.
\endproclaim
 
\proclaim{Sinai Theorem 11.2}
If $\rho < \frac 13 $ then there is $\alpha, 0 < \alpha < 1,$ such that
for any $c$
$$
\lim_{\delta\to 0}\delta^{-1}
\mu\left (\bigcup\{R\in\Cal G_\delta\;|\;R
\text{ is not } \alpha\text{-connecting }\}\right )=0,
$$
i.e., the union of rectangles which are not $\alpha$-connecting in either
the stable or the unstable direction has measure $o(\delta)$
\endproclaim
 
It is very important for the application of this theorem that given
$\rho < \frac 13 $ we get a certain $\alpha$  (which may be very small if
$\rho$ is close to $\frac 13$) and we are free to choose $c$ (which
determines the overlap of the rectangles in $\Cal G_\delta$) as small
as we may need.
 
We will prove Sinai Theorem in Sections 12 and 13. In the remainder of
this Section we will show how to obtain the Main Theorem in the
discontinuous case  from Sinai Theorem.
 
We start with some auxiliary abstract facts. The first one concerns
Measure Theory. For any finite subset $S$ we will denote by $|S|$ the
number of elements in $S$.
 
\proclaim{Lemma 11.3}
Let $\{A_s\;|\; s\in S\}$ be a finite family of measurable subsets
of equal measure $a$ in the measure space $\left(X,\,\nu\right)$
such that no point in $X$ belongs to more than $k$ elements of the
family. For any subfamily $\{A_s\;|\; s\in S_1\}, S_1 \subset S,$ we have
$$
\frac a k |S_1| \leq \nu\left(\bigcup_{s\in S_1}A_s \right) \leq a |S_1|.
$$
Further if for a measurable subset $Y \subset X$ and some
$\alpha, 0 <\alpha < 1,$
$$
\nu(A_s\cap Y) \geq \alpha \nu(A_s) \ \text { for } \ s \in S_1
$$
then
$$
\nu (Y)\geq\nu\left(\bigcup_{s\in S_1}A_s \cap Y \right) \geq
\frac \alpha k \nu\left(\bigcup_{s\in S_1}A_s \right).
$$
\endproclaim
\qed

The second fact concerns Combinatorics. Let us consider the lattice
$\Z^d$ and its finite pieces
$$
L_n =L_n(d) = \{0,1,\dots ,n-1\}^d \subset \Z^d.
$$
Let $K \subset L_n$ be an arbitrary subset
which we call a configuration. We think of elements of $K$ as occupied
sites and elements of $L_n \setminus K$ as empty sites.
 
For a given configuration $K \subset L_n$ we consider the graph
obtained by connecting by straight segments all pairs of occupied
sites which are nearest neighbors. Let $gK \subset K$ be the family
of sites in the largest connected component of the graph.
 
\proclaim{Proposition 11.4}
Let $K_n \subset L_n(d), n= 1,2,\dots ,$ be a sequence of configurations.
If
$$
n\frac {|L_n\setminus K_n|} {|L_n|} \rightarrow 0 \ \ \text {as} \ \
n \rightarrow +\infty
$$
then
$$
\frac {|gK_n|} {|L_n|} \rightarrow 1 \ \ \text {as} \ \
n \rightarrow +\infty.
$$
\endproclaim
\demo{Proof}
This proposition will follow immediately from the following
combinatorial Lemma.
\proclaim{Lemma 11.5}
Let $K \subset L_n(d)$ be an arbitrary configuration. If
$$
\frac {|L_n\setminus K|} {n^{d-1}} < a <1
$$
then
$$
\frac {|gK|} {n^d} \geq 1-(d-1)a.
$$
\endproclaim
\demo{Proof}
The proof is by induction on $d$. For $d=1$ the statement is obvious.
Suppose it is true for some $d$. We will establish it for $d+1$.
 
We partition $L_n(d+1)$ into subsets $L_n(d)\times\{ i \}, i =
0,\dots,n-1$ and we call them floors. We pick the floor with the
fewest number of empty sites.  Clearly the number of empty sites there
does not exceed $an^{d-1}$ so that we can apply to it the inductive
assumption. We obtain in this floor a connected graph with at least
$(1-(d-1)a)n^d$ elements.
 
Now we partition $L_n(d+1)$ into subsets $\{ z \} \times
\{0,\dots,n-1\}, z \in L_n(d)$ and we call them columns. A column is
called an elevator if all of its elements are occupied. The number of
elevators is at least $(1-a)n^d$. Hence the number of elevators which
intersect the connected graph in the floor considered above is at
least $(1-da)n^d$. Adding these elevators to the graph we obtain a
connected graph with at least $(1-da)n^{d+1}$ elements which ends the
proof of the inductive step.
\enddemo
\enddemo
 
\demo{Proof of Main Theorem (Discontinuous case)}
All the constructions of Sections
8 through 10 apply with some $\rho < \frac 13$.  We will be proving
that the neighborhood $\Cal U^1$ belongs to one ergodic component.
 
The Sinai Theorem gives us $\alpha < 1$ which depends only on $\rho$
and may have to be very small if $\rho$ is very close to $\frac 13$.
Let us consider the lattice $\Cal N(\delta,c)$ and the covering 
$\Cal G_\delta$. We choose $c$ so small that if the centers of two
rectangles in $\Cal G_\delta$ are nearest neighbors in $\Cal
N(\delta,c)$ then their unstable cores (and then automatically also
stable cores) overlap on more than $1-\alpha$ part of their measure.
Note that such a property depends on $c$ but is independent of the
value of $\delta$.  This choice of $c$ has the following consequence.
If two rectangles $R_1$ and $R_2$ with centers at nearest
neighbors in $\Cal N(\delta,c)$ are $\alpha$-connecting in the
unstable direction then $W^u(R_1)$ and $W^u(R_2)$ intersect
on a subset of positive measure. If in addition we also know that
$W^s(R_1)$ and $W^s(R_2)$ have positive measure then using
Proposition 10.2 we obtain that
$$
W^u(R_1) \cup W^u(R_2) \cup W^s(R_1) \cup W^s(R_2)
$$
belongs to one ergodic component.

We consider the configuration $\Cal K(\delta)$ in the lattice $\Cal
N(\delta,c)$ which consists of the centers of all rectangles in
$\Cal G_\delta$ which are $\alpha$-connecting both in the stable and
unstable directions. As in the discussion proceeding Proposition 11.4
we consider the graph obtained by connecting with straight segments
all pairs of nearest neighbors in $\Cal K(\delta)$. Let as before
$g\Cal K(\delta)$ be the collection of vertices in the largest
connected component of this graph. By our construction the set
$$
Y(\delta) = \bigcup \{W^u(R(z;\delta))\cup W^s(R(z;\delta)) \;|\;
z \in g\Cal K(\delta)\}
$$
belongs to one ergodic component. This set is crucial in our proof
that $\Cal U^1$ belongs to one ergodic component. It may be very small
in measure (if $\alpha$ is small) but it covers at least certain fixed
$\alpha'$ portion of the measure of each of the rectangles with centers
in $g\Cal K(\delta)$, i.e.,
$$
\mu \left( R(z;\delta) \cap Y(\delta)\right) \geq \alpha'
\mu \left( R(z;\delta)\right)
\tag{11.6}
$$
for any $z \in g\Cal K(\delta)$ ($\alpha'$ is smaller than $\alpha$
since $\alpha$ is only the part of the measure of the (un)stable core
covered by the connecting (un)stable manifolds). It remains to show
that the points in $g\Cal K(\delta)$ reach into all parts of $\Cal
U^1$. It will follow from Sinai Theorem.
 
By Sinai Theorem the total measure covered by rectangles
which are not $\alpha$-connecting is $o(\delta)$.
Using Lemma 11.3 we can translate this estimate as
$$
k(c)^{-1} |\Cal N(\delta,c) \setminus \Cal K(\delta)| \delta^{2d}
= o(\delta).
$$
Since in addition
$$
\frac {|\Cal N(\delta,c)|} {(c\delta)^{2d}} = O(1)
$$
we see that the assumptions of Proposition 11.4 are satisfied and we
can claim that
$$
\frac {|g\Cal K(\delta)|}{|\Cal N(\delta,c)|} \to 1 \ \text { as } \
\delta \to 0.
\tag {11.7}
$$
 
We are ready to finish the proof by a contradiction. Suppose there
are two $T$ invariant disjoint subsets $E_1$ and $E_2$ which have
intersections with $\Cal U^1$ of positive measure. Let us pick two
Lebesgue density points $p_1$ and $p_2$ for $E_1 \cap \Cal U^1$
and $E_2 \cap \Cal U^1$ respectively. Next we fix cubes $C_1$ and $C_2$
with centers at $p_1$ and $p_2$ so small that
$$
\mu(C_i \cap E_i)   \geq \left(1 - \frac {\alpha'}{2k(c)}\right)\mu(C_i),
\ i = 1,2.
$$
It follows from \thetag{11.7} that
$$
\frac {|\left(\Cal N(\delta,c) \setminus g\Cal K(\delta)\right)\cap
C_i|}{|\Cal N(\delta,c)|} \to 0 \ \text { as } \
\delta \to 0, \ i=1,2.
$$
Since
$$
\frac {|\Cal N(\delta,c)|}{|\Cal N(\delta,c)\cap C_i|} = O(1), \ i=1,2,
$$
we conclude that
$$
\frac {|\left(\Cal N(\delta,c)\cap C_i \right)
\setminus g\Cal K(\delta)|}
{|\Cal N(\delta,c)\cap C_i|} \to 0 \ \text { as } \
\delta \to 0, \ i=1,2.
$$
Now we get immediately that
$$
\mu \left(\left(\bigcup\{R(z;\delta) | z\in g\Cal K(\delta)\cap C_i\}\right)
\triangle C_i\right) \to 0  \ \text { as } \
\delta \to 0, \ i=1,2,
\tag{11.8}
$$
where $\triangle$ denotes the symmetric difference, i.e., for any two
sets $A$ and $B$ 
$$
A\triangle B = (A\setminus B)\cup (B\setminus A).
$$
By \thetag{11.6} and Lemma 11.3
$$
\mu \left(\bigcup\{R(z;\delta) | z\in g\Cal K(\delta)\cap C_i\}
\cap Y(\delta)
\right) \geq \frac {\alpha'}{k(c)}\mu \left(\bigcup\{R(z;\delta) |
z\in g\Cal K(\delta)\cap C_i\}\right),
$$
$ i =1,2$.
 
Comparing this with \thetag{11.8} and remembering how dense $E_i$ is in
$C_i, \ i=1,2,$ we conclude that for sufficiently small $\delta$ the
set $Y(\delta)$ must intersect both $E_1$ and $E_2$ over subsets of
positive measure which contradicts the fact that it belongs to one
ergodic component.
\enddemo

          \vskip.7cm
\subhead \S 12. PROOF OF SINAI THEOREM
\endsubhead
\vskip.4cm
We will be proving only the unstable version of the theorem, i.e.,
we will estimate the measure of the union of rectangles which are
not $\alpha$-connecting in the unstable direction. Everything can be then
repeated for the stable manifolds.
 
For a point $y = (y_1,\,y_2)$ in the core of a rectangle $R(z;\delta)$
there are two possibilities:

\roster
\item the point $y$ has an unstable manifold of size $\delta' >
\|y_1 -\pi_1z\| + \frac\delta 2 $ (which is connecting in
$R(z;\delta)$ by Lemma 9.12),
 
\item the point $y$ has an unstable manifold of size $\delta' \leq
\|y_1 -\pi_1z\| + \frac\delta 2 $ cut by $\bigcup_{i\geq 0}T^i\Cal S^-$.
\endroster

If a rectangle $R(z;\delta)$ is not connecting then
the second possibility must occur for at least $1-\alpha$ part of its
core.
 
The neighborhood $\Cal U$ was chosen so small that
$\Cal S^-_N = \bigcup_{i=0}^{N-1}T^i\Cal S^-$ is disjoint from
$\Cal U$. It follows that,
for points in $\Cal U^1$, the unstable manifolds of size $\delta' <
\delta_0$ cannot be cut by these singularities.
For any $M \geq N$ let us introduce the following special case of the
second property:
\vskip.2cm
 $(2_M)$ \vtop{\hsize=4.5in\noindent\it the point $y$ has an unstable
 manifold of size $\delta' \leq
\|y_1 -\pi_1z\| + \frac\delta 2 $ cut by $\bigcup_{i = N}^MT^i\Cal S^-$.}
\vskip.2cm
Further, we introduce the auxiliary notion of a
$M$-nonconnecting rectangle. Roughly speaking, it
is a rectangle which is not connecting because of
the singularity set $\bigcup_{i=N}^M T^i\Cal S^-$.

\proclaim{Definition 12.1} Given $\alpha <\frac 12$ 
we say that a rectangle $ R $ of size $\delta$ is $M$-nonconnec\-ting, if
at least $1- 2 \alpha $ part of the measure of the unstable
core of $ R$ consists of points which satisfy the property $(2_M)$.
\endproclaim
 
The plan of the proof is the following. We fix an arbitrary positive
$\varepsilon > 0$ and we divide the argument in two parts. In one part
we will prove that there is $M = M(\varepsilon)$
and $\delta_\varepsilon$ such that, for all
$\delta < \delta_\varepsilon$, the total measure of all rectangles
in $\Cal G_\delta$ which are not $\alpha$-connecting and are not
$M$-nonconnecting
is less than $\delta\frac \varepsilon 2$. This is the subject of the `tail
bound' (section 13) and it is by far the hardest part of the proof.
It will require
global considerations (i.e., outside of $\Cal U$). The  particular value
of $\alpha$ is immaterial there.
 
We will start with the easier part proving that, for a given
$\rho < \frac 13$ and any $M$,
there are $\alpha$ and $\delta_\varepsilon$ such that,
for all $\delta < \delta_\varepsilon$, the
total measure of all $M$-nonconnecting rectangles of size $\delta$
is less than $\delta\frac \varepsilon 2$. Let us
formulate it in a separate Proposition. Its proof will be completely
confined to the neighborhood $\Cal U$.
\proclaim{Proposition 12.2}
For any $\rho < \frac 13$, there is $\alpha, 0 < \alpha < 1$, such that,
for any $M \geq N$,
$$
\lim_{\delta\to 0}\delta^{-1}
\mu\left (\bigcup\{R\in\Cal G_\delta \;|\;R
\text{ is } M\text{-nonconnecting }
\}\right )=0.
$$
\endproclaim
\demo{Proof}
We rely on our assumption that $\Cal S^-$ and its images are sufficiently
`nice'. More precisely we have required that the singularity set
$\Cal S^-_{M+1} = \bigcup_{i=0}^MT^i \Cal S^- $ is regular.
The definition of regularity was tailored to the needs of this proof.
In particular the singularity set $\Cal S^-_{M+1}$
is a finite union of pieces of submanifolds $I_k$
of codimension one, with boundaries $\partial I_k,
k=1,\dots,p$. The boundaries $\partial I_k, k=1,\dots,p$ are
themselves also finite unions of compact subsets of submanifolds of 
codimension $2$ . What is more
$$
I_k\cap I_l \subset \partial I_k \cup \partial I_l\  \text{ for any }\ k, l.
$$
 
In each of the closed manifolds $I_k, k=1,\dots,p$, we consider the
open neighborhood of the boundary of radius $r$, and we denote
by $J_r$ the union of these neighborhoods, i.e.,
$$
J_r = \bigcup_{k=1}^p \{ p\in I_k \;|\; d (p,\, \partial I_k) < r \}.
$$
 
For each $\delta$ let $r(\delta)$ be the
smallest $r$ such that, for any $k \neq l$, the distance of $I_k \setminus
J_r$ and $I_l \setminus J_r$ is not less than $2\delta$.
(In other words, for any $k \neq l$, the sets  $I_k \setminus J_r$
and $I_l \setminus J_r$ are disjoint compact subsets, and their
distance is at least $2\delta$.) 
Clearly
$$
\lim_{\delta \to 0} r(\delta) = 0.
$$
Hence, by the property \thetag{7.3} 
$$
\lim_{\delta \to 0}\mu_{\Cal S}(J_{r(\delta)}) = 0
\tag {12.3}
$$
where $\mu_{\Cal S}$ is the natural volume element on $\Cal S_{M+1}^-$.
 
Let us note that, if a rectangle $R = R(z;\delta)$
contains a point with the unstable manifold of size $\delta'  <
\delta$ cut by $S^-_{M+1}$, then it intersects the $2\delta$-neighborhood
of $S^-_{M+1}$, but it does not necessarily intersect  the singularity set
itself. For technical reasons, we prefer to blow up every rectangle, so
that the blown up rectangle must intersect $S^-_{M+1}$ itself, and not
only its neighborhood.
For a fixed $b_0 < \frac 13$, to be chosen later, and for any rectangle
$R = R(z;\delta)$, we
introduce the blown up  rectangle
$$
\widetilde R=
\Cal B(\pi_1 z,\,(1+2b_0)\frac\delta 2) \times \Cal B(\pi_2 z,\,
\frac\delta 2).
$$
The diameter of $\widetilde R$ is less than $2\delta$, since we assume
that $b_0 < \frac 13$.
 
Let $y$ belong to the core of $R$, satisfy the property $(2_M)$, and
$$
\|\pi_1 y - \pi_1 z\| \leq b_0\frac\delta 2.
$$
This implies that the unstable manifold $W^u_{\delta'}(y)$
is contained in $\widetilde R$,
so that $\widetilde R$ intersects $\bigcup_{i=N}^M T^i\Cal S^-$.
We conclude that, for $\alpha$ sufficiently small, if a rectangle
$R $ of size $\delta$ is $M$-nonconnecting, then
$\widetilde R$ intersects at least one of the submanifolds $I_k,
k=1,\dots,p.$
If for a rectangle $R$ of size $\delta$ the blown up rectangle
$\widetilde R$
intersects two
submanifolds $I_k$ and $I_l, k \neq l$ then, by definition of
$r(\delta)$ it must intersect
$J_{r(\delta)}$, and so it must be contained in the neighborhood
of $J_{r(\delta)}$ of radius $2\delta$.
By \thetag {12.3} and Proposition 7.4
the measure of the neighborhood of $J_{r(\delta)}$ of radius $2\delta$
is $o(\delta)$ (i.e., when divided by $\delta$, it tends to zero as $\delta$
tends to zero). It remains to consider those blown up rectangles
which intersect only one of the submanifolds $ I_k, k=1,\dots,p.$

The proof  will be finished when we prove
that, for all sufficiently small $\delta$, if a blown up rectangle 
$\widetilde R$ intersects only one
of the submanifolds $ I_k, k=1,\dots,p,$ (and does not intersect
$\partial I_k$), then the rectangle
$R$ is {\bf not} $M$-nonconnecting.
 
Our first observation is that there is a constant $K$ depending only on
the manifolds $I_k, k=1,\dots,p,$ such that for any $x,x'\in I_k$
there is $v$ in the tangent space to $I_k$ at $x$ ($ v \in \Cal T_xI_k$)
for which
$$
\| x' - x - v \| \leq K \|x' - x \|^2
\tag 12.4
$$
Here we consider the tangent space $\Cal T_xI_k$ of $I_k$ at $x$ as a
subspace in
$\Re^{d} \times \Re^{d}$. This property is a formulation of the
fact that smooth submanifolds are locally close to their tangent subspaces
and follows easily from the Taylor expansion.
 
Further, in view of the proper alignment of the
singularity manifolds, the tangent subspaces $\Cal T_xI_k, x
\in I_k \cap \Cal U^1$ must have their characteristic lines in $\C_\rho$.
 
Let us now take a rectangle $R = R(z;\delta)$ such that the blown up
rectangle $\widetilde R$ intersects $I_k$. We will show that
$\pi_2(I_k \cap \widetilde R)$ is contained in a fairly narrow layer.
To show this, let
$x = (x_1,\,x_2),x' = (x'_1,\,x'_2) \in I_k \cap \widetilde R $ and let
$v = (\xi,\,\eta) \in T_xI_k$ be the vector for which \thetag {12.4}
holds.
We pick a nonzero vector $v_0 = (\xi_0,\,\eta_0) \in T_xI_k$
with the direction of the characteristic line. 
For convenience, we scale it so that $\|\xi_0\| = 1$.
We have, by the definition of a characteristic line,
$$
\omega (v,\,v_0) = \langle \xi,\,\eta_0 \rangle -
\langle \eta,\,\xi_0 \rangle = 0.
$$
It follows that
$$
\left | \langle \eta,\,\xi_0 \rangle \right | =
\left | \langle \xi,\,\eta_0 \rangle \right | \leq
\rho \|\xi_0\| \|\xi\| = \rho \|\xi\|.
$$
Replacing $v$ by $x - x'$ in the last inequality and using
\thetag{12.4}, we get
$$
\left | \langle \xi_0 ,\,x_2' - x_2 \rangle \right |
\leq \rho (\| x_1' - x_1 \| + K \|x' - x \|^2) + K \|x' - x \|^2.
$$
 
Since both $x$ and $x'$ are in $\widetilde R$, we have that
$$
\| x_1' - x_1 \| < (1+2b_0)\delta
$$
and
$$
\|x' - x \| < 2 \delta.
$$
Therefore, for any $x,x' \in I_k\cap \widetilde R$, we obtain the
inequality
$$
\left | \langle \xi_0,\,x_2' - x_2 \rangle \right |
\leq \rho (1+2b_0) \delta + const\ \delta^2
\tag {12.5}
$$
where the constant depends only on $\rho$ and $K$.
The inequality (12.5) shows that $\pi_2(I_k\cap\widetilde R)$ is contained
in a layer perpendicular to $\xi_0$ of width
$\rho (1+2b_0) \delta + const\ \delta^2$. Hence, there is $\bar x_2$
(in the `center' of the layer) such that every
$x = (x_1,\,x_2) \in I_k\cap \widetilde R$ 
must belong to the layer defined by the inequality 
$$
\left | \langle \xi_0,\,x_2 -\bar x_2 \rangle \right |
\leq \rho (1+2b_0)\frac \delta 2 + const\ \delta^2
\tag {12.6}
$$

We want to estimate the width of the layer where all the points from the
core of the rectangle with `short' unstable manifolds, cut by $I_k$,
must lie. To that end
let us take a point $y = (y_1,\,y_2)$ in the core of the rectangle
$R(z;\delta)$ and such that $\|y_1 - \pi_1z\| \leq b_0\frac \delta 2$.
If $y$ satisfies the property $(2_M)$ then by
Lemma 9.6 the projection $\pi_2W^u_{\delta'}(y)$ of the unstable
manifold lies in the ball
$$
\Cal B(y_2;\rho\delta') \subset \Cal B(y_2;\rho(1+b_0)\frac \delta 2).
$$
Assuming that $W^u_{\delta'}(y)$ is cut by $I_k$, there is 
$x = (x_1,\,x_2)\in I_k \cap \widetilde R$
for which
$$
\left | \langle \xi_0,\,y_2 - x_2 \rangle \right |
\leq  \rho(1+b_0)\frac \delta 2
$$
Hence, by \thetag{12.6}, the point $y$ must belong to the layer
defined by the inequality 
$$
\left | \langle \xi_0,\,y_2 -\bar x_2 \rangle \right |
\leq  \rho(1+b_0)\frac \delta 2+\rho(1+2b_0)\frac \delta 2 +const\
\delta^2
\tag {12.7}
$$


 The last step is to choose $b_0$ so small that this layer cannot
cover all of the core. We prefer, for convenience, to fit a Cartesian 
product into the unstable core, and to prove that a fixed part of this set
is cover by connecting manifolds.
We choose such set to be
$$
X(b_0) = \Cal B(\pi_1z;b_0\frac \delta 2) \times
\Cal B(\pi_2z; s(b_0)\frac \delta 2)
$$
where $s(b_0) = 1-\rho -\rho b_0$. By the definition of a core the set
$X(b_0)$ is contained in the core of $R(z;\delta)$, and its measure
 is not less than
certain fixed part of the measure of the core, 
depending on $b_0$ (and the dimension d)
but independent of $\delta$.
 
If the layer \thetag {12.7} is sufficiently narrow, it cannot cover
all of $X(b_0)$. The precise inequality, which guarantees that, is easily
transformed into
$$
3\rho + const\ \delta < 1 - 4\rho b_0.
\tag {12.8}
$$
After a moment of reflection the reader will realize that only if
$\rho < \frac 13$ we can choose $b_0$ so small that not only \thetag {12.8}
is satisfied, but also certain fixed part of $X(b_0)$ (depending on
$b_0$ but independent of $\delta$) is not covered by the
layer \thetag{12.7}. Thus, there is $\alpha$ sufficiently small, depending on
$\rho$ and $b_0$,
such that more than  $2\alpha$ part of the measure of the core is
free of points
satisfying the property $(2_M)$. Hence the rectangle $R$ is {\bf not}
$M$-nonconnecting.
\enddemo
 
If the reader finds it hard to follow the above argument, it is
because we strived to use as little hyperbolicity as possible on
our finite orbit. The amount of hyperbolicity is measured by the size
$\rho$ of the sector . We have managed to relax the condition on $\rho$ up
to $\rho < \frac 13$. It is not hard to see that if the last condition 
is relaxed further Proposition 12.2 will not hold in general.

\vskip.7cm
\subhead \S 13.  `TAIL BOUND'
\endsubhead
\vskip.4cm
 
We will be proving that for every $\varepsilon > 0$ there is $M$ such
that the measure of points $z \in \Cal U^1$ with the unstable
manifold of size $\delta' < \delta$ cut by
$\bigcup_{i\geq M+1}T^i\Cal S^-$ does not exceed $\varepsilon\delta$.
Comparing this set with the union of rectangles
in $\Cal G_\delta$ which are not $\alpha$-connecting and not
$M$-nonconnecting, we establish immediately that the measure of the
union can be bigger by at most an absolute (=independent of
$\delta$) factor, made up of $\rho, \alpha$ and the overlap coefficient
$k(c)$ (introduced prior to Definition 11.1). To arrive at this
conclusion it is important that we consider only the rectangles from 
the covering $\Cal G_\delta$ (and not all possible rectangles of size
$\delta$).

We start by exploring some of the consequences of the 
Sinai - Chernov Ansatz.
No reference to the neighborhood $\Cal U$ will be made at this stage.
So we have assumed that almost all points in $\Cal S^-$ (with respect
to the measure $\mu_{\Cal S}$) are strictly
unbounded in the future. It follows from Theorem 6.8 that, for almost
every point $p \in \Cal S^-$,
$$
\lim_{n\to +\infty}\inf_{0 \neq v\in \C(p)}\frac{\sqrt{\Q(D_pT^nv)}}
{\|v\|} = +\infty.
$$
For a linear monotone map, let us put
$$
\sigma_*(L) = \inf_{0 \neq v\in \C(p)}\frac{\sqrt{\Q(Lv)}}{\|v\|}.
$$
Consequently, for any (arbitrarily small) $h > 0$ and any (arbitrarily
large) $t >0$, there is $M = M(h,t)$ so large  that the subset
$$
\widetilde E_t = \{p \in \Cal S^-\;|\; \sigma_*(D_pT^M) \leq t + 1 \}
$$
has measure
$$\mu_{\Cal S}(\widetilde E_t) \leq h.$$
The map $T^M$ is, in general, not even continuous in all of $\Cal
S^-$.  The coefficient $\sigma_*(D_pT^M)$ is defined only for almost
every point $p \in \Cal S^-$. Hence, so far, the subset $\widetilde
E_t$ is defined modulo subsets of measure zero. We need a closed
subset, since we plan to use Proposition 7.4.

The map $T^M$ is discontinuous on $\Cal S^+_M,$ which was assumed to
be a regular set. Using the proper alignment of singularity sets and
monotonicity of the system, we conclude that $\Cal S^+_M$ is
transversal to $\Cal S^-$ (in the natural sense).  It follows that the
set $B_M = \left( \Cal S^+_M \cup \partial\M\right) \cap \Cal S^- $ is
a finite union of compact subsets of submanifolds of dimension $2d-2$.
Further, $\Cal S^-$ is decomposed into (possibly very large)
finite number of pieces of submanifolds of dimension $2d-1$ such that
$T^M$ is differentiable in the interior of every piece, and their
boundaries are subsets of $B_M$. It follows that the coefficient 
$\sigma_*(D_pT^M)$ is continuous in the interior of every piece.

Let us choose $\zeta$  so small that the closure of the
$\zeta$-neighborhood of $B_M$
in $\Cal S^-$
$$
B_M^\zeta = \{p\in \Cal S^-\;|\; d(p,\, B_M) < \zeta\}
$$
has small measure
$$
\mu_{\Cal S}\left(\overline{B_M^\zeta}\right) \leq h.
$$
Now the set $E_t$ defined by 
$$
E_t = \widetilde E_t \setminus B_M^\zeta  = 
\{p \in \Cal S^-\setminus B_M^\zeta \;|\; \sigma_*(D_pT^M) \leq t + 1 \}
$$
is {\bf closed}, and we have
$$
\mu_{\Cal S}\left(E_t\cup\overline{B_M^\zeta}\right) \leq 2h.
$$

Let
$$
\Cal S_t = 
\{p \in \Cal S^-\setminus B_M^\zeta \;|\; \sigma_*(D_pT^M) \geq t + 1 \}.
$$
$\Cal S_t$ is a compact set and the coefficient $\sigma_*(D_pT^M)$ 
is continuous in a neighborhood of $\Cal S_t$ in $\M$. Hence, there is $r >0$
such that 
$$
\sigma_*(D_pT^M) > t,
$$
for every point $p$ in the $r$-neighborhood of $\Cal S_t$ in $\M$, let
$$
\Cal S_t^r = \{p \in \M \;|\; d(p,\,\Cal S_t) < r\}.
$$

Now we look at our neighborhood $\Cal U$. Our goal is to estimate, for
given $\delta$, the measure of the set $Y(\delta, M)$ of points in
$\Cal U^1$ which have the unstable manifold of size
$\delta' < \delta$ cut by $\bigcup_{i\geq M+1}T^i\Cal S^-$.
We will achieve this
by splitting $Y(\delta, M)$ into convenient pieces and showing that
their preimages must end up in extremely small neighborhoods of $\Cal S^-$.
 
For $z \in Y(\delta, M)$ the unstable manifold $\Cal W^u_{\delta'}(z)$
may be cut by several  (possibly infinitely many) of the singularity sets
$T^i\Cal S^-$, $i = M+1,\dots .$
Let $m(z)$ be the smallest $i \geq M+1$ such that
$\Cal W^u_{\delta'}(z)$ is cut by  $T^i\Cal S^-$.
Let further
$$
k(z) = \#\{i\;|\;1 \leq i \leq m(z) - M,\ T^{-i}z \in \Cal U^1\}.
$$
Roughly speaking $k(z)$ is the number of times the point $z$
visits in $\Cal U^1$
in the past in the time frame bounded by $m(z)$.
We put for $k =0,1,\dots, \  m = M+1, \dots ,$
$$
Y^k_m = \{z \in Y(\delta, M)\;|\; m(z) = m , k(z) = k \}.
$$
We will now fix $k$ and estimate the measure of
$$
\bigcup_{m \geq M+1} Y^k_m.
$$
\proclaim{Lemma 13.1}
For $m \neq m'$
$$
T^{-m}Y^k_m \cap T^{-m'}Y^k_{m'} = \emptyset.
$$
\endproclaim
\demo{Proof}
Let $m < m'$. If $y \in T^{-m}Y^k_m \cap T^{-m'}Y^k_{m'}$
then for $z = T^my$ and $z' = T^{m'}y$ we have
$$
k(z') \geq k(z) + 1.
$$
It contradicts the fact that $z \in Y^k_m$ and $z' \in Y^k_{m'}$.
\enddemo
 
By Lemma 13.1 we have
$$
\mu(\bigcup_{m \geq M+1} Y^k_m) \leq \sum_{m \geq M+1}\mu(Y^k_m) =
\sum_{m \geq M+1}\mu(T^{-m}Y^k_m) = \mu(\bigcup_{m \geq M+1} T^{-m}Y^k_m).
$$

Let $z \in Y^k_m $ and $z' \in T^m \Cal S^-$ be a point in the boundary
of $\Cal W^u_{\delta'} (z)$. We connect $z$ and $z'$ by the curve $\gamma$
in $\Cal W^u_{\delta'} (z)$ which projects under $\pi_1$ onto the linear
segment from $\pi_1z$ to $\pi_1z'$. In the neighborhood $\Cal U$
we have three ways of measuring the length of $\gamma$. We can use
the quadratic form $\Q$, or the length of the projection onto the first
component, or finally, we can use the  Riemannian metric.
All these metrics are equivalent in $\Cal U$ and we will use the following
coefficients defined by their ratios
$$
\sup\left\{\frac {\| v \|}{\|\xi\|}\;|\;0 \neq v =(\xi,\,\eta)\in
\C_\rho \right\} = \sqrt{1 + \rho^2},
$$
$$
q=
\sup\left\{\frac {\sqrt{\Q(v)}}{\|\xi\|}\;|\;0 \neq v =(\xi,\,\eta)\in
\C_\rho \right\}
$$
where the last supremum is taken also over all of $\Cal U$.

Our goal is to estimate the distance of $T^{-m}z$ and $T^{-m}z'$
in the Riemannian metric, such a distance clearly does not exceed 
the length of the curve $T^{-m}\gamma$.
To that end, let $n \leq m-M$, be the time of the $k$-th visit
in the past by $z$ to $\Cal U^1$, i.e.,
$T^{-n}z \in \Cal U^1$.
By Proposition 8.4 on every spaced return to $\Cal U$ the projection
of the preimage of $\gamma$ is contracted by at least the coefficient
$\rho$. In the $k$ visits there must be at least $\frac kN -1$
spaced returns. Hence, the projection of $T^{-n}\gamma$ has the length
which, by \thetag{8.6u}, does not exceed 
$$
c_1\lambda^k\delta,
$$
where
$$
\lambda = \rho^{\frac 1N} \ \text{ and } \ c_1 = \frac 1{\rho b} = 
\frac 1 {\rho\sqrt{1-\rho^4}}.
$$
It follows that the Riemannian length of $T^{-n}\gamma$ does not exceed
$$
c_2\lambda^k\delta,
$$
where
$$
c_2 =  \frac 1 {\rho\sqrt{1-\rho^2}},
$$
and its length in the metric $\Q$ does not exceed
$$
c_3\lambda^k\delta,
$$
where
$$
c_3 =  \frac q {\rho\sqrt{1-\rho^2}},
$$

Now we apply $T^{-(m-n)}$ to $T^{-n}\gamma$ and we use the fact that
$m-n \geq M$.
There are two different cases.

{\it Case 1.}
$$
T^{-m}z' \in E_t\cup\overline{B_M^\zeta}
$$
 
We use the noncontraction property. Under the action of $T^{-(m-n)}$
the Riemannian length of $\gamma$ 
can expand at most by the factor $\frac 1a$.
We conclude that the Riemannian length of $T^{-m}\gamma$ does not exceed
$$
\frac {c_2}a \lambda^k \delta.
$$
Thus $T^{-m}z$ belongs to the neighborhood of $E_t\cup\overline{B_M^\zeta}$
in $\M$ of this radius.
By Proposition 7.4
its measure
does not exceed
$$
3h\frac{2c_2}a  \lambda^k \delta,
\tag 13.2
$$
if only $\delta$ is small enough ($\delta \leq \delta_0$ and $\delta_0$
does not depend on $k$ or $m$).
 
{\it Case 2.}
 
$$ T^{-m}z' \in \Cal S_t $$
 
We claim that, for sufficiently small $\delta$ the length of
$T^{-m}\gamma$ does not exceed
$$
\frac 1t c_3\lambda^k\delta.
$$
Indeed, it is so if $T^{-m}\gamma$ is contained in $\Cal S_t^r$ (the
$r$-neighborhood of $\Cal S_t$ in $\M$). Since $m-n \geq M$, we have 
$$
\sigma_*(D_pT^{m-n}) > t,
$$ 
for every point $p\in \Cal S_t^r$. Hence, the length in the metric $\Q$ of
$T^{-n}\gamma$  is longer
than the Riemannian length of $T^{-m}\gamma$ by at least the factor  $t$.
If $T^{-m}\gamma$ is not contained in $\Cal S_t^r$, then there must be
a segment of this curve in $\Cal S_t^r$ which has at least length $r$.
It follows that the image of this segment under $T^{m-n}$ has the
length in the metric $\Q$ 
not less than $tr$,
which is more than the total length in the metric $\Q$ 
of $T^{-n}\gamma$ for sufficiently
small $\delta$.
This contradiction
shows that, for sufficiently small $\delta$, $T^{-m}\gamma \subset 
\Cal S_t^r$. We have proven our claim.
It follows that $T^{-m}z$ belongs to
the neighborhood of $\Cal S^-$ of radius $\frac 1t c_3\lambda^k\delta$.
Using again Proposition 7.4, we can estimate the measure of
this neighborhood by
$$
2\mu_{\Cal S}(\Cal S^-)\frac {2c_3}{t} \lambda^k\delta,
\tag 13.3
$$
if only $\delta$ is sufficiently small
($\delta \leq \delta_0$ and $\delta_0$ does not depend on $k$ or $m$).
 
Combining the estimates \thetag {13.2} and \thetag {13.3} we obtain
that for any $k=0,1,\dots ,$
$$
\mu(\bigcup_{m \geq M+1} T^{-m}Y^k_m) \leq
\left(h \frac {6c_2}a  + \frac 1t 2c_3\mu_{\Cal S}(\Cal S^-)\right)
\  \lambda^k \ \delta.
$$
It follows that
$$
\mu(Y(\delta,M)) \leq \left(h \frac {6c_2}a  
+ \frac 1t 4c_3\mu_{\Cal S}(\Cal S^-)\right)\frac 1{1-\lambda} \ \delta.
$$
The last inequality tells us how we should choose a small $h$ and a
large $t$ at the
beginning of our argument to guarantee that
$$
\mu(Y(\delta,M)) \leq \varepsilon\delta.
$$
The `tail bound' is proven.

\vskip.7cm
\subhead \S 14. APPLICATIONS
\endsubhead
\vskip.4cm

{\bf A. Billiard systems in convex scattering domains. }
\bigskip
We assume that the reader is familiar with billiard systems.
If it is not the case, we recommend  \cite {W4} for a 
quick introduction into the subject. We will rely on the results
of that paper.

Let us consider a domain in the plane bounded by a locally 
convex closed curve given by the natural equation $r = r(s), 0 \leq s 
\leq l$
describing the radius of curvature $r$ as a function of the arc length
$s$. We assume that the radius of curvature satisfies the condition 
$$
\frac {d^2r}{ds^2} < 0,  \ \text { for all } \ s,\ 0\leq s \leq l.
\tag {14.1}
$$
Curves satisfying this condition were called in \cite {W4} strictly 
convex scattering.

\proclaim{Examples}
\endproclaim
1. Perturbation of a circle.

2. Cardioid.
\bigskip
Such a domain cannot be convex, and there is a singular point in the
boundary where the curve intersects itself. (If you do not like playing
billiards on a table which is not convex, you may take the convex hull
of our domain and everything below still applies.)

The following theorem is a fairly easy consequence of the Main Theorem.
\proclaim {Theorem 14.2}
The billiard system in a domain bounded by a strictly
convex scattering curve (i.e., satisfying \thetag{14.1}) is ergodic.
\endproclaim 

Let us consider the map $T$ describing the first return map to the 
boundary. $T$ is defined on the set $\Cal M$ of unit tangent vectors 
pointing inwards. We parametrize $\Cal M$ by the arc length parameter
of the foot point 
$s, 0  \leq s  \leq l,$ and the angle $\varphi, 0 \leq \varphi \leq \pi,$
which the unit vector makes with the boundary (oriented counterclockwise).
In these coordinates $\Cal M$ becomes the rectangle $[0,l] \times [0,\pi]$.
The symplectic form (the invariant area element) is given by
$\sin \varphi \ ds \wedge d\varphi.$
After we derive the formula for the derivative of $T$, 
we will be able to check immediately that $T$ preserves this area element.

The map $T$ is discontinuous at those billiard orbits which hit the
singular point of the boundary. They form a curve $\Cal S^+$ in $\Cal
M$ which is a graph of a strictly decreasing function, decreasing 
curve for short.  This curve
divides the rectangle $\Cal M$ into two curvilinear triangles, $\Cal
M^+_b$ with a side at the bottom and $\Cal M^+_t$ with a side at the
top.  

To find the images of $\Cal M^+_b$ and $\Cal M^+_t$ we  use the
reversibility of our system. Namely, let $S : \Cal M \to \Cal M$ be
defined by $S(s,\varphi) = (s,\pi-\varphi)$.  We have 
$$ 
T \circ S = S \circ T^{-1}.  
$$ 
We can  now claim that $T^{-1}$ is continuous except on $\Cal S^- =
S\Cal S^+$ which is an increasing curve (the graph of a strictly increasing
function).
$\Cal S^+$ divides  the rectangle $\Cal M$ into two curvilinear triangles 
$\Cal M^-_b = S\Cal M^+_t$ and $\Cal M^-_t = S\Cal M^+_b$. 
We have constructed our symplectic boxes. $T$ is a diffeomorphism 
on their interiors
and a homeomorphism on the closure. The derivative of $T$ does blow up 
at least at  one point of the boundary 
$\Cal S^+$ (different for $\Cal M^+_b$ and for $\Cal M^+_t$) 
corresponding to
the two billiard orbits tangent to one of the branches of the boundary at
the singular point.  In the case of the cardioid the derivative blows up
at any point of $\Cal S^+$ and also at the vertical boundaries 
because the curvature at the cusp is infinite
(see the formula for the derivative of $T$ below). It is very handy that 
we did not have to require in Section 7 that our map is a diffeomorphism 
on the closed symplectic boxes.

The derivative of $DT$ at $(s_0,\,\varphi_0)$ has the form 
$$
\left(\matrix 
\frac{\tau-d_0}{r_0\sin \varphi_1}& \frac{\tau}{\sin \varphi_1}\\ 
\frac{\tau-d_0-d_1}{r_0r_1\sin \varphi_1}&\frac{\tau-d_1}{r_1\sin \varphi_1}
\endmatrix\right) 
\tag {14.2}
$$
where $T(s_0,\,\varphi_0) = (s_1,\,\varphi_1)$,
$\tau$ is the time between consecutive hits 
(i.e., the length of the billiard orbit segment) and 
$d_i = r_i\sin \varphi_i,\ i=1,2.$
This derivative can be obtained by straightforward implicit differentiation
but we do not recommend it. There is a more geometric (and safer) way
to obtain the derivative by resorting to the description of billiard
orbit variations by Jacobi fields. In our two dimensional situation it
amounts to introducing coordinates $(J,\,J')$ in the tangent planes of 
$\Cal M$ 
$$
\aligned
J = & \sin \varphi ds,\\
J'= & -\frac 1r ds - d\varphi.
\endaligned
\tag {14.3}
$$
The evolution of $(J,\,J')$ between collisions is given by the matrix 
$$
\left(\matrix 
1 & \tau \\
0 & 1
\endmatrix\right).
\tag {14.4}
$$
At the collision $(J,\,J')$ is changed by 
$$
\left(\matrix 
-1 & 0 \\
\frac 2{d_1} & -1
\endmatrix\right).
\tag {14.5}
$$
Now the derivative \thetag{14.2} is obtained by multiplying the matrices
\thetag{14.4} and \thetag{14.5} and taking into account \thetag{14.3}.

The geometric meaning of $d_0,d_1,$ and the inequality 
$$
\tau > d_0 +d_1
\tag {14.6}
$$
is explained at length in \cite {W4}. It was proven there that \thetag
{14.6} holds for any billiard  orbit segment, if the boundary curve is
strictly convex scattering  (actually these two
properties are essentially equivalent). It follows from \thetag {14.6}
that for a strictly convex scattering curve all elements in 
\thetag {14.2} are positive. 

We choose as our family of sectors the constant sector between the 
horizontal line $\{d\varphi =0\}$  and the vertical line $\{ds =0\}$.
We see immediately that the derivative $DT$ is strictly monotone.

We are now ready to argue that the singularity sets $\Cal S^-_n =
\bigcup_{i=0}^n T^i\Cal S^- $ are regular. We claim that $\Cal S^-_n$
is a finite union of  increasing curves
which intersect each other only at the endpoints. It can be proven by
induction.  Indeed $\Cal S^-$ is an increasing curve and so it is also
properly aligned. The singularity set $\Cal S^+$ is a decreasing curve,
and as such it may intersect each of the increasing curves of $\Cal
S^-_n$ in at most one point. Hence both $\Cal M^+_b \cap \Cal S^-_n$
and $\Cal M^+_t \cap \Cal S^-_n$ are finite unions of increasing
curves with intersections only at the endpoints.  Hence in view of the
monotonicity of our system the images under $T$ are also finite
unions of increasing curves in $\Cal M^-_b $ and $\Cal M^-_t$
respectively. It is clear that we can safely add $\Cal S^-$ to these
images. We have thus checked that $\Cal S^-_{n+1} = \Cal S^- \cup
T\Cal S^-_n$ is also a finite union of increasing curves which
intersect only at the endpoints. Note that the assumptions of
Lemma 7.6 are too restrictive to allow its application in this case.

One can easily compute (and it was done explicitly in \cite {W4})
that 
$$
\sigma(DT) = \sqrt{1+\omega} +\sqrt{\omega}, \ \text { where } \ 
\omega = \frac {(\tau -d_0 -d_1)\tau}{d_0d_1}.
\tag {14.7}
$$ 
It follows from \thetag {14.7} and from the supermultiplicativity
of the coefficient of expansion $\sigma$ that the only way in which an
orbit can fail to be strictly unbounded is when the lengths of the
segments of the orbit go to zero. It was shown by Halpern \cite {Ha}
that there are no such billiard orbits, 
if $r(s)$ is a $C^1$ function bounded away from zero.  Hence,
under such an assumption, which excludes the cardioid, all orbits for
which arbitrary power of $T$ is differentiable are strictly unbounded.
To include the cardioid, or more generally the curves with the radius
of curvature $r(s)$ decreasing monotonously to zero at the endpoints
of the interval, $0 \leq s \leq l$, (at the singular point), we shall
argue that also for this class there is no accumulation of collisions
at the singular point. Indeed, if an arc of the boundary
between two consecutive hits by the billiard ball has monotone
curvature, then the angle of incidence(reflection) is smaller where 
the curvature is bigger. Hence, as an orbit gets closer to the
singularity point (the cusp for the cardioid), it is more and more 
perpendicular to the boundary, and so it cannot accumulate at the
singularity.

 This observation takes care of the Sinai - Chernov Ansatz.  We are
also guaranteed that the coefficient $\sigma (DT^n)$ can be made
arbitrarily large by increasing $n$, except possibly for points which
end up on the decreasing curve $\Cal S^+$ in the future and the
increasing curve $\Cal S^-$ in the past. These are the points in $\Cal
S^+_n \cap \Cal S^-_m$, for some $n$ and $m$, and so there are only
countably many such points. (The orbit of such a point `dies' both in
the future and in the past, and it may fail to pick up enough
hyperbolicity before then.)  We can apply the Main Theorem to all
other points, and they form a connected set.  Hence, the local
ergodicity obtained from the Main Theorem implies ergodicity.

It remains to check the noncontraction property. It was pointed out
to us by Donnay \cite{D1} that the derivative of $T$ increases
$|J'|^2$ on nonzero vectors from the sector. Indeed the interior of the
sector is defined by
$$
\frac {J'}{J} < - \frac 1d
$$
so that we have 
$$
\frac {|J'|}{|J|}  > \frac  1d.
$$
If $DT (J_0,\,J_0') = (J_1,\,J_1')$ then we have from \thetag {14.4} and
\thetag {14.5} that
$$
J_1 = -J_0 - \tau J_0'.
$$
It follows that 
$$
|J_1'| \geq \frac 1{d_1} |J_1| = \frac 1{d_1} | J_0 + \tau J_0'| \geq
\frac {\tau}{d_1}|J_0'| - \frac 1{d_1} |J_0| \geq 
\frac {\tau - d_0}{d_1}|J_0'|.
$$
In view of \thetag {14.6} $\frac {\tau - d_0}{d_1} > 1$. So indeed $|J'|^2$
gets increased.

Moreover, for all vectors in the sector we have the following estimates
$$
2(\frac 1{r^2}ds^2 + d\varphi^2) \geq |J'|^2 = |\frac 1r ds + d\varphi|^2 
\geq \frac 1{r^2}ds^2 + d\varphi^2.
$$

The metric $\frac 1{r^2}ds^2 + d\varphi^2$ is equivalent to the
standard Riemannian metric in the $(s,\,\varphi)$ coordinates 
($ds^2 + d\varphi^2$) if only $r$ is bounded away from zero. 
Thus noncontraction is established under this additional assumption,
which excludes the cardioid. 

To cover the case of the cardioid, we observe that the noncontraction
property is used only in the proof of the `tail bound'. In that proof
some subsets
of the neighborhood $\U$ are transported back to the neighborhood of
the singularity set $\Cal S^-$. We need the property that vectors from
the sector
$\C$ are not contracted too much, along the orbits from the vicinity of
the singularity set to the neighborhood $\U$, even if the orbit is
very long. We obtain readily this property from the observation that
although $|J'|^2$ is, in general, only bigger than the scaled standard
Riemannian metric, it is clearly equivalent to one locally in the
neighborhood $\U$.

The reader may be worried that the standard Riemannian metric in the
$(s,\,\varphi)$ coordinates does not generate the invariant area
element. However, the Riemannian area is not smaller than the
symplectic area. This is sufficient for the proof of Sinai Theorem.
We could also handle this complication by introducing from the  very 
beginning  coordinates in $\M$ in which the symplectic form is standard.

We can conclude that $T$ is ergodic, and so Theorem 14.2 is proven. 

It follows from the results of Katok and Strelcyn \cite {KS} that $T$ 
is a Bernoulli system.

The framework of this paper allows to cover also the class of billiard
systems in domains with more than one smooth piece in the boundary,
which are not necessarily convex scattering.  In the recent paper
\cite {D2} Donnay introduced a natural condition (focusing arc) on the
convex pieces of the boundary of a billiard table. He proves that if
two focusing arcs are connected by sufficiently long (extremely long
may be required) straight segments, then the billiard system in such a
(stadium like) domain has nonvanishing Lyapunov exponents. This work
puts the original stadium of Bunimovich \cite {B}, which had arcs of
circles in the boundary, into a large class of billiard systems with
nonuniform hyperbolic behavior, larger than the class with convex
scattering pieces introduced in \cite {W4}. 

All the properties listed in Section 7 are satisfied for the billiards
of Donnay in a straightforward fashion, with the notable exception of
the noncontraction property. The problem is that the construction of
the bundle of sectors depends heavily on the dynamics, and it is
unlikely that there is a geometrically defined Lyapunov metric (like
$|J'|^2$ for the convex scattering curves). Instead we use the
following two ideas.

We have remarked in Section 7 that if the map $T$ is differentiable up
to and including the boundary of symplectic boxes, and $DT$ is
strictly monotone, then the noncontraction property holds
automatically. In the billiards of Donnay the sectors are pushed
strictly inside at the time of crossing from one convex piece to the
other. Hence, we can use this observation on the compact part of the
phase space made up of orbits which cross over from one convex piece
to the other. We have the noncontraction property for the return map
to this set, where we measure vectors in $\C$ using the form $\Q$
defined by the bundle of sectors uniformly larger than $\C$.  The
construction of the bundle of sectors $\C$ by Donnay and his condition
on the separation of convex pieces allows to introduce immediately
these larger sectors with respect to which the derivative of the
return map is monotone.

It remains to check the noncontraction property along `grazing'
orbits which reflect many times in one convex piece. This is
essentially done in \cite {D2}, where Lazutkin coordinates are used 
to put the map $T$ in the vicinity of the boundary into a normal form.

These two observations, put together, give us the unconditional
noncontraction property, and thus our Main Theorem applies.

\bigskip
{\bf B. Piecewise linear standard map. }
\bigskip
Let $T:\Bbb T^2 \to \Bbb T^2$ be defined by
$$
T(x_1,\,x_2) = (x_1+x_2 + Af(x_1),x_2+Af(x_1))
$$
where  $(x_1,\,x_2)$ are taken modulo 1, $f$ is a periodic function
$$
f(t) = |t| -\frac 12,  \ \ \text { for }  -\frac12 \leq t \leq \frac 12,
$$
and $A$ is a real parameter. The mapping $T$ preserves the Lebesgue
measure.
For $A=1$ there is a simple invariant domain $\Cal D$ in the torus
shown in Figure 9.
It was proven in \cite {W5} that the Lyapunov exponents are different from 
zero almost everywhere in $\Cal D$.

\topinsert
\vskip 3in
\hsize=4.5in
\raggedright
\noindent{\bf Figure 9} The domain $\Cal D$
\endinsert

\proclaim{Theorem 14.8}
$T$  is ergodic in $\Cal D$.
\endproclaim

As in the previous application it follows that $T$ is a Bernoulli
system in $\Cal D$.

All the conditions of Section 7 are satisfied here in a very simple 
fashion. The reader can find all the necessary details in \cite {W5}
and \cite {W6}. In this piecewise linear case one does not have to rely 
on the general results of Katok and Strelcyn. The existence of stable
and unstable leaves can be obtained by the straightforward approach
of Sections 1-3. 

There are many other values of $A$ for which nonvanishing of 
Lyapunov exponents was established for $T$ in  some domains in the torus,
\cite {W5},\cite{W6}. The most interesting is the sequence of $A$'s
(roughly speaking) going to zero for which there is an invariant
domain, with similar geometry as $\Cal D$, where
$T$ has nonvanishing Lyapunov exponents.
It is a piecewise linear model for the unstable layer containing the
separatrices of the saddle fixed point $(0,\,\frac 14)$. One can apply Main
Theorem to all these
special domains , so that in each case the map $T$ is ergodic and hence
Bernoulli.  The reader should not have any difficulties in recovering
the details based on the two papers cited above (incidentally even the
noncontraction property was considered there).

\bigskip
{\bf C. The system of falling balls. }
\bigskip

One of the original motivations for our work was to prove ergodicity 
of the system of falling balls. 
This is a monotone system (\cite {W7}, \cite {W8}, \cite {W3}), and all
(semi-infinite)
smooth orbits are strictly unbounded. (The unboundedness of 
all orbits is obtained, under mild assumptions, by the application of 
Proposition 6.9) It follows that all Lyapunov exponents
are different from zero, and it looks like a prime candidate for the
application of Main Theorem. 
It turns out, however, that in this example  the 
singularity sets are not properly aligned, if the number of balls
is greater than two. We will show this, and briefly discuss the case of
two balls.

The system of falling balls is the system of point particles moving on
a vertical line, which also interact by elastic collisions, and are
subjected to a potential external field which forces the particles to
fall down.  To prevent the particles from falling into an abyss we
introduce the hard floor, and assume that the bottom particle bounces
back upon collision with it.  The masses of the particles are in
general different (the system of equal masses is completely integrable,
since the elastic collision of equal masses in one dimension amounts
to the exchanging of momenta).

The Hamiltonian of the system is 
$$
H = \sum_{i=1}^N\left(\frac {p_i^2}{2m_i} + m_iU\left(q_i\right)\right)
$$ 
where $q_i$ are the positions and
$p_i=m_iv_i$ the momenta of the particles, $q_i,p_i\in \Bbb{R},
i=1,\dots ,N$, and $U\left(q\right)$ is the potential of the
external field .
The differential equations of the system  are
$$ \aligned
\dot{q_i}&= \frac {p_i}{m_i}\\
\dot {p_i}&= -m_iU'\left(q_i\right),
\endaligned  $$
$i=1,\dots ,N$.

The description of the dynamics is completed by the assumptions that
 the particles are impenetrable, and that they 
collide elastically with each other and with the floor ${q =0}$.

We choose the following Lagrangian subspaces $$V_1 = \{ dp_1=\dots
=dp_N =0\} \ \ \  \text{and} \ \ \ V_2 = \{ dh_1=\dots =dh_N =0\},$$ where
$h_i=\frac {p_i^2}{2m_i} + m_iU\left(q_i\right), i=1,\dots ,N$, are 
individual energies of the particles.

We have $$ dh_i =\frac {p_idp_i}{m_i} + m_iU'\left(q_i\right)dq_i,$$ 
$ i=1,\dots ,N,$ so that $V_1$ and $V_2$ 
are indeed transversal if only $U' \neq 0$,
i.e., if the external field is actually present.

The form $\Cal{Q}$ is equal to 

$$\Cal{Q} = 
\sum_{i=1}^{N} \left(d q_i d p_i +
\frac{p_i}{m_i^2U'\left(q_i\right)}\left(d p_i\right)^2\right).$$

It was shown in the papers cited above that the system is
strictly monotone, provided that
$$
 U'\left(q\right) > 0  \ \ \text { and } \ \ U''\left(q\right) < 0 ,$$
and  
$$
m_1 > m_2 >\dots > m_N.
$$

The symplectic map $T$ that naturally arises in this system is the map 
``from collision to collision''. Our dynamical system is 
a suspension of the map.
So that the system is ergodic if and only if the map $T$ is ergodic.
As usual, the actual computations are easier done in the full phase space
of the flow.

Singularity set $\Cal S^-$ corresponds to triple collisions:
simultaneous collisions of three particles and the collision of two
particles with the floor.  Part of the first singularity set are not
properly aligned. The second set is. So the methods of this paper
apply only to the system of two particles. 

Let us show that indeed the triple collision of three particles
produces the singularity set which is not properly aligned.
We consider the manifold
$$
\{(q,\,p) | q_1 = q_2 = q_3 \}.
$$
Its tangent subspace is described by the equations 
$$
dq_1=dq_2=dq_3
$$
Its skew orthogonal complement is the two dimensional subspace given
by equations 
$$
\aligned
& d q = 0,\\
&d p_1 + d p_2 +d p_3 =0,\\
& d p_i = 0 \ \ \text { for } \ \ i \ge 4.
\endaligned
\tag {14.9}
$$
Restricting the form $\Q$ to this plane we get 
$$
\sum_{i=1}^3 \frac {p_i}{m_i^2 U'}(d p_i)^2.
\tag {14.10}
$$
We should assume that the particles emerge from collisions which
means that 
$$
\frac {p_1} {m_1} < \frac {p_2} {m_2} < \frac {p_3} {m_3}.
$$
But the momenta may, as well, be all negative which makes the quadratic form 
\thetag {14.10} negative definite. The actual characteristic line is 
obtained by intersecting the plane \thetag {14.9} by the tangent to 
the constant energy manifold. If all the momenta are negative,
it is guaranteed 
to be outside of the sector. It is not hard to compute that the precise 
condition for the characteristic line to be contained in the sector is 
$$
\frac {v_1}{m_1} \left( v_2 -v_3 \right)^2 +
\frac {v_2}{m_2} \left( v_3 -v_1 \right)^2 +
\frac {v_3}{m_3} \left( v_1 -v_2 \right)^2 \geq 0
$$
where $v_i = \frac {p_i}{m_i}, i \geq 1$ are the velocities.

We close with the discussion of the system of two balls. For clarity,
we restrict ourselves to the case of constant acceleration, 
$U(q) = q$. It was established in \cite {W7}, that also in this case 
all orbits are strictly monotone, if there are only two or three balls
and their masses decrease. (For more than three balls technical
problems arise, and it is an open problem to prove strict
monotonicity almost everywhere.) 

 Let us fix the value of the total energy of the system,
$H=\frac 12$. In this manifold we consider the two dimensional section $\M$ of
the flow, corresponding to the bottom particle emerging from the collision
with the floor; the surface $\M$ is given by 
 $\{H=\frac 12,\, q_1 = 0,\, v_1 \geq 0\}$. The state of
the system in $M$ is completely described by the velocities of the 
particles $(v_1,\,v_2)$; and we use the velocities as coordinates in
$\M$. Hence, our phase space $\M$ is the  domain bounded by the half-ellipse
$$
m_1v_1^2 + m_2v_2^2 \leq 1,\ \ , v_1 \geq 0.
$$
Let us calculate the symplectic form in these coordinates. We have 
$$
\omega = dp_1\wedge dq_1 + dp_2\wedge dq_2.
$$
On the surface of section $\M$ 
$$ 
dq_1 \equiv 0  \  \ \text{ and } \ \ dq_2 = -\frac {m_1}{m_2}v_1dv_1 -
v_2dv_2.
$$
Hence, we get
$$
\omega = m_1 v_1 dv_1\wedge dv_2.
$$

The map $T: \M \to \M$ is defined by the first return of the flow to
$\M$. Our symplectic box $\M$ is split into two symplectic boxes by
$\Cal S^+$, which is the arc of the ellipse $\{ m_1v_1^2 +
m_2(v_2-2v_1)^2 = 1\}$ contained in $\M$. The symplectic box $\M^+_f$,
above $\Cal S^+$, contains all the initial states for which the bottom
particle returns to the
floor without colliding with the top particle. The map $T$ in $\M^+_f$
is linear
$$
T(v_1,\,v_2) =(v_1,\,v_2-2v_1).
$$
The symplectic box $\M^+_c$, below $\Cal S^+$, contains all the initial
states for which there is a collision of the two particles before the 
bottom particle returns to the floor. The map $T$ in $\M^+_c$ is 
nonlinear and is best described in a coordinate system
$(h,\,z)$ where
$$
\aligned
h =& \frac 12 m_1 v_1^2\\
z= & v_2 -v_1.
\endaligned
$$
The symplectic form $\omega =  d h\wedge d z.$
(This coordinate system is derived from the canonical system 
of coordinates in the full phase space furnished by the individual
energies and velocities of the particles. The exceptional role of
these coordinates is well documented in \cite {W7}, \cite {CW}.)

Note that both the energy of the bottom particle and the difference of
velocities change only in collisions. Now $T = F_2\circ F_1$, where
$$
F_1(h,\,z) = ( - h - a z^2 + b,\, -z), 
\ \ a = \frac{m_1m_2(m_1-m_2)}{(m_1+m_2)^2} \ \ \text { and }
\ \ b = \frac {m_1}{m_1+m_2} ,
$$
describes the collision of the two particles, and 
$$
F_2(h,\,z) = (h,\, z + c\sqrt h), \ \ c = \sqrt{\frac {8}{m_1}},
$$
describes the collision of the bottom particle with the floor.

To find the image symplectic boxes $\M^-_f$ and $\M^-_c$  we can use 
the reversibility of our system. Namely, if we put $S(v_1,\,v_2)=
(v_1,\,-v_2)$ then $T \circ S = S \circ T^{-1}$, and so
$\M^-_f = S \M^+_f,\, \M^-_c = S \M^+_c$.

Our bundle of unstable sectors is constant in the coordinates $(h,\,z)$
and equal to the positive (and negative) quadrant; the form 
$\Q = d h d z.$ It is
immediate that $\Cal S^+$ and $\Cal S^- = S \Cal S^+$ are properly aligned.

We can now check that $T$ is monotone in $\M^+_f$ and strictly
monotone in $\M^+_c$ (both $F_1$ and $F_2$ are monotone).
Indeed, in the $(h,\,z)$ coordinates we have 
$$
DF_1 = 
\left(\matrix -1& -2a z \\ 
0 & -1
\endmatrix\right) 
\ \text { and } 
DF_2 = 
\left(\matrix 1& 0 \\ 
\frac c{2\sqrt h}  & 1
\endmatrix\right). 
$$
Moreover the map $T$ in $\M^+_f$ is equal in the coordinates $(h,\,z)$ to
$F_2$.

 Since the collision of the two particles must eventually occur, we
obtain strict monotonicity of all nondegenerate orbits. Unboundedness
of all nondegenerate orbits follows from Proposition 6.9.
 So the Sinai-Chernov
Ansatz holds.

To check the noncontraction property, we observe that the standard Riemannian
metric in the coordinates $(h,\,z)$ does not decrease on vectors from the
sector, when we apply one of the above matrices.

Finally, we are guaranteed that the coefficient $\sigma (DT^n)$ can be
made arbitrarily large by increasing $n$, except for points which end up
on the singularity set $\Cal S^+$ in the future and the singularity set
$\Cal S^-$ in the past. There are only countably many such points in
view of the proper alignment of singularity sets, and the Main Theorem 
applies to all other points. It follows that $T$ is ergodic and 
consequently, by the results of Katok and Strelcyn, it is a Bernoulli
system. 

The case of variable acceleration ($U'' < 0$) can be treated in a
similar fashion. It is not possible to write down the formulas for 
the return map $T$ but its derivative in the coordinates
$$
\aligned 
\delta h &= \frac {p_1}{m_1} \delta p_1 \\
\delta z &= \frac {1}{m_2U'(q_2)} \delta p_2 - 
\frac {1}{m_1U'(q_1)} \delta p_1,
\endaligned
$$
was essentially calculated in \cite {W8}. It is again a product of 
triangular matrices. 

\bigskip
\subheading{Afterword}
\bigskip
This paper was greatly improved thanks to many insightful comments and
corrections by the anonymous referees of the paper.

While we were writing this paper, several authors pursued similar
goals. There are the papers by Chernov \cite{Ch1}, \cite{Ch2}, the new version
of his old preprint by Katok, in collaboration with Burns \cite
{K2},  by Markarian \cite{M}, by Vaienti \cite {Va}, 
and the papers by Sim\'anyi \cite{S1}, \cite{S2}.


\par\newpage 

\Refs

\widestnumber\key{XXXX}

\ref\key{AW} \by R.L. Adler, B. Weiss 
\paper Entropy is a complete metric invariant for automorphisms of 
the torus
\jour Proc. Natl. Acad. Sci. USA
\vol 57
\pages 1537 -- 1576
\yr 1967
\endref


\ref\key{AS} \by D.V. Anosov, Ya.G.Sinai 
\paper Certain smooth ergodic systems
\jour Russ. Math. Surv. 
\vol 22 
\pages 103 -- 167
\yr 1982
\endref

\ref\key{B} \by L. A. Bunimovich
\paper On the ergodic properties of nowhere dispersing billiards
\jour Comm.\-Math.Phys. \vol 65  \yr 1979 \pages 295 -- 312 \endref 


\ref \key{BG} \by K. Burns, M. Gerber
\paper Continuous invariant cone families and ergodicity of flows in
dimension three
\jour Erg.Th.Dyn.Syst. \vol 9 \yr 1989 \pages 19 -- 25  \endref 

\ref \key{CW} \by  J.Cheng,  M.P.Wojtkowski 
\paper Linear stability of a periodic orbit in the system of falling
balls
\pages 53 -- 71 
\jour The Geometry of Hamiltonian Systems, Proceedings of a Workshop 
Held June 5-16,1989 MSRI Publications, Springer Verlag 1991 (ed. Tudor Ratiu)
\endref

\ref\key{Ch 1} \by N.I. Chernov \paper The ergodicity of a Hamiltonian
system of two particles in an external field
\jour Physica D \vol 53 \yr 1991 \pages 233 -- 239  \endref 

\ref\key{Ch 2} \by N.I. Chernov\paper On local ergodicity in
hyperbolic systems with singularities \paperinfo preprint \yr 1991
\endref



\ref \key{CS} \by N.I.Chernov, Ya.G.Sinai
\paper Ergodic properties of some systems of $2$-dimensional discs
and $3$-dimen- sional spheres
\jour Russ.Math.Surv. \yr 1987 \vol 42 \pages 181 -- 207 \endref


\ref\key{D1} \by V. Donnay 
\paperinfo private communication \yr 1988
\endref
\ref\key{D2} \by V. Donnay 
\paper Using integrability to produce chaos: billiards with positive
entropy
\jour Comm.\- Math.Phys. \vol 141  \yr 1991 \pages 225 - 257 \endref 

 
\ref\key{Ha}\by B. Halpern \paper Strange Billiard Tables
  \jour TAMS \vol 232 \yr 1977 \pages 297 -- 305
\endref


\ref\key{H}\by E. Hopf \paper Statistik der Geodatischen Linien in
Mannigfaltigkeiten Negativer Krummung \jour Ber. Verh. S\"achs. akad.wiss.,
Leipzig \vol 91\pages 261 -- 304\yr 1939\endref


\ref  \key{K1} \by A. Katok 
\paper Invariant cone families and stochastic properties of smooth
dynamical systems \paperinfo preprint \yr 1988 \endref

\ref  \key{K2} \by A. Katok in collaboration with K. Burns
\paper Infinitesimal Lyapunov functions, invariant cone families 
and stochastic properties of smooth
dynamical systems \paperinfo preprint \yr 1992 \endref


\ref \key{KS} \by A. Katok, J.-M. Strelcyn with the collaboration of
F. Ledrappier and F. Przytycki 
\book Invariant manifolds, entropy and billiards; smooth maps with
singularities \bookinfo Lecture Notes in Math.  1222
\publ Springer-Verlag \yr 1986 \endref


\ref\key{KSS} \by A. Kr\'amli, N. Sim\'anyi, D. Sz\'asz
\paper A \lq\lq Transversal" Fundamental Theorem for Semi-Dis\-persing 
Billiards
\jour Communications in Mathematical Physics\vol 129 \yr 1990
\pages 535 -- 560\paperinfo (see also Erratum)\endref



\ref\key{LW} \by C. Liverani, M.P. Wojtkowski
\paper Generalization of the Hilbert metric to the space of positive definite
matrices
\paperinfo to appear in Pac. J. Math. 
\endref


\ref\key{M} \by R. Markarian
\paper The Fundamental Theorem of Sinai -- Chernov for dynamical
systems with singularities
\paperinfo preprint  \yr 1991
\endref


 
\ref\key{O}\by V. I. Oseledets \paper A Multiplicative Ergodic Theorem:
Characteristic Lyapunov Exponents of Dynamical Systems\jour Trans. Moscow
Math. Soc.\vol 19\yr 1968\pages 197 -- 231\endref
 
\ref\key{P} \by Ya. B. Pesin\paper Lyapunov Characteristic Exponents
 and Smooth Ergodic Theory\jour Russ. Math. Surveys \vol 32, {\rm 4}
\yr 1977\pages 55 -- 114\endref

\ref \key{S} \by Ya.G.Sinai \paper Dynamical systems with
elastic reflections \jour Russ.Math.Surveys \vol 25 \yr 1970 \pages
137 -- 189
\endref

\ref\key{Si1} \by  N. Sim\'anyi
\paper The K-property of N billiard balls I
\paperinfo preprint
\yr 1991
\endref

\ref\key{Si} \by  N. Sim\'anyi
\paper The K-property of N billiard balls II: Computation of neutral
linear spaces
\paperinfo preprint
\yr 1991
\endref

\ref\key{Va}\by S. Vaienti \paper Ergodic properties of the discontinuous
sawtooth map \jour Jour. Stat. Phys. \vol 67 \yr 1992 \pages 251 -- 269
\endref

\ref\key{Ve}\by E. Vesentini \paper Invariant metrics on convex cones
\jour Ann. Sc. Norm. Sup. Pisa ser. 4 \vol 3  \yr 1976  \pages 671 -- 696 
\endref

\ref \key{W1} \by M.P.Wojtkowski \paper Invariant families of cones
 and Lyapunov exponents \jour Erg.Th.Dyn.Syst. \vol 5
\yr 1985 \pages 145 -- 161 \endref 

\ref \key{W2} \by M.P.Wojtkowski \paper Measure theoretic
entropy of the system of hard spheres \jour Erg.Th.Dyn.\-Syst. \vol 8
\yr 1988 \pages 133 -- 153 \endref 

\ref \key{W3} \by
M.P. Wojtkowski \paper Systems of classical interacting particles
with nonvanishing Lyapunov exponents \pages 243 -- 262
\yr 1991 \jour Lecture Notes in Math. 1486,
Springer-Verlag \paperinfo Lyapunov Exponents, Proceedings,
Oberwolfach  1990, L. Arnold, H. Crauel, J.-P. Eckmann (Eds)
 \endref

\ref \key{W4} \by M.P.Wojtkowski 
\paper  Principles for the design of billiards with nonvanishing Lyapunov
exponents \jour Comm. Math. Phys. \vol 105 \pages 391 -- 414 \yr 1986
\endref

\ref \key{W5} \by M.P.Wojtkowski 
\paper  A model problem with the coexistence of stochastic and integrable
behavior
\jour Comm.Math.Phys. \vol 80 \pages  453 -- 464 \yr 1981
\endref

\ref \key{W6} \by
M.P.Wojtkowski 
\paper  On the ergodic properties of piecewise linear perturbations of the
twist map \jour Ergodic Theory \& Dynamical Systems \vol 2
\pages 525 -- 542 \yr 1982
\endref


\ref \key{W7} \by
M.P.Wojtkowski \paper A system of one dimensional balls with gravity
\jour Comm.Math.Phys. \vol 126 \yr 1990 \pages 507 -- 533  \endref 

\ref \key{W8} \by
M.P.Wojtkowski \paper The system of one dimensional balls in an external 
field. II
\jour Comm.\-Math.Phys. \vol 127 \yr 1990 \pages 425 -- 432 \endref 


\endRefs
\enddocument